\documentclass[11pt,a4paper,twoside]{article}
\usepackage{amssymb}
\usepackage{amsfonts}
\usepackage{geometry}
\usepackage{color}
\usepackage{mathtools}
\usepackage{amsmath,amsthm,amssymb,amsfonts,enumitem}
\usepackage{latexsym}
\usepackage{graphicx}
\usepackage{pgfpages}
\usepackage{mathrsfs}
\usepackage{cite}
\usepackage[hyphens]{url}
\usepackage[bookmarks=false,hidelinks]{hyperref}
\usepackage[affil-it]{authblk}
\usepackage{appendix}
\usepackage{amsmath}
\usepackage{comment}
\allowdisplaybreaks[1]

\newtheorem{definition}{Definition}[section]
\newtheorem{theorem}[definition]{Theorem}
\newtheorem{lemma}[definition]{Lemma}
\newtheorem{proposition}[definition]{Proposition}
\newtheorem{corollary}[definition]{Corollary}
\newtheorem{remark}[definition]{Remark}
\newtheorem{condition}{Assumption}

\newtheorem{example}[definition]{Example}

\numberwithin{equation}{section}
\begin{document}

\title{Variational exponential estimates and stationary large deviations for SPDEs}
\author[a]{Rui Bai}
	\author [a] {Chunrong Feng}
	\author[a,b]{Huaizhong Zhao}
	\affil[a]{Department of Mathematical Sciences, Durham University, DH1 3LE, UK}
	\affil[b] {Research Centre for Mathematics and Interdisciplinary Sciences, Shandong University, Qingdao 266237, China}

	\affil[ ]{rui.bai@durham.ac.uk, chunrong.feng@durham.ac.uk,  huaizhong.zhao@durham.ac.uk}
	\date{}

\maketitle

\begin{abstract}
We establish an abstract criterion for large deviations of
invariant measures when the deterministic dynamics has a
single globally attracting equilibrium. Deterministic
regularisation and attraction yield an endpoint estimate
for long controlled paths of bounded action, uniformly
over each fixed bounded initial set. This allows the
finite time Freidlin-Wentzell uniform LDP on bounded
initial sets to be combined directly with stationary
exponential norm estimates. The proof requires neither
a separate exponential tightness argument nor a
Dembo-Zeitouni uniform LDP, and imposes no decay condition
on differences between controlled trajectories.
Finite path approximation identifies the quasipotential
as the good rate function. We also establish a variational criterion based on the
Bou\'e-Dupuis representation that derives exponential
estimates from bounds for controlled outputs, without
requiring their convergence. Uniform controlled terminal
moments yield estimates uniform in time, which pass to
Krylov-Bogoliubov invariant measures.

Together, these results establish stationary LDPs in four
regimes: two-dimensional Navier-Stokes equations with
possibly unbounded, globally Lipschitz, uniformly positive
scalar amplitudes and vanishing noise correlation;
reaction-diffusion systems without the additional
restriction linking noise growth to reaction strength;
Allen-Cahn equations with opposite boundary values and
globally Lipschitz, uniformly nondegenerate multiplicative
space-time white noise on intervals of arbitrary length;
and Burgers equations with fixed covariance
$(-\Delta)^\alpha$, $0\leq\alpha<1/2$.
For Burgers, we also obtain Gaussian norm estimates.

\medskip
\noindent\textbf{MSC2020 subject classifications:}
Primary 60H15, 60F10; secondary 35R60, 37A50, 37L55.

\end{abstract}

\tableofcontents
\clearpage

\section{Introduction}
\label{sec:introduction}

A large deviation principle for invariant measures describes rare
fluctuations of a small noise stochastic system around its
deterministic equilibrium (Freidlin and
Wentzell~\cite{freidlin2012random}). In the single equilibrium setting
considered here, the rate function is a quasipotential: the least
control energy required to reach a prescribed state from the
equilibrium, with the travel time allowed to vary.

Early work on stationary small noise asymptotics for SPDEs
includes Freidlin~\cite{freidlin1988random} and
Sowers~\cite{sowers1992large}. Subsequent results include
reaction-diffusion systems with multiplicative noise
\cite{cerrai2005large}, on unbounded domains~\cite{wang2024large}, two-dimensional Navier-Stokes equations
\cite{brzezniak2017large,cerrai2022large},  and
stochastic Burgers equations~\cite{BAI2026111284}.
For SPDEs, implementing the classical stationary theory requires
additional analysis of spatial compactness, controllability
and uniform finite time probability estimates. These issues
become particularly delicate when the noise is spatially rough
or its coefficient is unbounded.

This paper makes three contributions. 

(1) We prove a general stationary LDP theorem for dynamics
with a single globally attracting deterministic equilibrium.
The key step is a deterministic endpoint estimate for long
controlled paths with bounded total control energy, uniform
over each fixed bounded set of initial values.
Together with the FWULDP on bounded initial sets and stationary
exponential estimates outside norm balls, this estimate yields
the stationary upper bound.
The argument requires neither a separate proof of stochastic
exponential tightness nor a Dembo-Zeitouni uniform LDP,
and it imposes no contraction assumption.

(2) We establish a general variational criterion that
derives exponential estimates from bounds for controlled
outputs. When these bounds are uniform in time, the criterion
gives exponential estimates over arbitrarily long time horizons;
convergence of the controlled outputs is not required.

(3) As applications, these two abstract results provide the framework and
exponential estimates needed to establish stationary LDP and significantly improve the applicability of LDP for invariant measures to SPDEs. In particular, we study the previously untreated regimes for stochastic two-dimensional
Navier-Stokes equations, stochastic reaction-diffusion systems, and
stochastic Allen-Cahn and stochastic Burgers equations. These four classes of SPDEs have been subject to intensive studies in the literature, especially regarding their invariant measures and large deviations over the last 40 years.

\subsection{Abstract framework and main results}

Let $\mathsf B$ be a separable Banach state space and
$u_\ast\in\mathsf B$ the unique globally attracting
equilibrium of the deterministic equation. The skeleton equation is the controlled
deterministic equation obtained by replacing the noise with a
square integrable forcing through the same coefficient. Its
action is the least control energy required to produce a given
trajectory.

The stochastic input to
Theorem~\ref{thm:stationary-transfer} consists of the finite time
Freidlin-Wentzell uniform LDP (FWULDP) on bounded initial sets
and the following estimate: for every $r>0$, there exist
$R_r,\epsilon_r>0$ such that
\begin{equation}\label{eq:intro-stationary-norm-tail}
 \nu_\epsilon\bigl(\{x\in\mathsf B:
                    \|x\|_{\mathsf B}>R_r\}\bigr)
 \leq \exp\left(-\frac r\epsilon\right),
 \qquad 0<\epsilon\leq\epsilon_r,
\end{equation}
where $\nu_\epsilon$ is an invariant measure. The balls in
this estimate need not be compact.

The deterministic hypotheses are stated in
Assumption~\ref{ass:stationary-input}.
They require fixed time continuity of the skeleton in its
initial value and control, and compactness of finite time
action sublevels for each fixed initial value.
Controlled trajectories must remain bounded independently
of duration when the initial values and total control energies
are bounded. Their states must also form a relatively compact
set at a fixed positive time.
Attraction to $u_\ast$ is required uniformly on compact
initial sets; positive time compactness extends it to bounded
initial sets, as shown in
Lemma~\ref{lem:bounded-initial-attraction}.

Proposition~\ref{prop:deterministic-compactness} proves
compactness of the infinite past action sublevels.
Its part (b) establishes the endpoint estimate
\eqref{eq:stationary-endpoint-approximation}:
endpoints of sufficiently long controlled paths with
bounded action approach the corresponding infinite past
quasipotential sublevel, uniformly over each fixed
bounded initial set.
Assumption~\ref{ass:finite-path-approximation} then identifies
the infinite past and finite time quasipotentials.
It requires that each admissible history can be replaced by
a finite path from $u_\ast$ with the same endpoint and
arbitrarily small excess action.
The duration may depend on the history and the required
accuracy. This identification gives the finite time
quasipotential as a good rate function.

The endpoint estimate explains why
\eqref{eq:intro-stationary-norm-tail} suffices.
Invariance splits the initial values into a bounded ball
and its complement, whose contribution is controlled by
the stationary norm estimate.
Inside the ball, an endpoint deviation implies a fixed
positive distance from the corresponding finite time action
sublevel, so the bounded initial set FWULDP applies directly.
The ball need not be compact.
The proof therefore requires neither a Dembo-Zeitouni
uniform LDP nor the stochastic spatial exponential estimates
used to establish exponential tightness in
\cite{cerrai2022large,wang2024large}.
Remark~\ref{rem:stationary-mechanism} explains this mechanism
and its relation to earlier stationary upper bound arguments.
Our second abstract result,
Theorem~\ref{thm:abstract-exponential-estimate}, provides a
general method for obtaining exponential estimates through
the Bou\'e-Dupuis representation
\cite{budhiraja2000variational,budhiraja2008large}.
For each bound on control energy, it requires bounded
deterministic controlled outputs and a uniform first moment
bound on their differences from the stochastic controlled
outputs. These differences need not converge as
$\epsilon\to0$; separate first moment bounds suffice by
the triangle inequality. The theorem applies to general
measurable outputs in a separable Banach space, including
solution paths, values at a fixed time and stochastic
convolutions. The admissible noise interval is independent
of the exponential level.

For invariant measures, we apply the criterion to a single
solution value. Write $u_{\epsilon,f}^{x_0}$ for the controlled
solution starting from a fixed $x_0$, including the
deterministic case $\epsilon=0$. The central estimate is
\begin{equation}\label{eq:intro-central-bound}
 \sup_{T>0}\sup_{0\leq\epsilon\leq\epsilon_*}
 \sup_{f\in\mathcal P_T^N}
 \mathbb E\|u_{\epsilon,f}^{x_0}(T)\|_{\mathsf B}
 <\infty,\qquad N>0,
\end{equation}
where $\epsilon_*>0$ is independent of $N$ and $T$.
Here $\mathcal P_T^N$ consists of predictable controls
taking values in the noise Hilbert space $U$ and satisfying
\[
 \frac12\int_0^T\|f(t)\|_U^2\,dt\leq N
 \quad\text{almost surely}.
\]
The deterministic case of
\eqref{eq:intro-central-bound} and the triangle inequality
verify the criterion uniformly in $T$. The resulting
exponential estimate passes to weak limits of
Krylov-Bogoliubov averages, after enlarging the radius,
and gives \eqref{eq:intro-stationary-norm-tail}.
Only a moment of one solution value at a time is needed,
rather than a moment of the supremum over an infinite
time interval.

\subsection{Implications to SPDEs}

\paragraph{Stochastic 2-D Navier-Stokes.}
On the torus, we consider
\[
 du_\epsilon+[\nu Au_\epsilon+B(u_\epsilon)]\,dt
 =\sqrt\epsilon\,g(u_\epsilon)
       (I+\delta(\epsilon)A^\beta)^{-1/2}\,dW,
\]
where $A$ is the positive Stokes operator, $B$ is the
convection term, and $H$ is the space of square integrable,
divergence-free, mean zero velocity fields.
The function $g:H\to[g_*,\infty)$, with $g_*>0$, is a
globally Lipschitz scalar amplitude.

The covariance regularisation is that of
Cerrai-Paskal~\cite{cerrai2022large}.
Their stationary LDP result is for additive noise with $\beta>2$ and
$\epsilon\delta(\epsilon)^{-2/\beta}\to0$, using an
exponential estimate for the stochastic $H^1$ energy.
Theorem~\ref{thm:NS-stationary-LDP} permits $\beta>1$
under
\[
 \delta(\epsilon)\to0,
 \qquad
 \epsilon\delta(\epsilon)^{-1/\beta}\to0,
\]
and allows unbounded amplitudes such as
$g(u)=\sqrt{1+|u|_H^2}$.
The controlled $H$ energy estimate involves
$\epsilon\operatorname{Tr}(I+\delta(\epsilon)A^\beta)^{-1}$,
which tends to zero under this scaling.
Deterministic time integrability in $H^1$ provides
positive time compactness and a connection with
arbitrarily small action. The finite time FWULDP is
verified for the changing covariance.
Every fixed viscosity $\nu>0$ is covered.
\paragraph{Stochastic reaction-diffusion systems.}
Cerrai-R\"ockner~\cite{cerrai2005large} proved a stationary
LDP for reaction-diffusion systems with a monotone reaction term of polynomial growth (order $m$) and unbounded
multiplicative noise. Their additional Hypothesis~6 links
noise growth, reaction strength, dimension, and covariance
regularity. In one dimension with space-time white noise,
it requires
\[
 m>\max\{1+6\rho,2\}
 \quad\text{when}\quad
 |g(\xi,z)|\leq C(1+|z|^\rho).
\]
A cubic reaction therefore permits only $\rho<1/3$ under
this condition, whereas linear growth requires $m>7$.

Theorem~\ref{thm:RD-stationary-LDP} retains their
Hypotheses~1--5, but removes Hypothesis~6, so $m>1$ is sufficient.
In particular, it permits globally Lipschitz, linearly
growing noise with the cubic dissipative reactions. The dimension-dependent covariance
and controllability assumptions remain unchanged.
The controlled moment estimate follows by restarting
damped convolutions and absorbing their contributions
into the solution bound. The stationary theorem then
applies to the resulting norm estimates without a
stochastic spatial exponential estimate.

\paragraph{Stochastic Allen-Cahn with opposite boundary values.}
We consider
\[
 \partial_tu_\epsilon
 =\Delta u_\epsilon+u_\epsilon-u_\epsilon^3
       +\sqrt\epsilon\,g(\xi,u_\epsilon)\,\partial_tW,
 \qquad \xi\in(-L,L),\qquad
 u_\epsilon(t,\pm L)=\pm1.
\]
The double-well reaction is not globally monotone, so the
preceding reaction-diffusion result does not apply directly.
Nevertheless, the boundary values select a unique
Ginzburg-Landau energy minimiser $m_L$, and every
deterministic trajectory converges to it. The energy
also supplies deterministic spatial compactness and the
connection with small action required by the framework.

For additive space-time white noise, a Gibbs description
provides a direct route to stationary large deviations
\cite{reznikoff2005invariant}. Such a description is
generally unavailable for multiplicative noise.
Our earlier work~\cite{bai2025allencahn} established the
stationary LDP under the additional sublinear growth
condition
\[
 |g(\xi,z)|\leq C(1+|z|^\rho),
 \qquad 0\leq\rho<1.
\]
Although this allows unbounded coefficients, it excludes
globally Lipschitz coefficients with linear growth. Theorem~\ref{thm:Allen-Cahn-invariant-LDP} removes this
additional growth restriction. We assume that $g$ is
continuous, uniformly nondegenerate, and globally Lipschitz
in its state variable, uniformly in space.
The proof combines the deterministic energy and
quasipotential arguments of our earlier work with the
uniform controlled moment estimates developed here.
The interval length $L>0$ remains arbitrary.

\paragraph{Stochastic 1-D Burgers with fixed rough covariance.}
Our published 1-D Burgers result~\cite{BAI2026111284} uses
trace class covariances to approach a rough limit at
a rate linked to the noise magnitude.
Theorem~\ref{thm:Burgers-stationary-LDP} treats directly the fixed
covariance $(-\Delta)^\alpha$, $0\leq\alpha<1/2$, on
$(0,1)$ with zero Dirichlet boundary conditions.
This includes space-time white noise and more singular
forcing. Its infinite trace prevents a direct $L^2$
It\^o energy calculation for the full solution.
The finite time LDP and deterministic quasipotential
results from \cite{BAI2026111284} remain applicable;
the new controlled terminal bound supplies the stationary
norm estimate needed here.
\\

We also obtain Gaussian norm tails by applying the
quantitative variational criterion to these controlled
bounds. Independent and contemporaneous work of
De Vecchi, Jan\'ak and Priola~\cite{de2026exponential}
proves related exponential integrability results through
the Bou\'e-Dupuis formula and the damped
Ornstein-Uhlenbeck decomposition of
Da Prato-Debussche~\cite{da2007m}.
Although our approach shares the same variational
foundation, it provides an abstract criterion for
Gaussian tail estimates based on ordinary bounds
for controlled processes. In particular, the method
is not restricted to additive noise and can also
be applied to equations with multiplicative noise,
as demonstrated by the Allen-Cahn equation with
linear growth noise considered below.
\paragraph{Comparison with approaches based on contraction.}
Decay of differences between controlled paths
plays a central role in
Liu, Tang and Zhang~\cite{liu2026stationary} and
Wang and Yang~\cite{wangyang2026large}.
In \cite{liu2026stationary}, the authors prove an LDP
for stationary solutions constructed by moving the initial
time to $-\infty$. Their argument relies on decay of
differences between stochastic solutions driven by the
same control. Their Hilbert-Schmidt diffusion setting
excludes the space-time white noise allowed in our
reaction-diffusion and Allen-Cahn applications and the
fixed rough covariance of our Burgers equation.
Navier-Stokes amplitudes such as
$g(u)=\sqrt{1+|u|_H^2}$ also violate their additional
noise condition [A4]. In \cite{wangyang2026large}, the authors use a coming down
from infinity estimate to establish the FWULDP over the
entire state space. Together with controlled contraction
and a stationary second moment estimate, this yields
the stationary LDP without a separate exponential
tightness argument. Their condition $p>4$, with reaction
growth of order $p-1$, excludes our cubic reaction-diffusion
and Allen-Cahn examples. Their scalar reaction-diffusion
setting also excludes the Burgers and Navier-Stokes
equations considered here. 

Our framework uses deterministic attraction without requiring
contraction between distinct controlled trajectories. Allen-Cahn on long intervals
illustrates the distinction: the Poincar\'e estimate does
not yield global contraction when $L$ is large, but the
energy structure still gives convergence to $m_L$.
This attraction suffices for our stationary theorem,
without an upper restriction on $L$.

\paragraph{Organisation of the paper}
Section~\ref{sec:abstract-exponential-estimate} proves
the variational exponential estimate.
Section~\ref{sec:stationary-framework} derives stationary
norm estimates, establishes the deterministic compactness
and endpoint results, and proves the general stationary
LDP theorem.
Sections~\ref{sec:NS-application}, \ref{sec:reaction-diffusion},
\ref{sec:AC-application} and \ref{sec:burgers-application} give the four applications.

\section{An abstract exponential estimate}
\label{sec:abstract-exponential-estimate}

In this section, we establish an abstract criterion for exponential
estimates of stochastic dynamics driven by Gaussian noise. The argument
uses the Bou\'e-Dupuis variational representation and requires bounds
for bounded energy controlled processes, without requiring their
convergence as the noise vanishes. Subsequent sections apply this
criterion to stationary norm estimates and large deviations of
invariant measures.

Let $\beta=\{\beta_k\}_{k\geq1}$ be a sequence of independent
one-dimensional Brownian motions on a filtered probability space
$(\Omega,\mathcal F,\{\mathcal F_t\},\mathbb P)$ satisfying the usual
conditions. Equip $\mathbb R^\infty$ with the topology of coordinate-wise
convergence. Then $\mathbb R^\infty$ is Polish and $\beta$ is a random
variable taking values in $C([0,T];\mathbb R^\infty)$ for each $T>0$.
Define the separable Hilbert space
\[
U:=\left\{f=\{f_k\}_{k\geq1}\in\mathbb R^\infty:
                 \sum_{k\geq1}|f_k|^2<\infty\right\},
\qquad
\langle f,g\rangle_U:=\sum_{k\geq1}f_kg_k.
\]
Let $\mathcal P_T$ be the collection of predictable $U$-valued
processes $f$ satisfying
\[
    \mathbb P\left(\int_0^T\|f(s)\|_U^2\,ds<\infty\right)=1.
\]
For $N>0$, set
\[
\begin{aligned}
\mathcal S_T^N
&:=\left\{f\in L^2(0,T;U):
               \frac12\int_0^T\|f(s)\|_U^2\,ds\leq N\right\},\\
\mathcal P_T^N
&:=\left\{f\in\mathcal P_T:
               \frac12\int_0^T\|f(s)\|_U^2\,ds\leq N
               \quad\mathbb P\text{-a.s.}\right\}.
\end{aligned}
\]

Let $\mathsf X_0$ be a parameter set, and let $\mathsf X$ be a separable
Banach space. For $\epsilon>0$, $T>0$, and $x\in\mathsf X_0$, let
\[
\mathcal G_{x,T}^\epsilon,\mathcal G_{x,T}:
C([0,T];\mathbb R^\infty)\longrightarrow\mathsf X
\]
be measurable maps. Their output may be any quantity taking
values in the chosen separable Banach space $\mathsf X$, provided
the controlled bounds below hold. Each application specifies
the output and its norm. Examples include solution paths, values
at a fixed time, and stochastic convolutions.
In Section~\ref{sec:stationary-framework}, the superscripts
$\mathrm{path}$ and $\mathrm{term}$ will distinguish solution
paths from their values at the terminal time.
Define
\begin{equation}\label{eq:abstract-action}
I_{x,T}(u):=
\inf\left\{\frac12\int_0^T\|f(s)\|_U^2\,ds:
 f\in L^2(0,T;U),\quad
 u=\mathcal G_{x,T}\left(\int_0^\cdot f(s)\,ds\right)\right\},
\end{equation}
with $\inf\varnothing=\infty$, and set
\begin{equation}\label{eq:abstract-sublevel-set}
    \Phi_{x,T}(N):=\{u\in\mathsf X:I_{x,T}(u)\leq N\}.
\end{equation}
The functional $I_{x,T}$ is used to define these sublevel sets; no
large deviation principle or compactness of the sublevel sets is
assumed in the following criterion.

\begin{condition}\label{ass:abstract-controlled-bounds}
There is a fixed $\epsilon_0>0$ such that, for every $N>0$, there
are constants $C_N,\bar C_N>0$ satisfying
\begin{equation}\label{eq:abstract-skeleton-bound}
\sup_{f\in\mathcal S_T^N}
\left\|\mathcal G_{x,T}
       \left(\int_0^\cdot f(s)\,ds\right)\right\|_{\mathsf X}
\leq C_N
\end{equation}
and, for every $0<\epsilon\leq\epsilon_0$,
\begin{equation}\label{eq:abstract-difference-bound}
\sup_{f\in\mathcal P_T^N}
\mathbb E\left\|
\mathcal G_{x,T}^\epsilon
\left(\sqrt\epsilon\,\beta+\int_0^\cdot f(s)\,ds\right)
-\mathcal G_{x,T}\left(\int_0^\cdot f(s)\,ds\right)
\right\|_{\mathsf X}
\leq\bar C_N.
\end{equation}
The constants may depend on $x$ and $T$. When uniformity in these
parameters is asserted, the bounds and $\epsilon_0$ are assumed
uniform in the corresponding parameters.
\end{condition}

\begin{theorem}[Abstract exponential estimate]
\label{thm:abstract-exponential-estimate}
Under Assumption~\ref{ass:abstract-controlled-bounds}, for every
$T>0$, $x\in\mathsf X_0$, and $r>0$, there are $R,M>0$ such that,
for every $0<\epsilon\leq\epsilon_0$,
\begin{align}
&\mathbb P\left(
 \|\mathcal G_{x,T}^\epsilon(\sqrt\epsilon\,\beta)\|_{\mathsf X}
 \geq R\right)\nonumber\\
&\quad\leq
\mathbb P\left(\operatorname{dist}_{\mathsf X}
  \left(\mathcal G_{x,T}^\epsilon(\sqrt\epsilon\,\beta),
        \Phi_{x,T}(2r^2)\right)\geq M\right)
\leq\exp\left(-\frac{r^2}{\epsilon}\right).
\label{eq:abstract-exponential-tail}
\end{align}
If the constants in Assumption~\ref{ass:abstract-controlled-bounds}
are uniform in $T$, then $R$ and $M$ can be chosen independently
of $T$. The same assertion holds for uniformity in $x$.
Moreover, if
\[
    C_N+\bar C_N\leq C_x+C'\sqrt N,
    \qquad N>0,
\]
then $R$ can be chosen so that $R\leq C(C_x+r)$, with $C$
independent of $r$.
\end{theorem}

\begin{proof}
By the definition of $I_{x,T}$, if $u\in\Phi_{x,T}(2r^2)$,
there is $f\in\mathcal S_T^{3r^2}$ such that
\[
    u=\mathcal G_{x,T}\left(\int_0^\cdot f(s)\,ds\right).
\]
Here a control of energy at most $3r^2$ may be chosen even when
the infimum defining $I_{x,T}$ is not attained. Hence
\[
    \sup_{u\in\Phi_{x,T}(2r^2)}\|u\|_{\mathsf X}\leq C_{3r^2}.
\]
For $R=C_{3r^2}+M$, the triangle inequality gives
\[
\{u\in\mathsf X:\|u\|_{\mathsf X}\geq R\}
\subseteq
\{u\in\mathsf X:
       \operatorname{dist}_{\mathsf X}(u,\Phi_{x,T}(2r^2))\geq M\}.
\]
This proves the first inequality in
\eqref{eq:abstract-exponential-tail}. It remains to prove the
second inequality.

Define the bounded Lipschitz function
\begin{equation}\label{eq:abstract-test-function}
h_{r,T,M}(\varphi):=
2r^2-2r^2\min\left\{1,
\frac{\operatorname{dist}_{\mathsf X}
       (\varphi,\Phi_{x,T}(2r^2))}{M}\right\},
\qquad\varphi\in\mathsf X.
\end{equation}
Its Lipschitz constant is $2r^2/M$. Since
$\mathcal G_{x,T}^\epsilon$ is measurable,
$h_{r,T,M}(\mathcal G_{x,T}^\epsilon(\sqrt\epsilon\,\beta))$
is a measurable real-valued random variable. If
\[
\operatorname{dist}_{\mathsf X}
\left(\mathcal G_{x,T}^\epsilon(\sqrt\epsilon\,\beta),
      \Phi_{x,T}(2r^2)\right)\geq M,
\]
then this random variable is zero. Consequently,
\begin{align}
&\mathbb P\left(\operatorname{dist}_{\mathsf X}
 \left(\mathcal G_{x,T}^\epsilon(\sqrt\epsilon\,\beta),
       \Phi_{x,T}(2r^2)\right)\geq M\right)\nonumber\\
&\qquad\leq
\mathbb E\exp\left[
-\frac{h_{r,T,M}
 (\mathcal G_{x,T}^\epsilon(\sqrt\epsilon\,\beta))}{\epsilon}
\right].
\label{eq:abstract-cutoff-bound}
\end{align}
By the variational representation in
\cite[Theorem~2]{budhiraja2008large},
\begin{align}
&\epsilon\log\mathbb E\exp\left[
-\frac{h_{r,T,M}
 (\mathcal G_{x,T}^\epsilon(\sqrt\epsilon\,\beta))}{\epsilon}
\right]\nonumber\\
&\quad=-\inf_{f\in\mathcal P_T}\mathbb E\left\{
 \frac12\int_0^T\|f(s)\|_U^2\,ds
 +h_{r,T,M}\left[\mathcal G_{x,T}^\epsilon
       \left(\sqrt\epsilon\,\beta+\int_0^\cdot f(s)\,ds\right)
              \right]\right\}.
\label{eq:abstract-variational-formula}
\end{align}
For every $\eta>0$, choose an $\eta$-minimising control
$f_n\in\mathcal P_T$. Thus
\begin{align}
&\mathbb E\left\{
 \frac12\int_0^T\|f_n(s)\|_U^2\,ds
 +h_{r,T,M}\left[\mathcal G_{x,T}^\epsilon
       \left(\sqrt\epsilon\,\beta+\int_0^\cdot f_n(s)\,ds\right)
              \right]\right\}\nonumber\\
&\quad\leq
\inf_{f\in\mathcal P_T}\mathbb E\left\{
 \frac12\int_0^T\|f(s)\|_U^2\,ds
 +h_{r,T,M}\left[\mathcal G_{x,T}^\epsilon
       \left(\sqrt\epsilon\,\beta+\int_0^\cdot f(s)\,ds\right)
              \right]\right\}+\eta.
\label{eq:abstract-near-minimiser}
\end{align}
We first show that
\[
    \mathbb E\left[\frac12\int_0^T\|f_n(s)\|_U^2\,ds\right]
    \leq2r^2+\eta.
\]
Indeed, $0\leq h_{r,T,M}\leq2r^2$. If the expected control
cost were greater than $2r^2+\eta$, then
\eqref{eq:abstract-near-minimiser} would imply that the infimum
on its right-hand side is greater than $2r^2$. Taking $f=0$
shows instead that this infimum is at most $2r^2$, a contradiction.

We next localise the control to a bounded energy class. Stop
$f_n$ when its accumulated energy reaches $N$. The stopped
control belongs to $\mathcal P_T^N$, its energy cost is no larger
than the original cost, and the shifted inputs agree unless
$\frac12\int_0^T\|f_n(s)\|_U^2\,ds>N$. Markov's inequality
bounds the probability of this exceptional event by
$(2r^2+\eta)/N$. Since the test function is bounded by $2r^2$,
this yields the usual localisation estimate
\cite{budhiraja2000variational}. Denoting the localised control
again by $f_n$, we obtain
\begin{align}
&\mathbb E\left\{
 \frac12\int_0^T\|f_n(s)\|_U^2\,ds
 +h_{r,T,M}\left[\mathcal G_{x,T}^\epsilon
       \left(\sqrt\epsilon\,\beta+\int_0^\cdot f_n(s)\,ds\right)
              \right]\right\}\nonumber\\
&\quad\leq
\inf_{f\in\mathcal P_T}\mathbb E\left\{
 \frac12\int_0^T\|f(s)\|_U^2\,ds
 +h_{r,T,M}\left[\mathcal G_{x,T}^\epsilon
       \left(\sqrt\epsilon\,\beta+\int_0^\cdot f(s)\,ds\right)
              \right]\right\}\nonumber\\
&\qquad+\eta+\frac{4r^2(2r^2+\eta)}{N}.
\label{eq:abstract-control-localisation}
\end{align}
By the definition of $I_{x,T}$,
\begin{align}
&\inf_{\varphi\in\mathsf X}
    \{I_{x,T}(\varphi)+h_{r,T,M}(\varphi)\}\nonumber\\
&\quad=
\inf_{f\in L^2(0,T;U)}
\left\{\frac12\int_0^T\|f(s)\|_U^2\,ds
+h_{r,T,M}\left[\mathcal G_{x,T}
     \left(\int_0^\cdot f(s)\,ds\right)\right]\right\}\nonumber\\
&\quad\leq
\mathbb E\left\{\frac12\int_0^T\|f_n(s)\|_U^2\,ds
+h_{r,T,M}\left[\mathcal G_{x,T}
     \left(\int_0^\cdot f_n(s)\,ds\right)\right]\right\}.
\label{eq:abstract-deterministic-comparison}
\end{align}
The last inequality holds pathwise before taking expectations,
although $f_n$ is random.

Combining \eqref{eq:abstract-variational-formula},
\eqref{eq:abstract-control-localisation}, and
\eqref{eq:abstract-deterministic-comparison}, with $N=16r^2$,
gives
\begin{align*}
&\epsilon\log\mathbb E\exp\left[
-\frac{h_{r,T,M}
 (\mathcal G_{x,T}^\epsilon(\sqrt\epsilon\,\beta))}{\epsilon}
\right]
+\inf_{\varphi\in\mathsf X}
       \{I_{x,T}(\varphi)+h_{r,T,M}(\varphi)\}\\
&\quad\leq
-\mathbb E\left\{\frac12\int_0^T\|f_n(s)\|_U^2\,ds
+h_{r,T,M}\left[\mathcal G_{x,T}^\epsilon
\left(\sqrt\epsilon\,\beta+\int_0^\cdot f_n(s)\,ds\right)\right]
\right\}\\
&\qquad+
\mathbb E\left\{\frac12\int_0^T\|f_n(s)\|_U^2\,ds
+h_{r,T,M}\left[\mathcal G_{x,T}
             \left(\int_0^\cdot f_n(s)\,ds\right)\right]\right\}
+\eta+\frac{4r^2(2r^2+\eta)}{N}\\
&\quad\leq
\mathbb E\left\{
 h_{r,T,M}\left[\mathcal G_{x,T}
             \left(\int_0^\cdot f_n(s)\,ds\right)\right]
-h_{r,T,M}\left[\mathcal G_{x,T}^\epsilon
\left(\sqrt\epsilon\,\beta+\int_0^\cdot f_n(s)\,ds\right)\right]
\right\}
+\frac54\eta+\frac{r^2}{2}\\
&\quad\leq
\frac{2r^2}{M}\mathbb E\left\|
 \mathcal G_{x,T}\left(\int_0^\cdot f_n(s)\,ds\right)
-\mathcal G_{x,T}^\epsilon
\left(\sqrt\epsilon\,\beta+\int_0^\cdot f_n(s)\,ds\right)
\right\|_{\mathsf X}
+\frac54\eta+\frac{r^2}{2}.
\end{align*}
The last inequality uses the Lipschitz constant of $h_{r,T,M}$.
Since $f_n\in\mathcal P_T^N$, Assumption~\ref{ass:abstract-controlled-bounds}
bounds the expectation in the last line by $\bar C_N$ for
$0<\epsilon\leq\epsilon_0$. Choosing $M=8\bar C_N$ and
$0<\eta\leq r^2/5$ gives
\begin{align}
&\epsilon\log\mathbb E\exp\left[
-\frac{h_{r,T,M}
 (\mathcal G_{x,T}^\epsilon(\sqrt\epsilon\,\beta))}{\epsilon}
\right]
+\inf_{\varphi\in\mathsf X}
       \{I_{x,T}(\varphi)+h_{r,T,M}(\varphi)\}\nonumber\\
&\quad\leq
\frac{2\bar C_N}{M}r^2+\frac54\eta+\frac{r^2}{2}
\leq\frac34r^2+\frac54\eta\leq r^2.
\label{eq:abstract-laplace-bound}
\end{align}
For every $\varphi\in\mathsf X$, either
$\varphi\in\Phi_{x,T}(2r^2)$, in which case
$h_{r,T,M}(\varphi)=2r^2$, or
$\varphi\notin\Phi_{x,T}(2r^2)$, in which case
$I_{x,T}(\varphi)>2r^2$. Therefore
\begin{equation}\label{eq:abstract-action-plus-cutoff}
\inf_{\varphi\in\mathsf X}
       \{I_{x,T}(\varphi)+h_{r,T,M}(\varphi)\}\geq2r^2.
\end{equation}
Finally, combining \eqref{eq:abstract-cutoff-bound},
\eqref{eq:abstract-laplace-bound}, and
\eqref{eq:abstract-action-plus-cutoff}, we obtain
\begin{align*}
&\epsilon\log\mathbb P\left(
\operatorname{dist}_{\mathsf X}
\left(\mathcal G_{x,T}^\epsilon(\sqrt\epsilon\,\beta),
      \Phi_{x,T}(2r^2)\right)\geq M\right)\\
&\quad\leq
\epsilon\log\mathbb E\exp\left[
-\frac{h_{r,T,M}
 (\mathcal G_{x,T}^\epsilon(\sqrt\epsilon\,\beta))}{\epsilon}
\right]
+\inf_{\varphi\in\mathsf X}
       \{I_{x,T}(\varphi)+h_{r,T,M}(\varphi)\}-2r^2\\
&\quad\leq-r^2.
\end{align*}
This proves \eqref{eq:abstract-exponential-tail}.
The choices $M=8\bar C_{16r^2}$ and $R=C_{3r^2}+M$
show that uniformity of the hypotheses in $T$ gives uniform
radii. Under the additional square-root growth assumption,
\[
R=C_{3r^2}+8\bar C_{16r^2}
\leq9C_x+(\sqrt3+32)C'r,
\]
which proves the quantitative assertion.
\end{proof}

\begin{remark}\label{rem:boundedness-versus-convergence}
The philosophy of Theorem~\ref{thm:abstract-exponential-estimate}
is related to the weak-convergence approach to large deviations
\cite{budhiraja2008large}. Both use the variational representation,
localisation of controls, and estimates for controlled processes.
A weak-convergence criterion for a full large deviation principle
requires convergence of the controlled outputs. Here only boundedness
is required; the difference in \eqref{eq:abstract-difference-bound}
need not converge to zero as $\epsilon\to0$.

This distinction is useful for long time problems, where establishing
controlled convergence uniformly over unbounded time horizons can be
difficult. When the bounds in the assumption are uniform in time,
the theorem yields uniform in time exponential estimates. For
example, when its output is the single value $u_\epsilon(T)$,
the target space is the state space and $T$ is a parameter.

Classical exponential estimates are also often obtained by applying
It\^o's formula to a suitable Lyapunov functional. The present criterion
instead derives the exponential estimate from ordinary controlled
bounds. Separate bounds on the stochastic and deterministic controlled
outputs verify \eqref{eq:abstract-difference-bound} by the triangle
inequality. The statement with a radius of order $C_x+r$ uses the
additional $\sqrt N$ growth assumption explicitly stated in the theorem.
\end{remark}

\section{From controlled bounds to large deviations for invariant measures}
\label{sec:stationary-framework}

Theorem~\ref{thm:abstract-exponential-estimate} permits an
arbitrary measurable output in a separable Banach space.
The applications below use solution paths or values at a fixed
time according to the estimate needed.
For the invariant measure norm bound, we choose the terminal value
$u_\epsilon(T)$, so the target is the state space $\mathsf B$ and
$T$ is a parameter. A bound on $\mathbb E\|u_\epsilon(T)\|$
with a constant independent of $T$ suffices for that application.
Uniform controlled bounds for these terminal values provide the
stationary probability estimate used in this section.

\subsection{Controlled terminal values and invariant measure tails}
\label{subsec:solution-output-maps}

Let $\mathsf B$ be a separable Banach space. We specialise the
measurable maps of Section~\ref{sec:abstract-exponential-estimate}
to solution paths, with target $\mathsf X=C([0,T];\mathsf B)$.
Write $\mathcal G_{x,T}^{\epsilon,\mathrm{path}}$ and
$\mathcal G_{x,T}^{\mathrm{path}}$ for these maps, with initial
value $x\in\mathsf B$. The superscript records their target;
the maps in the abstract theorem remain generic.
For a control $f\in L^2(0,T;U)$, write
$h_f(t)=\int_0^t f(s)\,ds$. Define the controlled stochastic path
and the deterministic controlled path by
\begin{equation}\label{eq:controlled-skeleton-definition}
\begin{aligned}
 u_{\epsilon,f}^x
 &:=\mathcal G_{x,T}^{\epsilon,\mathrm{path}}(\sqrt\epsilon\,\beta+h_f),
       &&\epsilon>0,\\
 z_f^x:=u_{0,f}^x
 &:=\mathcal G_{x,T}^{\mathrm{path}}(h_f),
       &&z_f^x(0)=x.
\end{aligned}
\end{equation}
Both outputs are paths in $C([0,T];\mathsf B)$ beginning at $x$.
We assume that
these maps give the unique controlled solutions, consistently
under restriction to shorter time intervals. The uncontrolled
stochastic solution is $u_\epsilon^x:=u_{\epsilon,0}^x$. For each
$\epsilon>0$, these uncontrolled solutions are assumed to form
a time-homogeneous Markov family.
The path $z_f^x$ is called the \emph{skeleton}: it is obtained by
removing the stochastic integral and inserting $f$ through the
same noise coefficient. For deterministic $f$ it is deterministic;
for a predictable random control it is evaluated for each realised
control path and may be random. We write $z_f$ when the initial
value is understood, and $z_0^x$ when the control is zero.

For the terminal value, define $\mathrm{ev}_T(v):=v(T)$ on
$C([0,T];\mathsf B)$ and set
\begin{equation}\label{eq:terminal-solution-maps}
\begin{aligned}
 \mathcal G_{x,T}^{\epsilon,\mathrm{term}}
 &:=\mathrm{ev}_T\circ\mathcal G_{x,T}^{\epsilon,\mathrm{path}},\\
 \mathcal G_{x,T}^{\mathrm{term}}
 &:=\mathrm{ev}_T\circ\mathcal G_{x,T}^{\mathrm{path}}.
\end{aligned}
\end{equation}
These maps take values in $\mathsf B$ and are measurable because
$\mathrm{ev}_T$ is continuous. In particular,
\[
 \mathcal G_{x,T}^{\epsilon,\mathrm{term}}
       (\sqrt\epsilon\,\beta+h_f)=u_{\epsilon,f}^x(T),
 \qquad
 \mathcal G_{x,T}^{\mathrm{term}}(h_f)=z_f^x(T).
\]
The invariant measure criterion below applies
Theorem~\ref{thm:abstract-exponential-estimate} to the terminal
maps with $\mathsf X=\mathsf B$. The path maps, with
$\mathsf X=C([0,T];\mathsf B)$, are used for finite time
trajectory estimates. These are two choices of output for
the abstract theorem; its general measurable maps also cover
other measurable functions of the driving noise.

The control space $U$ and the classes $\mathcal S_T^N,\mathcal P_T^N$
are those of Section~\ref{sec:abstract-exponential-estimate}; in the
applications $U$ is identified with the noise Hilbert space using an
orthonormal basis. The first line of
\eqref{eq:controlled-skeleton-definition} is the Brownian shift
used in the variational representation.
In the SPDE applications, the solution outputs
admit measurable versions as functions of the driving Brownian
paths. For bounded energy predictable controls, Girsanov's
theorem and pathwise uniqueness identify their Brownian shifts
with the controlled solutions. The same identification is used
for the controlled convergence that proves the finite time LDP.

\begin{proposition}
\label{prop:terminal-criterion}
Fix $x_0\in\mathsf B$ and $\epsilon_*>0$. Suppose that for every $N>0$,
\begin{equation}\label{eq:terminal-controlled-bound}
 \sup_{T>0}\sup_{0\leq\epsilon\leq\epsilon_*}
 \sup_{f\in\mathcal P_T^N}
 \mathbb E\|u_{\epsilon,f}^{x_0}(T)\|_{\mathsf B}
 \leq K_N<\infty.
\end{equation}
For $\epsilon=0$ the equation is the controlled deterministic equation;
it may nevertheless be random when its control is random.
For every $r>0$ there is $R_r>0$ such that
\begin{equation}\label{eq:terminal-criterion-tail}
 \sup_{T>0}\mathbb P\bigl(\|u_\epsilon^{x_0}(T)\|_{\mathsf B}
                  \geq R_r\bigr)\leq e^{-r/\epsilon},
 \qquad0<\epsilon\leq\epsilon_*.
\end{equation}
Every invariant measure $\nu_\epsilon$ obtained as a weak limit of
Krylov-Bogoliubov averages starting from $x_0$ satisfies, after
enlarging $R_r$ if necessary,
\begin{equation}\label{eq:terminal-criterion-invariant}
 \nu_\epsilon\bigl(\|x\|_{\mathsf B}\geq R_r\bigr)
 \leq e^{-r/\epsilon},\qquad0<\epsilon\leq\epsilon_*.
\end{equation}
\end{proposition}
\begin{proof}
Use the terminal maps \eqref{eq:terminal-solution-maps}
with target $\mathsf X=\mathsf B$.
For a deterministic control the $\epsilon=0$ instance of
\eqref{eq:terminal-controlled-bound} gives the first hypothesis of
Theorem~\ref{thm:abstract-exponential-estimate}, with $C_N=K_N$.
For random bounded energy controls, the triangle inequality gives
\[
 \sup_{T>0}\sup_{0<\epsilon\leq\epsilon_*}\sup_{f\in\mathcal P_T^N}
 \mathbb E\|u_{\epsilon,f}^{x_0}(T)-u_{0,f}^{x_0}(T)\|_{\mathsf B}
 \leq2K_N.
\]
Thus its second hypothesis holds with $\bar C_N=2K_N$, uniformly
in $T$. Apply that theorem with $\sqrt r$ in place of $r$ to obtain
\eqref{eq:terminal-criterion-tail}.

Put $\nu_{\epsilon,T}=T^{-1}\int_0^T
\mathcal L(u_\epsilon^{x_0}(t))\,dt$ and let
$\nu_{\epsilon,T_n}\Rightarrow\nu_\epsilon$.
The open set $O_r=\{x:\|x\|_{\mathsf B}>R_r\}$ satisfies
\[
 \nu_\epsilon(O_r)
 \leq\liminf_{n\to\infty}\nu_{\epsilon,T_n}(O_r)
 \leq e^{-r/\epsilon}.
\]
The closed tail $\{\|x\|_{\mathsf B}\geq2R_r\}$ is contained in
$O_r$. This proves \eqref{eq:terminal-criterion-invariant} with
the enlarged radius.
\end{proof}

The estimate in Proposition~\ref{prop:terminal-criterion} concerns
the complements of bounded balls. The stationary theorem below
requires precisely this estimate. Its proof uses compactness of
the deterministic controlled dynamics, rather than an exponential
estimate outside compact subsets of the stochastic state space.

\subsection{Basic properties of the deterministic controlled equation}

The next argument is common to the four applications. We first
obtain compactness on each fixed time interval from the continuity of
the skeleton and the compactness of the finite time rate function
sublevels. Spatial estimates then allow us to take limits of paths
defined on longer and longer intervals. This gives the endpoint
estimate used in the stationary upper bound.

We use the skeleton $z_f^x$ defined in
\eqref{eq:controlled-skeleton-definition}. In the semilinear
applications of this paper it is the unique mild solution of
\begin{equation}\label{eq:general-skeleton-mild}
\begin{aligned}
 z_f^x(t)
 &=S(t)x+\int_0^t S(t-s)F(z_f^x(s))\,ds\\
 &\quad+\int_0^t S(t-s)\Sigma(z_f^x(s))f(s)\,ds,
 \qquad 0\leq t\leq T.
\end{aligned}
\end{equation}
Here $S$ is the linear semigroup, $F$ is the nonlinear drift,
and $\Sigma$ is the operator multiplying the Brownian increment
in the stochastic equation. The coefficients are specified in
each application. The abstract statements use the solution-map
definition \eqref{eq:controlled-skeleton-definition} and the
properties below.

We assume that the equation is autonomous and that solutions can
be restarted and concatenated with their controls. Let $u_\ast$ be an
equilibrium of the uncontrolled equation. For $x\in\mathsf B$ and
$z\in C([0,T];\mathsf B)$, define the finite time action by
\begin{equation}\label{eq:finite-time-action}
 I_T^x(z):=\frac12\inf\left\{
     \int_0^T\|f(s)\|_U^2\,ds:
     f\in L^2(0,T;U),\ z=z_f^x\right\},
\end{equation}
with $\inf\varnothing=+\infty$. On an interval $[s,t]$, the
action $I_{s,t}$ is defined by translating time and using the
initial value $z(s)$. When applied to a path on a larger interval,
it denotes the action of its restriction to $[s,t]$.
Restriction and concatenation of controls give
\begin{equation}\label{eq:finite-action-additivity}
 I_{s,t}(z)=I_{s,u}(z)+I_{u,t}(z),\qquad s<u<t,
\end{equation}
when these actions are finite. For the reverse inequality, choose
controls on the two subintervals with arbitrarily small excess
cost and concatenate them; their costs need not attain the infima.

\begin{definition}
\label{def:quasipotentials}
For $z\in C((-\infty,0];\mathsf B)$, define
\begin{equation}\label{eq:infinite-past-action}
 I_{-\infty}(z):=\sup_{T>0}I_{-T,0}(z).
\end{equation}
For the skeleton equations considered below,
$I_{-\infty}(z)$, when finite, is equivalently the minimum
energy $\frac12\|f\|_{L^2(-\infty,0;U)}^2$ over controls
generating $z$ on the whole past, as follows by restriction
and weak compactness of bounded controls.
For $y\in\mathsf B$, set
\begin{equation}\label{eq:finite-time-quasipotential}
 V(y):=\inf_{\substack{T>0,\ z\in C([0,T];\mathsf B)\\
                         z(0)=u_\ast,\ z(T)=y}} I_T^{u_\ast}(z)
\end{equation}
and
\begin{equation}\label{eq:infinite-past-quasipotential}
 \widehat V(y):=
 \inf_{\substack{z\in C((-\infty,0];\mathsf B),\ z(0)=y\\
       \lim_{t\to-\infty}\|z(t)-u_\ast\|_{\mathsf B}=0}}
       I_{-\infty}(z).
\end{equation}
In both definitions, $\inf\varnothing=+\infty$. The function $V$
uses finite paths starting exactly at $u_\ast$; $\widehat V$ uses paths
that approach $u_\ast$ in the remote past and end at time zero.
Their sublevel sets are
\[
 \Phi_V(r):=\{y\in\mathsf B:V(y)\leq r\},\qquad
 \widehat\Phi(r):=\{y\in\mathsf B:\widehat V(y)\leq r\},
 \qquad r\geq0.
\]
Their equality will be established under
Assumption~\ref{ass:finite-path-approximation}; it is not part
of the definition.
\end{definition}

Write $\|z\|_T=\sup_{0\leq t\leq T}\|z(t)\|_{\mathsf B}$ and
\[
 K_T^x(r)=\{z\in C([0,T];\mathsf B):I_T^x(z)\leq r\}.
\]
For a compact initial set $\Lambda$, let
\[
 K_{\Lambda,T}(r)
 =\{z\in C([0,T];\mathsf B):
             z(0)\in\Lambda,\ I_T^{z(0)}(z)\leq r\}.
\]
Distance between paths is measured by $\|\cdot\|_T$.

\begin{condition}
\label{ass:stationary-input}
The controlled equation has the following properties.
\begin{enumerate}[label=(\roman*)]
\item
For every $T>0$ and $x\in\mathsf B$, the sets $K_T^x(r)$ are
compact for all $r\geq0$. The solution depends continuously on
its initial value and control, uniformly on bounded initial sets
and bounded energy control sets: for every $T,R,N>0$,
\begin{equation}\label{eq:stationary-skeleton-continuity}
 \lim_{h\downarrow0}
 \sup_{\substack{\|x\|_{\mathsf B},\|y\|_{\mathsf B}\leq R\\
                  f,g\in\mathcal S_T^N\\
                  \|x-y\|_{\mathsf B}+\|f-g\|_{L^2(0,T;U)}\leq h}}
          \|z_f^x-z_g^y\|_T=0.
\end{equation}
This is a fixed time statement; its constants may depend on $T$.

\item
For every $R,N>0$,
\begin{equation}\label{eq:det-uniform-bound}
 \sup_{T>0}\sup_{\|x\|_{\mathsf B}\leq R}
       \sup_{f\in\mathcal S_T^N}\|z_f^x\|_T<\infty.
\end{equation}
Moreover, for some fixed $t_0>0$, the set
\begin{equation}\label{eq:det-positive-time-compactness}
 \left\{z_f^x(t_0):\|x\|_{\mathsf B}\leq R,\
                     f\in\mathcal S_{t_0}^N\right\}
\end{equation}
is relatively compact in $\mathsf B$ for every $R,N>0$.

\item\label{item:stationary-attraction}
The uncontrolled solution approaches $u_\ast$ uniformly on each
nonempty compact set $\Lambda\subset\mathsf B$:
\begin{equation}\label{eq:stationary-uniform-attraction}
 \lim_{T\to\infty}
       \sup_{x\in\Lambda}\|z_0^x(T)-u_\ast\|_{\mathsf B}=0.
\end{equation}
\end{enumerate}
\end{condition}

Condition (ii) can be checked by a bound in a stronger spatial norm.
An integral energy estimate may first select an intermediate time
depending on the path; a forward energy estimate then propagates
the bound to the prescribed time $t_0$. Restarting over intervals
of length $t_0$ yields a common compact set at all later times,
as shown in Lemma~\ref{lem:det-common-compact} below.

\begin{condition}
\label{ass:finite-path-approximation}
For every path
$z\in C((-\infty,0];\mathsf B)$ with
$\lim_{t\to-\infty}z(t)=u_\ast$ and
$I_{-\infty}(z)<\infty$,
and every $\eta>0$, there are $S_\eta>0$ and
$\bar v_\eta\in C([-S_\eta,0];\mathsf B)$ such that
\begin{equation}\label{eq:finite-path-approximation}
 \bar v_\eta(-S_\eta)=u_\ast,\qquad
 \bar v_\eta(0)=z(0),\qquad
 I_{-S_\eta,0}(\bar v_\eta)\leq I_{-\infty}(z)+\eta.
\end{equation}
The endpoint agrees exactly with that of $z$. The duration
may depend on $z$ and $\eta$; no uniform duration is assumed.
This approximation is an additional input, not a consequence
of finite time action compactness. It identifies the two
quasipotentials.
\end{condition}

\begin{lemma}
\label{lem:compact-initial-action}
Under Assumption~\ref{ass:stationary-input}(i),
$K_{\Lambda,T}(r)$ is compact whenever $\Lambda$ is compact.
Moreover, if $x_n\to x$ and $z_n\to z$ in
$C([0,T];\mathsf B)$, then
\begin{equation}\label{eq:joint-action-lsc}
 I_T^x(z)\leq\liminf_{n\to\infty}I_T^{x_n}(z_n).
\end{equation}
\end{lemma}
\begin{proof}
Take $z_n\in K_{\Lambda,T}(r)$ and put $x_n=z_n(0)$.
After taking a subsequence, $x_n\to x\in\Lambda$. Choose controls
$f_n$ producing $z_n$ with
\[
 \frac12\|f_n\|_{L^2(0,T;U)}^2\leq r+\frac1n.
\]
Drive the equation from $x$ with the same controls, writing
$\widetilde z_n=z_{f_n}^x$. By
\eqref{eq:stationary-skeleton-continuity},
\[
 \|z_n-\widetilde z_n\|_T\longrightarrow0.
\]
The paths $\widetilde z_n$ belong to the compact set $K_T^x(r+1)$.
A subsequence therefore converges to a path $z$, and $z_n\to z$
along the same subsequence. Compact sublevels of $I_T^x$ imply
lower semicontinuity, so
\[
 I_T^x(z)\leq\liminf_n I_T^x(\widetilde z_n)
             \leq\liminf_n\frac12\|f_n\|_{L^2}^2\leq r.
\]
Thus $z\in K_{\Lambda,T}(r)$.

For \eqref{eq:joint-action-lsc}, there is nothing to prove if the
liminf on the right is infinite. Otherwise choose a subsequence
along which the actions converge to that finite liminf and
controls with energy at most $I_T^{x_n}(z_n)+1/n$. The same
comparison with $\widetilde z_n=z_{f_n}^x$ gives
$\widetilde z_n\to z$. Lower semicontinuity of $I_T^x$ proves
the assertion. No strongly convergent subsequence of the controls
is required.
\end{proof}

\begin{lemma}
\label{lem:det-common-compact}
Under Assumption~\ref{ass:stationary-input}(ii), let $t_0>0$
be the duration in \eqref{eq:det-positive-time-compactness}.
For every $R>0$ and $r\geq0$ there is a compact
$C_{R,r}\subset\mathsf B$ containing $z(t)$ whenever
\[
 \|z(0)\|_{\mathsf B}\leq R,\qquad
 I_T^{z(0)}(z)\leq r,\qquad t_0\leq t\leq T.
\]
The same compact set works for every such $t$ and $T$.
\end{lemma}
\begin{proof}
Choose a control producing $z$ with energy at most $r+1$.
Let $M$ be a bound in \eqref{eq:det-uniform-bound} for
initial radius $R$ and energy budget $r+1$. Restart on $[t-t_0,t]$, using elapsed time
$s\in[0,t_0]$. Autonomy of the equation gives
\[
 z(t)=z_h^{z(t-t_0)}(t_0),\qquad h(s)=f(t-t_0+s).
\]
Its starting value has norm at most $M$, and the shifted control
has energy at most $r+1$. Therefore $z(t)$ belongs to the closure
of the set in \eqref{eq:det-positive-time-compactness} with
$R=M$ and $N=r+1$. That closure is a suitable $C_{R,r}$,
independent of $t$ and $T$. The absolute restart time varies,
but every restarted path has the same duration $t_0$.
\end{proof}

\begin{lemma}
\label{lem:bounded-initial-attraction}
Under Assumption~\ref{ass:stationary-input}(ii),(iii), for every
$R>0$,
\begin{equation}\label{eq:bounded-initial-attraction}
 \lim_{T\to\infty}\sup_{\|x\|_{\mathsf B}\leq R}
       \|z_0^x(T)-u_\ast\|_{\mathsf B}=0.
\end{equation}
\end{lemma}
\begin{proof}
Let $t_0$ be the time in
\eqref{eq:det-positive-time-compactness}, and let $C_R$ be the
compact closure of
$\{z_0^x(t_0):\|x\|_{\mathsf B}\leq R\}$.
Autonomy and uniqueness give, for $T\geq t_0$,
\[
 \sup_{\|x\|_{\mathsf B}\leq R}\|z_0^x(T)-u_\ast\|_{\mathsf B}
 \leq\sup_{y\in C_R}\|z_0^y(T-t_0)-u_\ast\|_{\mathsf B}.
\]
The right side tends to zero by
Assumption~\ref{ass:stationary-input}(iii).
\end{proof}

Using the infinite past action in
\eqref{eq:infinite-past-action}, set
\[
 \mathcal K_{-\infty}(r):=
 \left\{z\in C((-\infty,0];\mathsf B):
 I_{-\infty}(z)\leq r,\ \lim_{t\to-\infty}z(t)=u_\ast\right\}.
\]
The path space is equipped with uniform convergence on every
interval $[-k,0]$. This topology is metrized by
\[
 d_-(z,w):=\sum_{k=1}^{\infty}2^{-k}
       \left(1\wedge\|z-w\|_{C([-k,0];\mathsf B)}\right).
\]

\begin{lemma}
\label{lem:ancient-convergence}
Suppose that the continuity condition
\eqref{eq:stationary-skeleton-continuity} and the attraction
condition \eqref{eq:stationary-uniform-attraction} hold.
Let $z\in C((-\infty,0];\mathsf B)$ satisfy
\[
 I_{-\infty}(z)<\infty,
 \qquad
 \overline{\{z(t):t\leq0\}}
 \text{ is compact in }\mathsf B.
\]
Then
\[
 \lim_{t\to-\infty}\|z(t)-u_\ast\|_{\mathsf B}=0.
\]
\end{lemma}

\begin{proof}
Put $C=\overline{\{z(t):t\leq0\}}$ and fix $\eta>0$.
By \eqref{eq:stationary-uniform-attraction}, choose
$\tau>0$ such that
\[
 \sup_{x\in C}\|z_0^x(\tau)-u_\ast\|_{\mathsf B}<\eta/2.
\]
On this fixed interval, continuity
\eqref{eq:stationary-skeleton-continuity} gives
$b\in(0,1)$ such that
\[
 x\in C,\qquad
 \frac12\|f\|_{L^2(0,\tau;U)}^2<b
 \quad\Longrightarrow\quad
 \|z_f^x(\tau)-z_0^x(\tau)\|_{\mathsf B}<\eta/2.
\]

By additivity and nonnegativity of the action,
\[
 0\leq I_{t-\tau,t}(z)
 =I_{t-\tau,0}(z)-I_{t,0}(z)
 \leq I_{-\infty}(z)-I_{t,0}(z).
\]
The right-hand side tends to zero as $t\to-\infty$,
by the definition of $I_{-\infty}$.
Thus $I_{t-\tau,t}(z)<b$ for all sufficiently negative $t$.
For each such $t$, the definition of the action provides
a control $f_t\in L^2(0,\tau;U)$ satisfying
\[
 z(t-\tau+s)=z_{f_t}^{z(t-\tau)}(s),
 \qquad 0\leq s\leq\tau,
 \qquad
 \frac12\|f_t\|_{L^2(0,\tau;U)}^2<b.
\]
Since $z(t-\tau)\in C$, the preceding continuity estimate
and the choice of $\tau$ give
\begin{align*}
 \|z(t)-u_\ast\|_{\mathsf B}
 &\leq
 \|z_{f_t}^{z(t-\tau)}(\tau)
       -z_0^{z(t-\tau)}(\tau)\|_{\mathsf B}\\
 &\quad+\|z_0^{z(t-\tau)}(\tau)-u_\ast\|_{\mathsf B}
 <\eta.
\end{align*}
Since $\eta>0$ is arbitrary, the conclusion follows.
\end{proof}

\begin{proposition}
\label{prop:deterministic-compactness}
Let $V$ and $\widehat V$ be the quasipotentials of
Definition~\ref{def:quasipotentials}.
Under Assumption~\ref{ass:stationary-input}, the following
statements hold.
\begin{enumerate}[label=(\alph*)]
\item
For every $r\geq0$, $\mathcal K_{-\infty}(r)$ is compact
for uniform convergence on compact time intervals.
The infimum defining $\widehat V(y)$ is attained whenever
it is finite, and $\widehat\Phi(r)$ is compact in $\mathsf B$.

\item\label{item:stationary-quasipotential}
For every $R,r,\delta>0$, there exists
$T_0=T_0(R,r,\delta)>0$ such that, for every $T\geq T_0$,
\begin{equation}\label{eq:stationary-endpoint-approximation}
 \|x\|_{\mathsf B}\leq R,\qquad I_T^x(z)\leq r
 \quad\Longrightarrow\quad
 \operatorname{dist}_{\mathsf B}
       (z(T),\widehat\Phi(r))<\delta.
\end{equation}
The initial radius $R$ is arbitrary; the waiting time
may depend on it.
\end{enumerate}
Moreover, $\widehat V\leq V$ and
$\Phi_V(r)\subset\widehat\Phi(r)$, so every sublevel
$\Phi_V(r)$ has compact closure.
If Assumption~\ref{ass:finite-path-approximation} also holds,
then $V=\widehat V$.
Consequently, $V$ is a good rate function, and
\eqref{eq:stationary-endpoint-approximation} holds with
$\Phi_V(r)$ in place of $\widehat\Phi(r)$.
\end{proposition}

\begin{proof}
Extending a finite path from $u_\ast$ by the constant
$u_\ast$ in the earlier past, with zero control there,
gives $\widehat V\leq V$.
We first prove (a) and (b) using only
Assumption~\ref{ass:stationary-input}.

\emph{(a) Compactness on the infinite past.}
Fix $r\geq0$ and $R_0>\|u_\ast\|_{\mathsf B}+1$.
For any $z\in\mathcal K_{-\infty}(r)$ and $t\leq0$,
convergence to $u_\ast$ in the remote past allows us to
choose $s<t-t_0$ such that
$\|z(s)\|_{\mathsf B}\leq R_0$.
The restriction to $[s,t]$ has action at most $r$.
Lemma~\ref{lem:det-common-compact}, applied after translating
this interval, therefore gives
\[
 z(t)\in C_r:=C_{R_0,r}.
\]
The compact set $C_r$ is independent of $z$ and $t$.
It contains the range of every path in
$\mathcal K_{-\infty}(r)$.

Let $z_n\in\mathcal K_{-\infty}(r)$.
For each integer $k\geq1$, the restrictions
$z_n|_{[-k,0]}$, translated to $[0,k]$, belong to the
compact set $K_{C_r,k}(r)$ of
Lemma~\ref{lem:compact-initial-action}.
Successive extraction and a diagonal subsequence yield
$z\in C((-\infty,0];\mathsf B)$ such that
\[
 z_n\longrightarrow z
 \quad\text{in }C([-k,0];\mathsf B)
 \quad\text{for every }k\geq1.
\]
The limit takes values in $C_r$, and
\eqref{eq:joint-action-lsc} gives
\[
 I_{-k,0}(z)
 \leq\liminf_{n\to\infty}I_{-k,0}(z_n)\leq r,
 \qquad k\geq1.
\]
Hence $I_{-\infty}(z)\leq r$.
Lemma~\ref{lem:ancient-convergence} implies
$z(t)\to u_\ast$ as $t\to-\infty$.
Thus $z\in\mathcal K_{-\infty}(r)$, proving compactness.

Suppose $\widehat V(y)<\infty$.
Choose admissible histories $z_n$ ending at $y$ with
\[
 I_{-\infty}(z_n)\leq\widehat V(y)+\frac1n.
\]
They belong to the compact set
$\mathcal K_{-\infty}(\widehat V(y)+1)$.
Along a convergent subsequence, the limiting history $z$
still ends at $y$ and satisfies
\[
 I_{-k,0}(z)
 \leq\liminf_{n\to\infty}I_{-k,0}(z_n)
 \leq\widehat V(y),
 \qquad k\geq1.
\]
Therefore $I_{-\infty}(z)\leq\widehat V(y)$.
The reverse inequality follows from the definition of
$\widehat V$, so the infimum is attained.

Consequently,
\[
 \widehat\Phi(r)
 =\{z(0):z\in\mathcal K_{-\infty}(r)\}
\]
is compact, since evaluation at time zero is continuous.
The inequality $\widehat V\leq V$ also gives
$\Phi_V(r)\subset\widehat\Phi(r)$, so $\Phi_V(r)$ has
compact closure.

\emph{(b) Endpoints of long paths from bounded initial sets.}
Fix $R,r,\delta>0$.
Lemma~\ref{lem:det-common-compact} provides a compact set
$C_{R,r}$ such that
\[
 \|x\|_{\mathsf B}\leq R,\qquad I_T^x(z)\leq r
 \quad\Longrightarrow\quad
 z(t)\in C_{R,r},\qquad t_0\leq t\leq T.
\]
The same compact set works for every duration and every
initial value in the ball of radius $R$.
Thus, after the fixed time $t_0$, all these controlled
paths take values in a common compact set.

Suppose \eqref{eq:stationary-endpoint-approximation} fails.
There are $T_n\to\infty$, initial values
$\|x_n\|_{\mathsf B}\leq R$, and paths $z_n$ satisfying
\[
 z_n(0)=x_n,\qquad I_{T_n}^{x_n}(z_n)\leq r,\qquad
 \operatorname{dist}_{\mathsf B}
       (z_n(T_n),\widehat\Phi(r))\geq\delta.
\]
Translate their terminal times to zero:
\[
 \bar z_n(t):=z_n(T_n+t),\qquad -T_n\leq t\leq0.
\]
For each fixed integer $k\geq1$, once $T_n\geq k+t_0$,
\[
 \bar z_n(t)\in C_{R,r}\quad(-k\leq t\leq0),
 \qquad I_{-k,0}(\bar z_n)\leq r.
\]
In particular, the initial values $\bar z_n(-k)$
of these restrictions lie in $C_{R,r}$.
After time translation, the restrictions therefore belong
to the compact set $K_{C_{R,r},k}(r)$.

A diagonal subsequence gives
$\bar z\in C((-\infty,0];\mathsf B)$ such that
\[
 \bar z_n\longrightarrow\bar z
 \quad\text{in }C([-k,0];\mathsf B)
 \quad\text{for every }k\geq1.
\]
The limit has range in $C_{R,r}$, and
\eqref{eq:joint-action-lsc} yields
\[
 I_{-k,0}(\bar z)\leq r,\qquad k\geq1.
\]
Thus $I_{-\infty}(\bar z)\leq r$.
By Lemma~\ref{lem:ancient-convergence},
$\bar z(t)\to u_\ast$ as $t\to-\infty$.
Hence $\bar z\in\mathcal K_{-\infty}(r)$ and
$\bar z(0)\in\widehat\Phi(r)$.
Convergence at time zero now gives
\[
 \delta
 \leq\operatorname{dist}_{\mathsf B}
          (z_n(T_n),\widehat\Phi(r))
 \leq\|\bar z_n(0)-\bar z(0)\|_{\mathsf B}
 \longrightarrow0,
\]
a contradiction.
This proves (b).

\emph{Identification under finite path approximation.}
Suppose Assumption~\ref{ass:finite-path-approximation}
also holds.
For every admissible finite-action history $z$ ending
at $y$ and every $\eta>0$, it provides a finite path
from $u_\ast$ to $y$ satisfying
\[
 V(y)\leq I_{-S_\eta,0}(\bar v_\eta)
       \leq I_{-\infty}(z)+\eta.
\]
Letting $\eta\downarrow0$ and taking the infimum over
such histories gives $V\leq\widehat V$.
Together with the reverse inequality proved above,
this yields $V=\widehat V$.
Their sublevels coincide and are compact by (a).
Thus $V$ is a good rate function, and (b) holds with
$\Phi_V(r)$ in place of $\widehat\Phi(r)$.
\end{proof}

\subsection{Identification of the finite time quasipotential}

The compactness and path estimates (a) and (b) of
Proposition~\ref{prop:deterministic-compactness} use only
Assumption~\ref{ass:stationary-input}. Its identification
$V=\widehat V$, and hence the closedness of the finite time
quasipotential sublevels, additionally use
Assumption~\ref{ass:finite-path-approximation}. The following
criterion verifies this additional approximation in an
application. It uses only the autonomous skeleton and the
restriction and concatenation properties of its action,
not compactness or the approximation property itself.

\begin{proposition}
\label{prop:local-connection-identification}
For the autonomous skeleton and additive actions defined above,
suppose that for every $\eta>0$ there is $\delta>0$ with the
following property: whenever a path $w$ on $[0,1]$ satisfies
\[
 \|w(0)-u_\ast\|_{\mathsf B}<\delta,\qquad
 I_1^{w(0)}(w)<\delta,
\]
there are $s\in[1/2,1]$, $S>0$, and a finite-action path $v$
such that
\[
 v(0)=u_\ast,\qquad v(S)=w(s),\qquad I_S^{u_\ast}(v)\leq\eta.
\]
Then \eqref{eq:finite-path-approximation} holds. In particular,
\begin{equation}\label{eq:finite-infinite-identification}
 V(y)=\widehat V(y),\qquad y\in\mathsf B.
\end{equation}
\end{proposition}
\begin{proof}
Extending a finite path by $u_\ast$ gives $\widehat V\leq V$.
For the opposite inequality, take an infinite past path $z$
of finite action ending at $y$ and converging to $u_\ast$.
Fix $\eta>0$ and the corresponding $\delta$.
For $T$ sufficiently large,
\[
 \|z(-T)-u_\ast\|_{\mathsf B}<\delta,\qquad
 I_{-T,-T+1}(z)<\delta.
\]
Applying the local condition to $w(s)=z(-T+s)$ gives a
path from $u_\ast$ to $z(-T+s)$ with action at most $\eta$.
Concatenate it with $z|_{[-T+s,0]}$, placing the connecting
path just before time $-T+s$. For its resulting finite duration
$S_\eta$, the concatenated path $\bar v_\eta$ satisfies
\[
 \bar v_\eta(-S_\eta)=u_\ast,\qquad \bar v_\eta(0)=y,
 \qquad
 I_{-S_\eta,0}(\bar v_\eta)
       \leq\eta+I_{-T+s,0}(z)
       \leq\eta+I_{-\infty}(z).
\]
This proves \eqref{eq:finite-path-approximation}, and also
$V(y)\leq\eta+I_{-\infty}(z)$.
Let $\eta\downarrow0$ and take the infimum over $z$.
This proves $V\leq\widehat V$.
\end{proof}

The local condition concerns reachability after spatial
regularisation; it does not follow from continuity alone.
An application can instead verify
\eqref{eq:finite-path-approximation} directly. A published
identity $V=\widehat V$ also supplies it: for any admissible
$z$, choose a finite path to $z(0)$ with action at most
$V(z(0))+\eta\leq I_{-\infty}(z)+\eta$.
Thus, the approximation is checked separately from the
compactness of infinite past action sublevels and the endpoint estimate, which follow from
Assumption~\ref{ass:stationary-input} alone.

\subsection{Finite time probability estimates and the stationary LDP}

We now add the finite time uniform LDP. Assume that for every
$T,R,r,\delta,\gamma>0$ there is $\epsilon_0>0$ such that,
for $0<\epsilon\leq\epsilon_0$ and $\|x\|_{\mathsf B}\leq R$,
\begin{align}
 \mathbb P(\|u_\epsilon^x-z\|_T<\delta)
 &\geq\exp\left[-\frac{I_T^x(z)+\gamma}{\epsilon}\right],
 \qquad z\in K_T^x(r),
 \label{eq:stationary-finite-lower}\\
 \mathbb P(d_T(u_\epsilon^x,K_T^x(r))\geq\delta)
 &\leq\exp\left[-\frac{r-\gamma}{\epsilon}\right].
 \label{eq:stationary-finite-upper}
\end{align}
The same $\epsilon_0$ works over these initial values and paths.
Here $d_T$ is distance in the path norm $\|\cdot\|_T$.
The uniformity over bounded initial sets is used directly in
the theorem below. A uniform LDP available only on compact initial
sets would not suffice for the same argument.

\begin{theorem}[Stationary LDP from exponential norm estimates]
\label{thm:stationary-transfer}
Assume Assumption~\ref{ass:stationary-input} and the finite time
estimates \eqref{eq:stationary-finite-lower}-\eqref{eq:stationary-finite-upper},
uniformly over bounded initial sets. Let $\{\nu_\epsilon\}$ be
invariant measures satisfying the following norm estimate:
for every $b>0$, there are $R_b,\epsilon_b>0$ such that
\begin{equation}\label{eq:stationary-norm-tail}
 \nu_\epsilon\bigl(\|x\|_{\mathsf B}>R_b\bigr)
       \leq e^{-b/\epsilon},\qquad0<\epsilon\leq\epsilon_b.
\end{equation}
Then, for every open $O\subset\mathsf B$ and closed
$F\subset\mathsf B$,
\begin{align}
 \liminf_{\epsilon\downarrow0}\epsilon\log\nu_\epsilon(O)
 &\geq-\inf_{y\in O}V(y),
 \label{eq:stationary-lower-without-approximation}\\
 \limsup_{\epsilon\downarrow0}\epsilon\log\nu_\epsilon(F)
 &\leq-\inf_{y\in F}\widehat V(y).
 \label{eq:stationary-upper-without-approximation}
\end{align}
If Assumption~\ref{ass:finite-path-approximation} also holds,
then $V=\widehat V$ and $\{\nu_\epsilon\}$ satisfies an LDP
in $\mathsf B$, with speed $1/\epsilon$ and good rate function $V$.
\end{theorem}

\begin{proof}
Proposition~\ref{prop:deterministic-compactness} gives compact
sublevels of $\widehat V$ and the endpoint estimate
\eqref{eq:stationary-endpoint-approximation} for arbitrary fixed
initial radius. Lemma~\ref{lem:bounded-initial-attraction}
gives uniform attraction on those same bounded initial sets.
We first prove the probability bounds without identifying $V$
and $\widehat V$.

\emph{Lower bound.}
Fix $y$ with $V(y)<\infty$, $\delta>0$, and $\gamma>0$.
Use \eqref{eq:stationary-norm-tail}, for instance with $b=1$,
to choose $R$ such that
\[
 \nu_\epsilon(B_R)\geq\frac12,
 \qquad B_R=\{x:\|x\|_{\mathsf B}\leq R\},
\]
for all sufficiently small $\epsilon$. Choose a control
$\bar f$ taking $u_\ast$ to $y$ in time $T_2$, with energy at most
$V(y)+\gamma/4$. The continuity estimate
\eqref{eq:stationary-skeleton-continuity} gives $\eta>0$ such
that
\[
 \|w-u_\ast\|_{\mathsf B}<\eta
 \quad\Longrightarrow\quad
 \|z_{\bar f}^{w}(T_2)-y\|_{\mathsf B}<\delta/2.
\]
By \eqref{eq:bounded-initial-attraction}, choose a common
$T_1$ such that
$\sup_{x\in B_R}\|z_0^x(T_1)-u_\ast\|_{\mathsf B}<\eta$.
Starting from each $x\in B_R$, use zero control until $T_1$
and then $\bar f$. The resulting paths $z^x$ have the same
duration $T=T_1+T_2$, have action at most $V(y)+\gamma/4$,
and end within $\delta/2$ of $y$. The uniform lower bound
\eqref{eq:stationary-finite-lower} gives
\[
 \inf_{x\in B_R}\mathbb P\bigl(
        \|u_\epsilon^x-z^x\|_T<\delta/2\bigr)
       \geq e^{-(V(y)+\gamma/2)/\epsilon}.
\]
By invariance,
\[
 \nu_\epsilon(B(y,\delta))
       \geq\frac12 e^{-(V(y)+\gamma/2)/\epsilon}
       \geq e^{-(V(y)+\gamma)/\epsilon}
\]
for sufficiently small $\epsilon$. Taking balls contained in
$O$ proves \eqref{eq:stationary-lower-without-approximation}.

\emph{Upper bound near sublevels.}
Fix $r,\delta,\gamma>0$. The assertion is immediate if
$\gamma\geq r$, so assume $0<\gamma<r$.
Choose $R$ in \eqref{eq:stationary-norm-tail} with level $r+1$.
Next choose a fixed $T$ large enough that
\eqref{eq:stationary-endpoint-approximation} holds for this
$R$, action level $r$, and tolerance $\delta/2$.
If a stochastic path starting from $x\in B_R$ is within
$\delta/2$ of $K_T^x(r)$ in the path norm, its endpoint is
within $\delta$ of $\widehat\Phi(r)$. Therefore
\eqref{eq:stationary-finite-upper}, with error $\gamma/2$,
implies
\[
 \sup_{x\in B_R}\mathbb P\bigl(
   \operatorname{dist}_{\mathsf B}
         (u_\epsilon^x(T),\widehat\Phi(r))\geq\delta\bigr)
       \leq e^{-(r-\gamma/2)/\epsilon}.
\]
Invariance and the split into $B_R$ and its complement give
\begin{align}
 \nu_\epsilon\bigl(
   \operatorname{dist}_{\mathsf B}(x,\widehat\Phi(r))\geq\delta\bigr)
 &\leq e^{-(r+1)/\epsilon}+e^{-(r-\gamma/2)/\epsilon}
       \nonumber\\
 &\leq e^{-(r-\gamma)/\epsilon}
 \label{eq:stationary-level-upper-bounded}
\end{align}
after decreasing $\epsilon_0$ so that
$\epsilon_0\log2\leq\gamma/2$.

If $F$ is closed and $0<s<\inf_F\widehat V$, the compact
set $\widehat\Phi(s)$ is disjoint from $F$, and hence their
distance is positive. Apply
\eqref{eq:stationary-level-upper-bounded} with $r=s$, let
$\gamma\downarrow0$, and then let $s\uparrow\inf_F\widehat V$.
If this infimum is infinite, let $s\to\infty$; if it is zero,
the bound follows from probabilities being at most one.
This proves \eqref{eq:stationary-upper-without-approximation}.

\emph{Identification and goodness.}
Under Assumption~\ref{ass:finite-path-approximation},
Proposition~\ref{prop:deterministic-compactness} gives
$V=\widehat V$ and compact sublevels of $V$. The two bounds
then prove the LDP with this good rate function.
\end{proof}

\begin{remark}[Why stationary norm estimates suffice]
\label{rem:stationary-mechanism}
A common stationary upper bound argument combines two
deterministic estimates.
The first controls endpoints of bounded-action paths
starting sufficiently close to the equilibrium.
The second gives an action lower bound for paths that
remain outside that neighbourhood at successive sampling
times; see
\cite[Lemmas~7.2 and~7.3]{brzezniak2017large}
and
\cite[Lemmas~4.2 and~4.3]{cerrai2022large}.
These estimates are used to decompose trajectories
according to whether they visit the equilibrium
neighbourhood.

Proposition~\ref{prop:deterministic-compactness}(b)
gives the endpoint estimate directly for every fixed
bounded initial set.
After translating terminal times to zero, positive time
compactness and finite time action compactness produce
a limiting history with finite action and compact range.
Lemma~\ref{lem:ancient-convergence} shows that this history
approaches $u_\ast$ in the remote past.
Its endpoint therefore belongs to the corresponding
infinite past quasipotential sublevel.
This argument establishes the bounded initial set
endpoint estimate without a separate avoidance estimate
or decomposition according to visits near the equilibrium.

In particular, fix $R,r,\delta>0$ and choose $T$
sufficiently large that
\eqref{eq:stationary-endpoint-approximation}
holds with tolerance $\delta/2$.
Uniformly for $\|x\|_{\mathsf B}\leq R$,
\[
 \left\{
 \operatorname{dist}_{\mathsf B}
       (u_\epsilon^x(T),\widehat\Phi(r))\geq\delta
 \right\}
 \subset
 \left\{
 d_T(u_\epsilon^x,K_T^x(r))\geq\delta/2
 \right\}.
\]
The event on the right is precisely the distance event
controlled by the FWULDP.
Invariance then allows us to split the initial values
into the ball of radius $R$ and its complement.
The bounded initial set FWULDP controls the contribution
from inside the ball, and the stationary exponential
norm estimate controls its complement.

The ball need not be compact.
Compactness enters through the deterministic controlled
histories and the quasipotential sublevels.
Thus the stationary proof requires neither a separate
exponential tightness argument nor a Dembo-Zeitouni
uniform LDP on compact initial sets.
Assumption~\ref{ass:finite-path-approximation}
retains its separate role in identifying
$V=\widehat V$, giving the finite time quasipotential
as the good stationary rate function.
\end{remark}

\begin{remark}
For Allen-Cahn, $u_\ast=m_L-\psi$ and $V=U_L$.
Section~\ref{subsec:AC-abstract-verification} verifies the
deterministic assumptions and the small-cost connection.
Quasipotential compactness and the end estimate
\eqref{eq:stationary-endpoint-approximation} then follow from the
abstract results, without a second model specific proof.
No continuity estimate uniform over all time horizons
is required.
\end{remark}

\section{Two-dimensional Navier-Stokes equations with multiplicative noise}
\label{sec:NS-application}

Brze\'zniak, Cerrai and Freidlin~
\cite{brzezniak2015quasipotential} analysed quasipotentials
for the two-dimensional Navier-Stokes equation with additive
noise, established their convergence as the noise correlation
vanishes, and applied these results to exit time asymptotics.
Brze\'zniak and Cerrai~\cite{brzezniak2017large} proved
the stationary LDP on the torus for fixed smooth additive
noise. Cerrai-Paskal~\cite{cerrai2022large} subsequently
established the stationary LDP in the joint limit of vanishing
noise magnitude and correlation scale.

Here we consider the covariance regularisation of
\cite{cerrai2022large} with a globally Lipschitz, uniformly
positive scalar amplitude depending on the velocity.
The amplitude may be unbounded, and the result holds for
every $\beta>1$ under the scaling specified below.
Our stochastic argument uses the controlled energy estimate
in $H$. Spatial compactness is established for the skeleton,
while its integral $H^1$ estimate allows us to select a state
that can be reached from the equilibrium with small action.
These deterministic estimates provide the compactness and
quasipotential identification required by the abstract
stationary theorem, without a uniform stochastic $H^1$
moment bound.
\subsection{Setting and statement of the result}

Let $\mathbb T^2$ be a fixed two-dimensional flat torus, and set
\[
 H=\left\{
 u\in L^2(\mathbb T^2;\mathbb R^2):
 \operatorname{div}u=0,\quad
 \int_{\mathbb T^2}u\,d\xi=0
 \right\}.
\]
We write $|\cdot|$ for the norm in $H$, $P$ for the Leray
projection, and $A=-P\Delta$ for the positive Stokes operator.
For $s\geq0$, let $H^s=D(A^{s/2})$, endowed with the norm
$|A^{s/2}\cdot|$. In particular,
\[
 \mathcal V=H^1,\qquad
 \|u\|_{\mathcal V}=|A^{1/2}u|,
\]
and $\mathcal V$ is compactly embedded in $H$.
The operator $A$ has an orthonormal eigenbasis
$\{e_j\}_{j\geq1}$, with
\[
 Ae_j=\lambda_j e_j,\qquad
 0<\lambda_1\leq\lambda_2\leq\cdots,\qquad
 \lambda_j\asymp j.
\]
We shall use Poincar\'e's inequality
\[
 \lambda_1|u|^2\leq\|u\|_{\mathcal V}^2.
\]

For sufficiently regular divergence-free vector fields, put
\[
 B(u,v)=P((u\cdot\nabla)v),\qquad
 B(u)=B(u,u),\qquad
 b(u,v,w)=\langle B(u,v),w\rangle.
\]
The periodic boundary conditions give the cancellations
\begin{equation}\label{eq:NS-cancellations}
 b(u,v,v)=0,\qquad
 \langle B(u),Au\rangle=0,
\end{equation}
whenever the expressions are defined.
Fix $\nu>0$ and write $S(t)=e^{-\nu tA}$.

\begin{condition}\label{ass:NS}
The function $g:H\to\mathbb R$ is globally Lipschitz and satisfies
\begin{equation}\label{eq:NS-amplitude}
 g(u)\geq g_*>0,\qquad u\in H.
\end{equation}
Let $\beta>1$, and let $\delta(\epsilon)>0$ satisfy
\begin{equation}\label{eq:NS-correlation-scale}
 \delta(\epsilon)\longrightarrow0,\qquad
 \epsilon\delta(\epsilon)^{-1/\beta}\longrightarrow0
 \quad\text{as }\epsilon\downarrow0.
\end{equation}
Define
\begin{equation}\label{eq:NS-covariance}
 C_\epsilon=(I+\delta(\epsilon)A^\beta)^{-1},
 \qquad R_\epsilon=C_\epsilon^{1/2}.
\end{equation}
\end{condition}

Here $g(u)$ is a scalar multiplying every noise coordinate,
so the diffusion operator is $g(u)R_\epsilon$.
Its dependence on $u$ is through the velocity as an element
of $H$, rather than through pointwise multiplication by a
function of $u(\xi)$.
If $L_g$ is a global Lipschitz constant of $g$, then
\[
 |g(u)|\leq |g(0)|+L_g|u|,
 \qquad
 |g(u)|^2\leq C_g(1+|u|^2)
\]
for some $C_g>0$. No additional growth condition is imposed.

Let $W$ be a cylindrical Wiener process on $H$. We study
\begin{equation}\label{eq:NS-stochastic}
 du_\epsilon^x+
 \bigl[\nu Au_\epsilon^x+B(u_\epsilon^x)\bigr]\,dt
 =
 \sqrt\epsilon\,g(u_\epsilon^x)R_\epsilon\,dW,
 \qquad u_\epsilon^x(0)=x\in H.
\end{equation}

We first record the spectral estimates used below.
Since $\beta>1$ and $\lambda_j\asymp j$, $C_\epsilon$
is trace class for every fixed $\epsilon>0$. Moreover,
for $0<\delta\leq1$,
\begin{equation}\label{eq:NS-spectral-bounds}
 \operatorname{Tr}(I+\delta A^\beta)^{-1}
 =
 \sum_{j\geq1}\frac{1}{1+\delta\lambda_j^\beta}
 \leq C_\beta\delta^{-1/\beta},
 \qquad
 \|R_\epsilon\|_{\mathcal L(H)}\leq1.
\end{equation}
Consequently,
\[
 \rho_\epsilon:=\epsilon\operatorname{Tr}C_\epsilon
 \leq C_\beta\epsilon\delta(\epsilon)^{-1/\beta}
 \longrightarrow0.
\]

For the later comparison with the skeleton, we also estimate
$R_\epsilon-I$ as an operator from $H$ to $\mathcal V'$.
Put $\delta=\delta(\epsilon)$. Since
\[
 R_\epsilon e_j
 =(1+\delta\lambda_j^\beta)^{-1/2}e_j,
\]
the spectral representation gives
\begin{equation}\label{eq:NS-rho-theta}
 \theta_\epsilon^2
 :=
 \|A^{-1/2}(R_\epsilon-I)\|_{\mathcal L(H)}^2
 =
 \sup_{j\geq1}
 \frac{\left[1-(1+\delta\lambda_j^\beta)^{-1/2}\right]^2}
      {\lambda_j}.
\end{equation}
For every $j$,
\[
 0\leq1-(1+\delta\lambda_j^\beta)^{-1/2}
 =
 \frac{\delta\lambda_j^\beta}
 {\sqrt{1+\delta\lambda_j^\beta}
  \bigl(1+\sqrt{1+\delta\lambda_j^\beta}\bigr)}
 \leq\frac12\delta\lambda_j^\beta.
\]
If $\lambda_j\leq\delta^{-1/\beta}$, then
$2\beta-1>0$ gives
\begin{align*}
 \frac{\left[1-(1+\delta\lambda_j^\beta)^{-1/2}\right]^2}
      {\lambda_j}
 &\leq\frac{\delta^2}{4}\lambda_j^{2\beta-1}\\
 &\leq\frac{\delta^2}{4}
          \left(\delta^{-1/\beta}\right)^{2\beta-1}
 =\frac14\delta^{1/\beta}.
\end{align*}
If $\lambda_j>\delta^{-1/\beta}$, we instead use
$0\leq1-(1+\delta\lambda_j^\beta)^{-1/2}\leq1$ to obtain
\[
 \frac{\left[1-(1+\delta\lambda_j^\beta)^{-1/2}\right]^2}
      {\lambda_j}
 \leq\lambda_j^{-1}
 \leq\delta^{1/\beta}.
\]
Taking the supremum over $j$ therefore proves
\begin{equation}\label{eq:NS-theta-bound}
 \theta_\epsilon^2
 \leq\delta(\epsilon)^{1/\beta}
 \longrightarrow0.
\end{equation}
In particular, for $f\in H$ and $w\in\mathcal V$,
\[
 |\langle(R_\epsilon-I)f,w\rangle|
 \leq\theta_\epsilon|f|\,\|w\|_{\mathcal V}.
\]

For each fixed $\epsilon>0$, the diffusion coefficient
satisfies
\begin{align*}
 \|\sqrt\epsilon\,[g(u)-g(v)]R_\epsilon\|_{\mathrm{HS}}
 &\leq L_g\sqrt{\rho_\epsilon}\,|u-v|,\\
 \|\sqrt\epsilon\,g(u)R_\epsilon\|_{\mathrm{HS}}
 &\leq
 \sqrt{\rho_\epsilon}\bigl(|g(0)|+L_g|u|\bigr).
\end{align*}
It is therefore globally Lipschitz and has linear growth
as a map from $H$ into the Hilbert-Schmidt operators on $H$.
Theorem~1.1 of \cite{liu2010spde}, with the verification for
the two-dimensional Navier-Stokes equation in its
Example~3.3, gives a unique global solution on the given
stochastic basis, satisfying
\[
 u_\epsilon^x\in
 C([0,T];H)\cap L^2(0,T;\mathcal V)
 \qquad\text{almost surely},\qquad T>0.
\]
The same verification applies with periodic boundary
conditions. The solution defines a Markov semigroup,
which is Feller by continuous dependence in probability
on the initial value; see also the proof of
\cite[Proposition~5.2]{ferrario2024uniqueness}.

For all sufficiently small $\epsilon>0$, there is also
a unique invariant measure. We verify the noise conditions
used in \cite[Proposition~5.2 and Theorem~5.6]
{ferrario2024uniqueness}.
The linear growth coefficient in the Hilbert-Schmidt
bound above is $L_g\sqrt{\rho_\epsilon}$.
Their sufficient condition for uniqueness therefore becomes
\[
 \nu>\frac{3L_g^2\rho_\epsilon}{2\lambda_1}.
\]
Since $\rho_\epsilon\to0$, we can choose $\epsilon_*>0$
such that this holds whenever $0<\epsilon\leq\epsilon_*$.

To check their nondegeneracy condition, let $P_m$ be the
orthogonal projection onto the first $m$ Stokes modes.
For each fixed $\epsilon>0$ and $m\geq1$,
\[
 \bigl(\sqrt\epsilon\,g(u)R_\epsilon\bigr)
 \frac{R_\epsilon^{-1}P_m}{\sqrt\epsilon\,g(u)}
 =P_m,
\]
and
\[
 \sup_{u\in H}
 \left\|
 \frac{R_\epsilon^{-1}P_m}{\sqrt\epsilon\,g(u)}
 \right\|_{\mathcal L(H)}
 \leq
 \frac{\sqrt{1+\delta(\epsilon)\lambda_m^\beta}}
      {\sqrt\epsilon\,g_*}
 <\infty.
\]
Here $R_\epsilon^{-1}P_m$ is bounded because $P_m$ has
finite dimensional range.
Thus the noise has the required bounded right inverse
on every finite collection of Stokes modes.
Although the cited results are formulated with Dirichlet
boundary conditions, their existence and uniqueness
arguments extend to the present periodic setting using
the same energy cancellation, Poincar\'e inequality and
trilinear estimates.

Consequently, for $0<\epsilon\leq\epsilon_*$, the
Krylov-Bogoliubov averages
\[
 \mu_{\epsilon,T}
 :=
 \frac1T\int_0^T\mathcal L(u_\epsilon^0(t))\,dt
\]
are tight as $T\to\infty$, and every weak limit equals the
unique invariant measure, denoted by $\mu_\epsilon$.
The value of $\epsilon_*$ may be decreased in the estimates
below.

The limiting controlled equation is
\begin{equation}\label{eq:NS-skeleton}
 \partial_tz_f^x+\nu Az_f^x+B(z_f^x)
 =g(z_f^x)f,\qquad z_f^x(0)=x,
 \qquad f\in L^2(0,T;H).
\end{equation}
We denote its unique solution in
$C([0,T];H)\cap L^2(0,T;\mathcal V)$ by $z_f^x$.
The limiting noise amplitude defines this controlled
deterministic equation; the argument does not require an
$H$-valued stochastic equation driven by unregularised
space-time white noise.

For $x\in H$ and $z\in C([0,T];H)$, define the action
\begin{equation}\label{eq:NS-action}
 I_T^x(z)
 :=
 \frac12\inf\left\{
 \int_0^T|f(t)|^2\,dt:
 f\in L^2(0,T;H),\quad z=z_f^x
 \right\},
 \qquad\inf\varnothing=+\infty.
\end{equation}
In particular, $I_T^x(z)=+\infty$ if $z(0)\ne x$.
The quasipotential based at the deterministic equilibrium
$0$ is
\begin{equation}\label{eq:NS-quasipotential}
\begin{aligned}
 V_{\rm NS}(y)
 &:=
 \inf_{\substack{T>0,\ z\in C([0,T];H)\\
                  z(0)=0,\ z(T)=y}} I_T^0(z)\\
 &=
 \frac12\inf\left\{
 \int_0^T|f(t)|^2\,dt:
 T>0,\ f\in L^2(0,T;H),\quad z_f^0(T)=y
 \right\},
 \qquad y\in H.
\end{aligned}
\end{equation}
Thus $V_{\rm NS}:H\to[0,\infty]$ and $V_{\rm NS}(0)=0$.

\begin{theorem}\label{thm:NS-stationary-LDP}
Under Assumption~\ref{ass:NS}, there exists $\epsilon_*>0$
such that \eqref{eq:NS-stochastic} has a unique invariant
measure $\mu_\epsilon$ for $0<\epsilon\leq\epsilon_*$.
The family $\{\mu_\epsilon\}_{0<\epsilon\leq\epsilon_*}$
satisfies an LDP in $H$ as $\epsilon\downarrow0$, with
speed $\epsilon^{-1}$ and good rate function $V_{\rm NS}$
defined in \eqref{eq:NS-quasipotential}.
More precisely, for every open set $O\subset H$ and
closed set $F\subset H$,
\begin{align*}
 \liminf_{\epsilon\downarrow0}
 \epsilon\log\mu_\epsilon(O)
 &\geq-\inf_{y\in O}V_{\rm NS}(y),\\
 \limsup_{\epsilon\downarrow0}
 \epsilon\log\mu_\epsilon(F)
 &\leq-\inf_{y\in F}V_{\rm NS}(y).
\end{align*}
The viscosity $\nu>0$ is arbitrary; the threshold
$\epsilon_*$ may depend on $\nu$.
\end{theorem}

We prove the LDP by verifying the hypotheses of
Theorem~\ref{thm:stationary-transfer}.
Uniqueness of the invariant measure is recorded above
but is not needed for that argument.
Throughout the remainder of this section,
$\mathcal P_T^N$ and $\mathcal S_T^N$ are defined with
noise Hilbert space $H$.
Constants may depend on $\nu$, the torus and $g$;
they are independent of $\epsilon$ unless stated otherwise.

\subsection{Uniform controlled moments and exponential estimates}

Fix $T>0$. For $f\in\mathcal P_T^N$, let
$u_{\epsilon,f}^x$ denote the solution of
\begin{equation}\label{eq:NS-controlled}
\begin{aligned}
 du_{\epsilon,f}^x
 &+\bigl[\nu Au_{\epsilon,f}^x
              +B(u_{\epsilon,f}^x)\bigr]\,dt\\
 &=g(u_{\epsilon,f}^x)R_\epsilon f\,dt
   +\sqrt\epsilon\,g(u_{\epsilon,f}^x)R_\epsilon\,dW,
 \qquad u_{\epsilon,f}^x(0)=x.
\end{aligned}
\end{equation}
For $\epsilon=0$, we set $u_{0,f}^x=z_f^x$, where
$z_f^x$ solves \eqref{eq:NS-skeleton}.
When $f$ is random, this deterministic equation is solved
for each realisation of the control.
We use the conventions $R_0=I$ and $\rho_0=0$.
\begin{lemma}\label{lem:NS-controlled-moments}
There exists $\epsilon_*>0$, independent of $N$ and $T$, such that
\begin{equation}\label{eq:NS-terminal-moment}
 \sup_{T>0}\sup_{0\leq\epsilon\leq\epsilon_*}
 \sup_{f\in\mathcal P_T^N}
       \mathbb E|u_{\epsilon,f}^x(T)|^2
       \leq C e^{CN}(1+|x|^2+N),\qquad N>0.
\end{equation}
For deterministic controls, there is also a bound independent of
the duration:
\begin{equation}\label{eq:NS-skeleton-global-bound}
 \sup_{T>0}\sup_{|x|\leq R}\sup_{f\in\mathcal S_T^N}
 \left[\sup_{t\leq T}|z_f^x(t)|^2
       +\int_0^T\|z_f^x(t)\|_{\mathcal V}^2\,dt\right]
       \leq C_{R,N}.
\end{equation}
\end{lemma}

\begin{proof}
Fix $T>0$ and $f\in\mathcal P_T^N$, and write
$u=u_{\epsilon,f}^x$. We use the conventions
$\rho_0=0$ and $R_0=I$.

It\^o's formula and $\langle B(u),u\rangle=0$ give
\[
 d|u|^2+2\nu\|u\|_{\mathcal V}^2\,dt
 =
 2g(u)\langle R_\epsilon f,u\rangle\,dt
 +\rho_\epsilon g(u)^2\,dt+dM,
\]
where
\[
 M(t)=2\sqrt\epsilon\int_0^t
       g(u(s))\langle R_\epsilon^*u(s),dW(s)\rangle.
\]
Global Lipschitz continuity implies
$|g(u)|^2\leq C(1+|u|^2)$.
Since $\|R_\epsilon\|_{\mathcal L(H)}\leq1$,
Young's and Poincar\'e's inequalities yield
\[
 2|g(u)\langle R_\epsilon f,u\rangle|
 \leq\frac{\nu}{2}\|u\|_{\mathcal V}^2
       +C|f|^2(1+|u|^2).
\]
Also,
\[
 \rho_\epsilon g(u)^2
 \leq C\rho_\epsilon(1+|u|^2).
\]

Let $\lambda_1>0$ be the first eigenvalue of $A$.
Since $\rho_\epsilon\to0$, choose
$0<\epsilon_*\leq1$, independently of $N$ and $T$, such that
\[
 C\rho_\epsilon\leq\frac{\nu\lambda_1}{2},
 \qquad0<\epsilon\leq\epsilon_*.
\]
Then $\lambda_1|u|^2\leq\|u\|_{\mathcal V}^2$ gives
\[
 \rho_\epsilon g(u)^2
 \leq C\rho_\epsilon
       +\frac{\nu}{2}\|u\|_{\mathcal V}^2.
\]
Substituting these estimates into the energy identity, we obtain
\[
 d|u|^2+\nu\|u\|_{\mathcal V}^2\,dt
 \leq C|f|^2(1+|u|^2)\,dt
       +C\rho_\epsilon\,dt+dM.
\]
Moreover,
\[
 \nu\|u\|_{\mathcal V}^2
 \geq\frac{\nu}{2}\|u\|_{\mathcal V}^2
      +\frac{\nu\lambda_1}{2}|u|^2.
\]
Consequently, with
$c=\frac{\nu}{2}\min\{1,\lambda_1\}>0$,
\begin{equation}\label{eq:NS-weighted-energy-start}
 d|u|^2+c|u|^2\,dt+c\|u\|_{\mathcal V}^2\,dt
 \leq C|f|^2|u|^2\,dt
       +C(|f|^2+\rho_\epsilon)\,dt+dM.
\end{equation}

Fix the coefficient $C$ in this inequality and define
\[
 p(t)=\exp\left(-C\int_0^t|f(s)|^2\,ds\right).
\]
Because $f\in\mathcal P_T^N$,
\[
 e^{-2CN}\leq p(t)\leq1,\qquad0\leq t\leq T.
\]
The identity $p'(t)=-C|f(t)|^2p(t)$ and the product rule give
\[
 d\bigl(p|u|^2\bigr)
 +cp|u|^2\,dt+cp\|u\|_{\mathcal V}^2\,dt
 \leq Cp(|f|^2+\rho_\epsilon)\,dt+p\,dM.
\]
Thus the derivative of $p$ cancels the term containing
$|f|^2|u|^2$.

We first stop the martingale before taking expectations
and subsequently remove the stopping using the resulting
bounds. Taking expectations and applying Gronwall's
inequality gives
\begin{align*}
 \mathbb E[p(t)|u(t)|^2]
 &\leq e^{-ct}|x|^2
 +C\mathbb E\int_0^t
       e^{-c(t-s)}p(s)|f(s)|^2\,ds\\
 &\quad+
 C\rho_\epsilon\mathbb E\int_0^t
       e^{-c(t-s)}p(s)\,ds.
\end{align*}
Since $p\leq1$, the control-energy bound gives
\[
 \int_0^t e^{-c(t-s)}p(s)|f(s)|^2\,ds
 \leq\int_0^T|f(s)|^2\,ds
 \leq2N
 \qquad\text{almost surely}.
\]
The other integral satisfies
\[
 \int_0^t e^{-c(t-s)}p(s)\,ds
 \leq\frac{1-e^{-ct}}{c}
 \leq\frac1c.
\]
Hence
\[
 \mathbb E[p(t)|u(t)|^2]
 \leq e^{-ct}|x|^2+2CN+\frac{C\rho_\epsilon}{c}.
\]
Using $p(t)\geq e^{-2CN}$, we conclude that
\[
 \mathbb E|u(t)|^2
 \leq e^{2CN}
       \left(e^{-ct}|x|^2+2CN+\frac{C\rho_\epsilon}{c}\right).
\]
The choice of $\epsilon_*$ bounds $\rho_\epsilon$ uniformly
for $0\leq\epsilon\leq\epsilon_*$. The last estimate therefore
proves \eqref{eq:NS-terminal-moment}, after enlarging the
constant in its statement.

For the deterministic estimate, take $\epsilon=0$ and
$f\in\mathcal S_T^N$. Then $u=z_f^x$, $\rho_0=0$ and $M=0$.
Integrating the same weighted inequality gives
\[
 p(t)|z_f^x(t)|^2
 +c\int_0^t p(s)\|z_f^x(s)\|_{\mathcal V}^2\,ds
 \leq |x|^2+C\int_0^t p(s)|f(s)|^2\,ds
 \leq |x|^2+2CN.
\]
Applying the lower bound on $p$ separately to the solution
norm and its time integral yields
\[
 \sup_{t\leq T}|z_f^x(t)|^2
 +\int_0^T\|z_f^x(t)\|_{\mathcal V}^2\,dt
 \leq
 \left(1+\frac1c\right)e^{2CN}
       \bigl(|x|^2+2CN\bigr).
\]
This bound is independent of $T$. Taking the supremum over
$|x|\leq R$ and $f\in\mathcal S_T^N$ proves
\eqref{eq:NS-skeleton-global-bound}.
\end{proof}

In the uncontrolled equation, the same calculation gives
\[
 \mathbb E|u_\epsilon^0(T)|^2\leq C,\qquad
 \frac1T\mathbb E\int_0^T\|u_\epsilon^0(t)\|_{\mathcal V}^2\,dt
       \leq C,\qquad T\geq1,
\]
uniformly for $\epsilon\leq\epsilon_*$. The compactness of
$\mathcal V\hookrightarrow H$ gives tightness of the time
averages. The Feller property and the Krylov-Bogoliubov
construction yield the invariant measures in
Theorem~\ref{thm:NS-stationary-LDP}.

Applying Theorem~\ref{thm:abstract-exponential-estimate} to the
terminal outputs $u_{\epsilon,f}^0(T)$, with
\eqref{eq:NS-terminal-moment} and the triangle inequality,
gives, for every $r>0$, a radius $R_r$ independent of $T$ such that
\begin{equation}\label{eq:NS-uniform-H-tail}
 \sup_{T>0}\mathbb P(|u_\epsilon^0(T)|\geq R_r)
       \leq e^{-r/\epsilon},\qquad0<\epsilon\leq\epsilon_*.
\end{equation}
Passing to the time-average limits, using the open set
$\{|x|>R_r\}$ and then enlarging the radius, gives
\begin{equation}\label{eq:NS-invariant-H-tail}
 \mu_\epsilon(|x|\geq2R_r)\leq e^{-r/\epsilon}.
\end{equation}
The controlled moment constants may grow faster than $\sqrt N$.
Thus \eqref{eq:NS-invariant-H-tail} is an estimate at every
exponential level, without a claimed Gaussian dependence on
the norm radius. It is precisely the stochastic tail hypothesis
of Theorem~\ref{thm:stationary-transfer}.

\subsection{Verification of the stationary LDP}
\paragraph{Deterministic properties.}
We first record a deterministic continuity fact used to prove
compactness of the finite time action sublevels. In this section,
$L^4$ denotes the closed subspace of divergence free, mean zero
vector fields in $L^4(\mathbb T^2;\mathbb R^2)$. For
$Z\in C([0,T];L^4)$, let
$\mathcal R_T(x,Z)=v+Z$, where
\begin{equation}\label{eq:NS-path-solution-map}
 \partial_t v+\nu Av+B(v+Z)=0,\qquad v(0)=x.
\end{equation}
If $Z(0)=0$, $x$ is the initial value of $v+Z$.

\begin{lemma}\label{lem:NS-path-map}
The map $\mathcal R_T$ is continuous from
$H\times C([0,T];L^4)$ to $C([0,T];H)$, uniformly on bounded
sets.
\end{lemma}

\begin{proof}
Testing \eqref{eq:NS-path-solution-map} with $v$ and using
$\|v\|_{L^4}^2\leq C|v|\|v\|_{\mathcal V}$ gives
\[
 \frac{d}{dt}|v|^2+\nu\|v\|_{\mathcal V}^2
     \leq C\|Z\|_{L^4}^4|v|^2+C\|Z\|_{L^4}^4.
\]
Thus $v$ is bounded in $L^\infty(0,T;H)\cap L^2(0,T;\mathcal V)$
on bounded sets of inputs. In particular, $v+Z$ is bounded in
$L^4(0,T;L^4)$ and $\partial_t v$ in $L^2(0,T;\mathcal V')$.
These estimates give existence by Galerkin approximation and
compactness in $L^2(0,T;H)$. The difference estimate below gives
uniqueness, and the usual energy-space continuity gives
$v\in C([0,T];H)$.
For two inputs, write $d=v_1-v_2$ and $\zeta=Z_1-Z_2$.
Skew symmetry of $b$ and the same interpolation inequality give
\begin{align*}
 \frac{d}{dt}|d|^2+\nu\|d\|_{\mathcal V}^2
 &\leq C(\|v_1+Z_1\|_{L^4}^4+
                    \|v_2+Z_2\|_{L^4}^4)|d|^2\\
 &\quad+C\|\zeta\|_{L^4}^2
       (\|v_1+Z_1\|_{L^4}^2+\|v_2+Z_2\|_{L^4}^2).
\end{align*}
Gronwall's inequality proves
\begin{equation}\label{eq:NS-path-map-continuity}
 \|\mathcal R_T(x_1,Z_1)-\mathcal R_T(x_2,Z_2)\|_{C([0,T];H)}
       \leq C_{T,R,K}\bigl(|x_1-x_2|
                             +\|Z_1-Z_2\|_{C([0,T];L^4)}\bigr)
\end{equation}
when $|x_i|\leq R$ and $\|Z_i\|_{C L^4}\leq K$.

\end{proof}

\begin{lemma}
\label{lem:NS-skeleton-continuity}
For every $T,R,N>0$, there exists $C_{T,R,N}>0$ such that
\begin{equation}\label{eq:NS-skeleton-continuity}
 \|z_f^x-z_h^y\|_{C([0,T];H)}
 \leq C_{T,R,N}
       \bigl(|x-y|+\|f-h\|_{L^2(0,T;H)}\bigr)
\end{equation}
for $|x|,|y|\leq R$ and $f,h\in\mathcal S_T^N$.
\end{lemma}

\begin{proof}
Put $w=z_f^x-z_h^y$. Subtracting the two skeleton equations
and taking the inner product with $w$, we obtain
\begin{align*}
 \frac12\frac{d}{dt}|w|^2+\nu\|w\|_{\mathcal V}^2
 &=-b(w,z_h^y,w)\\
 &\quad+
 \bigl(g(z_f^x)-g(z_h^y)\bigr)\langle f,w\rangle
 +g(z_h^y)\langle f-h,w\rangle.
\end{align*}
Here we used
\[
 B(z_f^x)-B(z_h^y)
 =B(z_f^x,w)+B(w,z_h^y),
 \qquad b(z_f^x,w,w)=0.
\]
The two-dimensional interpolation inequality gives
\[
 |b(w,z_h^y,w)|
 \leq C\|z_h^y\|_{\mathcal V}|w|\|w\|_{\mathcal V}
 \leq\frac{\nu}{2}\|w\|_{\mathcal V}^2
       +C_\nu\|z_h^y\|_{\mathcal V}^2|w|^2.
\]
Global Lipschitz continuity of $g$ gives
\[
 \left|
 \bigl(g(z_f^x)-g(z_h^y)\bigr)\langle f,w\rangle
 \right|
 \leq C|f||w|^2.
\]
Moreover, \eqref{eq:NS-skeleton-global-bound} implies
\[
 \sup_{t\leq T}|g(z_h^y(t))|^2
 \leq C_g(1+C_{R,N}).
\]
Consequently, Young's inequality yields
\[
 |g(z_h^y)\langle f-h,w\rangle|
 \leq\frac12|w|^2+C(1+C_{R,N})|f-h|^2.
\]
Combining these estimates and using
$2|f|\leq1+|f|^2$, we obtain
\[
 \frac{d}{dt}|w|^2+\nu\|w\|_{\mathcal V}^2
 \leq
 C\bigl(1+\|z_h^y\|_{\mathcal V}^2+|f|^2\bigr)|w|^2
 +C(1+C_{R,N})|f-h|^2.
\]
The coefficient of $|w|^2$ has a uniform integral bound:
\[
 \int_0^T
 \bigl(1+\|z_h^y(t)\|_{\mathcal V}^2+|f(t)|^2\bigr)\,dt
 \leq T+C_{R,N}+2N.
\]
Since $w(0)=x-y$, Gronwall's inequality gives
\[
 \sup_{t\leq T}|w(t)|^2
 \leq C_{T,R,N}
       \left(|x-y|^2+\|f-h\|_{L^2(0,T;H)}^2\right).
\]
Taking square roots proves
\eqref{eq:NS-skeleton-continuity}.
\end{proof}

\begin{lemma}
\label{lem:NS-action-compactness}
For every $T>0$, $x\in H$ and $r\geq0$, the set
\[
 K_T^x(r)
 :=\{z\in C([0,T];H):I_T^x(z)\leq r\}
\]
is compact in $C([0,T];H)$.
\end{lemma}

\begin{proof}
We use the compactness of additive controlled Navier-Stokes
paths with bounded $L^2(0,T;H)$ forcing, established through
the compact linear control map and the continuous solution
map in \cite[Section~3]{cerrai2022large}. This deterministic
argument applies for every fixed $\nu>0$.

Take a sequence $z_n\in K_T^x(r)$. By the definition of
the action, choose controls $f_n$ such that
\[
 z_n=z_{f_n}^x,\qquad
 \frac12\|f_n\|_{L^2(0,T;H)}^2\leq r+\frac1n.
\]
Since $f_n\in\mathcal S_T^{r+1}$,
\eqref{eq:NS-skeleton-global-bound} gives
\begin{equation}\label{eq:NS-compactness-energy}
 \sup_n\left[
 \sup_{0\leq t\leq T}|z_n(t)|^2
 +\int_0^T\|z_n(t)\|_{\mathcal V}^2\,dt
 \right]\leq C_{x,r}.
\end{equation}
The growth bound $|g(v)|^2\leq C_g(1+|v|^2)$ therefore implies
\[
 \int_0^T|g(z_n(t))f_n(t)|^2\,dt
 \leq
 C_g\left(1+\sup_{t\leq T}|z_n(t)|^2\right)
       \|f_n\|_{L^2(0,T;H)}^2
 \leq C_{x,r}.
\]
Thus each $z_n$ solves an additive controlled equation with
the same initial value $x$ and uniformly bounded forcing.
The cited compactness result gives a subsequence such that
\[
 z_n\longrightarrow z
       \quad\text{in }C([0,T];H).
\]
The controls are bounded in $L^2(0,T;H)$, so, after passing
to a further subsequence,
\[
 f_n\rightharpoonup f
       \quad\text{in }L^2(0,T;H).
\]

We first verify the regularity of the limiting path.
Estimate \eqref{eq:NS-compactness-energy} also bounds $z_n$
in $L^2(0,T;\mathcal V)$. Weak compactness gives a further
weakly convergent subsequence in that space. Its weak limit
is $z$, because $\mathcal V\hookrightarrow H$ continuously
and $z_n\to z$ strongly in $L^2(0,T;H)$. Consequently,
\[
 z\in C([0,T];H)\cap L^2(0,T;\mathcal V),
 \qquad
 z_n\rightharpoonup z
       \quad\text{in }L^2(0,T;\mathcal V).
\]

We next identify the limiting forcing. Global Lipschitz
continuity gives
\[
 \sup_{t\leq T}|g(z_n(t))-g(z(t))|
 \leq L_g\|z_n-z\|_{C([0,T];H)}
 \longrightarrow0.
\]
For any $\eta\in L^2(0,T;H)$,
\begin{align*}
 &\int_0^T
       \langle g(z_n)f_n-g(z)f,\eta\rangle\,dt\\
 &\quad=
 \int_0^T[g(z_n)-g(z)]\langle f_n,\eta\rangle\,dt
 +\int_0^T\langle f_n-f,g(z)\eta\rangle\,dt.
\end{align*}
The first integral tends to zero, since its absolute value
is bounded by
\[
 \sup_{t\leq T}|g(z_n(t))-g(z(t))|\,
 \|f_n\|_{L^2(0,T;H)}
 \|\eta\|_{L^2(0,T;H)}.
\]
The second tends to zero by weak convergence of $f_n$:
the function $g(z)$ is bounded on $[0,T]$, so
$g(z)\eta\in L^2(0,T;H)$. Hence
\[
 g(z_n)f_n\rightharpoonup g(z)f
       \quad\text{in }L^2(0,T;H).
\]

We now pass to the limit in the equation for $z_n$,
tested against a smooth divergence-free periodic function
$\phi\in H$. The linear terms converge by the strong
convergence in $C([0,T];H)$ and the weak convergence in
$L^2(0,T;\mathcal V)$. For the nonlinear term, fix a smooth divergence-free test
function $\phi$. Using
$b(u,v,w)=-b(u,w,v)$, we obtain
\begin{align*}
 &|b(z_n,z_n,\phi)-b(z,z,\phi)|\\
 &\quad\leq
 |b(z_n-z,\phi,z_n)|+|b(z,\phi,z_n-z)|\\
 &\quad\leq
 \|\nabla\phi\|_{L^\infty}
 (|z_n|+|z|)|z_n-z|.
\end{align*}
Since $z_n\to z$ in $C([0,T];H)$, it follows that
\[
 \int_0^t b(z_n(s),z_n(s),\phi)\,ds
 \longrightarrow
 \int_0^t b(z(s),z(s),\phi)\,ds
\]
uniformly for $t\in[0,T]$.
Together with the convergence of the forcing, this yields,
for every $t\in[0,T]$,
\begin{align*}
 \langle z(t),\phi\rangle
 &+\nu\int_0^t
       \langle A^{1/2}z(s),A^{1/2}\phi\rangle\,ds\\
 &+\int_0^t b(z(s),z(s),\phi)\,ds\\
 &=\langle x,\phi\rangle
   +\int_0^t g(z(s))\langle f(s),\phi\rangle\,ds.
\end{align*}
Also $z(0)=x$. Thus $z$ belongs to the energy class and
satisfies the weak formulation of the skeleton equation
with initial value $x$ and control $f$. Uniqueness gives
$z=z_f^x$.

Finally, the definition of the action and weak lower
semicontinuity of the control norm give
\[
 I_T^x(z)
 \leq\frac12\|f\|_{L^2(0,T;H)}^2
 \leq\liminf_{n\to\infty}
       \frac12\|f_n\|_{L^2(0,T;H)}^2
 \leq r.
\]
Every sequence in $K_T^x(r)$ therefore has a subsequence
converging in $C([0,T];H)$ to an element of $K_T^x(r)$.
This proves compactness.
\end{proof}

\begin{proposition}
\label{prop:NS-stationary-input}
The skeleton \eqref{eq:NS-skeleton} satisfies
Assumption~\ref{ass:stationary-input} with state space $H$
and equilibrium $0$.

In particular, for every $t_0,R,N>0$,
\begin{equation}\label{eq:NS-positive-time-bound}
 \sup_{|x|\leq R}\sup_{f\in\mathcal S_{t_0}^N}
       \|z_f^x(t_0)\|_{\mathcal V}^2
 \leq
 \frac{2C_{R,N}}{t_0}
 +\frac{2C_g}{\nu}(1+C_{R,N})N,
\end{equation}
where $C_{R,N}$ is a constant in
\eqref{eq:NS-skeleton-global-bound}, independent of $t_0$.
\end{proposition}

\begin{proof}
Lemmas~\ref{lem:NS-skeleton-continuity} and
\ref{lem:NS-action-compactness} verify
Assumption~\ref{ass:stationary-input}(i).
The bound independent of duration in part~(ii) follows from
\eqref{eq:NS-skeleton-global-bound}.

We prove the remaining spatial compactness requirement.
Fix $t_0,R,N>0$, $|x|\leq R$ and
$f\in\mathcal S_{t_0}^N$. By
\eqref{eq:NS-skeleton-global-bound},
\[
 \sup_{t\leq t_0}|z_f^x(t)|^2
 +\int_0^{t_0}\|z_f^x(t)\|_{\mathcal V}^2\,dt
 \leq C_{R,N}.
\]
Therefore there exists $a\in(0,t_0/2)$ such that
\[
 z_f^x(a)\in\mathcal V,\qquad
 \|z_f^x(a)\|_{\mathcal V}^2
 \leq\frac{2C_{R,N}}{t_0}.
\]

On $[a,t_0]$, the forcing $g(z_f^x)f$ belongs to
$L^2(a,t_0;H)$. Starting from the value in $\mathcal V$
at time $a$, the following calculation is justified by
finite Stokes projections with this given forcing.

Taking the $H$ inner product of the equation with $Az_f^x$
and using
\[
 \langle\partial_tz_f^x,Az_f^x\rangle
 =\frac12\frac{d}{dt}\|z_f^x\|_{\mathcal V}^2,
 \qquad
 \langle B(z_f^x),Az_f^x\rangle=0,
\]
we obtain
\begin{align*}
 \frac12\frac{d}{dt}\|z_f^x\|_{\mathcal V}^2
 +\nu|Az_f^x|^2
 &=g(z_f^x)\langle f,Az_f^x\rangle\\
 &\leq |g(z_f^x)|\,|f|\,|Az_f^x|\\
 &\leq\frac{\nu}{2}|Az_f^x|^2
       +\frac{|g(z_f^x)|^2}{2\nu}|f|^2.
\end{align*}
Subtracting $\frac{\nu}{2}|Az_f^x|^2$ and multiplying by
$2$ gives
\[
 \frac{d}{dt}\|z_f^x\|_{\mathcal V}^2
 +\nu|Az_f^x|^2
 \leq\frac{|g(z_f^x)|^2}{\nu}|f|^2.
\]
Since
\[
 |g(z_f^x(t))|^2
 \leq C_g(1+|z_f^x(t)|^2)
 \leq C_g(1+C_{R,N}),
\]
we conclude that
\begin{equation}\label{eq:NS-enstrophy-skeleton}
 \frac{d}{dt}\|z_f^x\|_{\mathcal V}^2
 +\nu|Az_f^x|^2
 \leq\frac{C_g}{\nu}(1+C_{R,N})|f|^2.
\end{equation}
Integration from $a$ to any $t\in[a,t_0]$ yields
\begin{align*}
 \|z_f^x(t)\|_{\mathcal V}^2
 +\nu\int_a^t|Az_f^x(s)|^2\,ds
 &\leq
 \|z_f^x(a)\|_{\mathcal V}^2
 +\frac{C_g}{\nu}(1+C_{R,N})
       \int_a^t|f(s)|^2\,ds\\
 &\leq
 \frac{2C_{R,N}}{t_0}
 +\frac{2C_g}{\nu}(1+C_{R,N})N.
\end{align*}
Taking $t=t_0$ proves \eqref{eq:NS-positive-time-bound}.
Its right-hand side is independent of $x$, $f$ and the
selected time $a$. Hence the endpoints belong to a common
bounded subset of $\mathcal V$.
The compact embedding $\mathcal V\Subset H$ proves
\eqref{eq:det-positive-time-compactness}.

Finally, with zero control, the energy identity gives
\[
 \frac{d}{dt}|z_0^x(t)|^2
 +2\nu\|z_0^x(t)\|_{\mathcal V}^2=0.
\]
Poincar\'e's and Gronwall's inequalities imply
\[
 |z_0^x(t)|\leq e^{-\nu\lambda_1t}|x|.
\]
Thus the uncontrolled flow converges to $0$ uniformly on
bounded initial sets, and therefore on compact initial sets.
This verifies Assumption~\ref{ass:stationary-input}(iii)
and completes the proof.
\end{proof}

\paragraph{Finite time uniform LDP.}
\begin{lemma}
\label{lem:NS-controlled-convergence}
Under Assumption~\ref{ass:NS}, for every $T,R,N>0$,
there exists $C_{T,R,N}>0$ such that
\begin{equation}\label{eq:NS-controlled-convergence}
\begin{aligned}
 &\sup_{|x|\leq R}
 \sup_{f\in\mathcal P_T^N}
 \mathbb E\left[
 \sup_{0\leq t\leq T}
       |u_{\epsilon,f}^x(t)-z_f^x(t)|^2
 +
 \int_0^T
       \|u_{\epsilon,f}^x(t)-z_f^x(t)\|_{\mathcal V}^2\,dt
 \right]
 \\
 &\qquad\leq
 C_{T,R,N}\bigl(\rho_\epsilon+\theta_\epsilon^2\bigr),
 \qquad 0<\epsilon\leq\epsilon_*,
\end{aligned}
\end{equation}
where $\epsilon_*$ is as in
Lemma~\ref{lem:NS-controlled-moments}.
\end{lemma}

\begin{proof}
Fix $|x|\leq R$ and $f\in\mathcal P_T^N$, and write
\[
 u=u_{\epsilon,f}^x,\qquad z=z_f^x,\qquad w=u-z.
\]
For a random control, the deterministic skeleton estimates
apply to each realisation. Thus
\begin{equation}\label{eq:NS-comparison-skeleton-bound}
 \sup_{0\leq t\leq T}|z(t)|^2
 +\int_0^T\|z(t)\|_{\mathcal V}^2\,dt
 \leq C_{R,N}
 \qquad\text{almost surely}.
\end{equation}
Moreover, the growth bound on $g$ and
\eqref{eq:NS-terminal-moment} imply
\begin{equation}\label{eq:NS-comparison-diffusion-bound}
 \mathbb E\int_0^T|g(u(t))|^2\,dt
 \leq C_g\int_0^T
             \bigl(1+\mathbb E|u(t)|^2\bigr)\,dt
 \leq C_{T,R,N}.
\end{equation}
Indeed, the restriction of $f$ to $[0,t]$ belongs to
$\mathcal P_t^N$, so the terminal estimate applies at every
$t\leq T$.

Subtracting the skeleton equation from the controlled stochastic
equation gives
\[
\begin{aligned}
 dw+\bigl[\nu Aw+B(u)-B(z)\bigr]\,dt
 &=
 \bigl[g(u)-g(z)\bigr]R_\epsilon f\,dt\\
 &\quad+
 g(z)(R_\epsilon-I)f\,dt
 +\sqrt\epsilon\,g(u)R_\epsilon\,dW,
 \qquad w(0)=0.
\end{aligned}
\]
Since
\[
 B(u)-B(z)=B(u,w)+B(w,z),
 \qquad b(u,w,w)=0,
\]
It\^o's formula yields
\begin{align*}
 d|w|^2+2\nu\|w\|_{\mathcal V}^2\,dt
 &=-2b(w,z,w)\,dt\\
 &\quad+
 2\bigl[g(u)-g(z)\bigr]
                  \langle R_\epsilon f,w\rangle\,dt\\
 &\quad+
 2g(z)\langle(R_\epsilon-I)f,w\rangle\,dt\\
 &\quad+\rho_\epsilon g(u)^2\,dt+dM_w,
\end{align*}
where
\[
 M_w(t)=2\sqrt\epsilon\int_0^t
       g(u(s))
       \langle R_\epsilon^*w(s),dW(s)\rangle.
\]

The trilinear estimate gives
\[
 2|b(w,z,w)|
 \leq\frac{\nu}{2}\|w\|_{\mathcal V}^2
       +C_\nu\|z\|_{\mathcal V}^2|w|^2.
\]
Global Lipschitz continuity of $g$ and
$\|R_\epsilon\|_{\mathcal L(H)}\leq1$ give
\[
 2|g(u)-g(z)|\,|\langle R_\epsilon f,w\rangle|
 \leq C|f|\,|w|^2
 \leq C(1+|f|^2)|w|^2.
\]
For the covariance difference, the definition of
$\theta_\epsilon$ gives
\begin{align*}
 |\langle(R_\epsilon-I)f,w\rangle|
 &=
 |\langle A^{-1/2}(R_\epsilon-I)f,
                         A^{1/2}w\rangle|\\
 &\leq\theta_\epsilon|f|\,
                         \|w\|_{\mathcal V}.
\end{align*}
Consequently, by
\eqref{eq:NS-comparison-skeleton-bound},
\[
 2|g(z)\langle(R_\epsilon-I)f,w\rangle|
 \leq\frac{\nu}{2}\|w\|_{\mathcal V}^2
       +C_{R,N}\theta_\epsilon^2|f|^2.
\]
Finally, since $u=z+w$,
\[
 \rho_\epsilon g(u)^2
 \leq C\rho_\epsilon|w|^2+C_{R,N}\rho_\epsilon.
\]
Combining these inequalities yields
\begin{equation}\label{eq:NS-comparison-energy}
\begin{aligned}
 d|w|^2+\nu\|w\|_{\mathcal V}^2\,dt
 &\leq
 C\bigl(1+\|z\|_{\mathcal V}^2
                 +|f|^2+\rho_\epsilon\bigr)|w|^2\,dt\\
 &\quad+
 C_{R,N}\bigl(\theta_\epsilon^2|f|^2
                         +\rho_\epsilon\bigr)\,dt+dM_w.
\end{aligned}
\end{equation}

Define the adapted weight
\[
 p(t)=\exp\left[
 -C\int_0^t
 \bigl(1+\|z(s)\|_{\mathcal V}^2
                 +|f(s)|^2+\rho_\epsilon\bigr)\,ds
 \right],
\]
using the coefficient $C$ in
\eqref{eq:NS-comparison-energy}.
The control energy bound and
\eqref{eq:NS-comparison-skeleton-bound} imply
\[
 e^{-C_{T,R,N}}\leq p(t)\leq1,
 \qquad 0\leq t\leq T,
\]
uniformly for $0<\epsilon\leq\epsilon_*$.
The product rule therefore gives
\begin{equation}\label{eq:NS-comparison-weighted}
\begin{aligned}
 p(t)|w(t)|^2
 +\nu\int_0^t p(s)\|w(s)\|_{\mathcal V}^2\,ds
 &\leq
 C_{T,R,N}\bigl(\theta_\epsilon^2+\rho_\epsilon\bigr)
 +\int_0^t p(s)\,dM_w(s).
\end{aligned}
\end{equation}

We estimate the martingale using
\[
 \epsilon|R_\epsilon^*w|^2
 \leq\epsilon\operatorname{Tr}C_\epsilon\,|w|^2
 =\rho_\epsilon|w|^2.
\]
Hence
\[
 d\langle M_w\rangle(t)
 \leq4\rho_\epsilon|g(u(t))|^2|w(t)|^2\,dt.
\]
The BDG inequality and Young's inequality give
\begin{align*}
 &\mathbb E\sup_{0\leq t\leq T}
       \left|\int_0^t p(s)\,dM_w(s)\right|\\
 &\quad\leq
 C\sqrt{\rho_\epsilon}\,
 \mathbb E\left[
   \sup_{0\leq t\leq T}p(t)|w(t)|^2
   \int_0^T p(s)|g(u(s))|^2\,ds
 \right]^{1/2}\\
 &\quad\leq
 \frac14\mathbb E\sup_{0\leq t\leq T}p(t)|w(t)|^2
 +C\rho_\epsilon
       \mathbb E\int_0^T p(s)|g(u(s))|^2\,ds\\
 &\quad\leq
 \frac14\mathbb E\sup_{0\leq t\leq T}p(t)|w(t)|^2
 +C_{T,R,N}\rho_\epsilon,
\end{align*}
where the last inequality follows from
\eqref{eq:NS-comparison-diffusion-bound}.

Taking the supremum in
\eqref{eq:NS-comparison-weighted}, and also using that
inequality at $t=T$, we obtain
\begin{align*}
 &\mathbb E\left[
   \sup_{0\leq t\leq T}p(t)|w(t)|^2
   +\nu\int_0^T p(t)\|w(t)\|_{\mathcal V}^2\,dt
 \right]\\
 &\quad\leq
 C_{T,R,N}\bigl(\theta_\epsilon^2+\rho_\epsilon\bigr)
 +\frac12\mathbb E
       \sup_{0\leq t\leq T}p(t)|w(t)|^2.
\end{align*}
Absorbing the final term and using the lower bound on $p$
proves \eqref{eq:NS-controlled-convergence}.

The martingale estimates are first applied after stopping.
Their bounds are independent of the stopping level, which
allows the stopping to be removed.
\end{proof}

\begin{proposition}\label{prop:NS-finite-time-LDP}
Under Assumption~\ref{ass:NS}, for every $T,R>0$, the laws
of $u_\epsilon^x$ satisfy the Freidlin-Wentzell uniform LDP
in $C([0,T];H)$ as $\epsilon\downarrow0$, uniformly for
$|x|\leq R$, with speed $1/\epsilon$ and good rate functions
$I_T^x$ defined in \eqref{eq:NS-action}.
\end{proposition}

\begin{proof}
We prove the equicontinuous uniform Laplace principle and
then use its equivalence with the Freidlin-Wentzell uniform
LDP.

Put $\mathsf X_T=C([0,T];H)$.
Lemma~\ref{lem:NS-controlled-convergence} and Markov's
inequality imply, for every $N>0$ and $a>0$,
\begin{equation}\label{eq:NS-uniform-controlled-probability}
\begin{aligned}
 &\sup_{|x|\leq R}\sup_{f\in\mathcal P_T^N}
 \mathbb P\left(
       \|u_{\epsilon,f}^x-z_f^x\|_{\mathsf X_T}>a
 \right)\\
 &\qquad\leq
 \frac{C_{T,R,N}}{a^2}
       \bigl(\rho_\epsilon+\theta_\epsilon^2\bigr)
 \longrightarrow0.
\end{aligned}
\end{equation}
Thus the controlled convergence holds uniformly over both
bounded initial values and controls with bounded energy.

Represent $W$ by its coordinate Brownian motions.
The strong solution admits a measurable solution map, and
shifting the Brownian input by the integral of $f$ produces
$u_{\epsilon,f}^x$. The corresponding deterministic controlled
output is $z_f^x$. Therefore
\eqref{eq:NS-uniform-controlled-probability} verifies the
controlled convergence hypothesis of the EULP criterion
in \cite[Theorem~2.13]{salins2019equivalences}.

To make the uniformity over test functions explicit, let
$\mathscr H\subset C_b(\mathsf X_T)$ be equibounded and
uniformly equicontinuous. Write
\[
 L=\sup_{h\in\mathscr H}\|h\|_\infty<\infty
\]
and let
\[
 \omega_{\mathscr H}(a)
 :=
 \sup_{h\in\mathscr H}
 \sup_{\|\varphi-\psi\|_{\mathsf X_T}\leq a}
       |h(\varphi)-h(\psi)|.
\]
Then $\omega_{\mathscr H}(a)\to0$ as $a\downarrow0$.
For every fixed $N>0$,
\begin{align*}
 &\sup_{|x|\leq R}
 \sup_{f\in\mathcal P_T^N}
 \sup_{h\in\mathscr H}
 \mathbb E\left|
       h(u_{\epsilon,f}^x)-h(z_f^x)
 \right|\\
 &\qquad\leq
 \omega_{\mathscr H}(a)
 +\frac{2LC_{T,R,N}}{a^2}
       \bigl(\rho_\epsilon+\theta_\epsilon^2\bigr).
\end{align*}
First letting $\epsilon\downarrow0$ and then $a\downarrow0$
shows that this expected difference tends to zero uniformly
in $x$, $f$ and $h$.

The cited EULP criterion consequently gives
\begin{equation}\label{eq:NS-EULP}
\begin{aligned}
 \lim_{\epsilon\downarrow0}
 \sup_{|x|\leq R}\sup_{h\in\mathscr H}
 \left|
 -\epsilon\log
       \mathbb E\exp\left[-\frac{h(u_\epsilon^x)}{\epsilon}\right]
 -
 \inf_{\varphi\in\mathsf X_T}
       \bigl\{I_T^x(\varphi)+h(\varphi)\bigr\}
 \right|
 =0.
\end{aligned}
\end{equation}
Here the rate function is $I_T^x$, since its defining
variational problem is precisely
\[
 \inf_{f\in L^2(0,T;H)}
 \left\{
       \frac12\int_0^T|f(t)|^2\,dt+h(z_f^x)
 \right\}
 =
 \inf_{\varphi\in\mathsf X_T}
       \bigl\{I_T^x(\varphi)+h(\varphi)\bigr\}.
\]

By \cite[Theorem~2.10]{salins2019equivalences},
\eqref{eq:NS-EULP} is equivalent to the Freidlin-Wentzell
uniform LDP on the same bounded initial set.
The finite time action sublevels were proved compact above,
so each $I_T^x$ is a good rate function.
\end{proof}

\paragraph{Finite path approximation.}
The compactness argument above also shows that the infimum in
\eqref{eq:NS-action} is attained whenever it is finite. Since
$g\geq g_*$, the path determines its control uniquely through
\eqref{eq:NS-skeleton}, interpreted in distributions. Controls
on overlapping intervals are therefore consistent. In particular,
a finite-action infinite past path has a single control in
$L^2(-\infty,0;H)$, and its action is half the squared norm
of that control.

\begin{lemma}\label{lem:NS-small-connection}
For every $y\in\mathcal V$, there is a path $w_y$ from $0$
to $y$ in time one such that
\begin{equation}\label{eq:NS-small-connection}
 I_1^0(w_y)
 \leq C\left(\|y\|_{\mathcal V}^2
                    +\|y\|_{\mathcal V}^4\right),
\end{equation}
where $C$ is independent of $y$.
Consequently, Assumption~\ref{ass:finite-path-approximation}
holds for \eqref{eq:NS-skeleton}.
\end{lemma}

\begin{proof}
For $y\in\mathcal V$, set
\[
 w_y(t):=tS(1-t)y,\qquad0\leq t\leq1.
\]
Then $w_y(0)=0$ and $w_y(1)=y$.
For $t>0$, the semigroup satisfies
\[
 \frac{d}{dt}\|S(t)y\|_{\mathcal V}^2
 =-2\nu|AS(t)y|^2.
\]
Integrating and using continuity in $\mathcal V$ at
$t=0$, we obtain
\[
 \int_0^1|AS(t)y|^2\,dt
 =
 \frac{\|y\|_{\mathcal V}^2
       -\|S(1)y\|_{\mathcal V}^2}{2\nu}
 \leq\frac1{2\nu}\|y\|_{\mathcal V}^2.
\]
Consequently,
\begin{align*}
 \int_0^1|Aw_y(t)|^2\,dt
 &=\int_0^1t^2|AS(1-t)y|^2\,dt\\
 &\leq\int_0^1|AS(t)y|^2\,dt
 \leq\frac1{2\nu}\|y\|_{\mathcal V}^2.
\end{align*}
Semigroup contractivity in $\mathcal V$ gives
\[
 \sup_{0\leq t\leq1}\|w_y(t)\|_{\mathcal V}
 \leq\|y\|_{\mathcal V}.
\]
Moreover, for $0<t<1$,
\[
 w_y'(t)=S(1-t)y+\nu Aw_y(t).
\]
It follows that
\begin{align*}
 \|w_y'\|_{L^2(0,1;H)}
 &\leq
 \|S(1-\cdot)y\|_{L^2(0,1;H)}
 +\nu\|Aw_y\|_{L^2(0,1;H)}\\
 &\leq |y|+\sqrt{\frac{\nu}{2}}\|y\|_{\mathcal V}
 \leq C_\nu\|y\|_{\mathcal V},
\end{align*}
where the last inequality uses Poincar\'e's inequality.
Thus
\[
 w_y\in W^{1,2}(0,1;H)\cap L^2(0,1;D(A))
                  \cap C([0,1];\mathcal V).
\]

For $w\in D(A)$, the two-dimensional estimate
\[
 |B(w)|
 \leq C|w|^{1/2}\|w\|_{\mathcal V}|Aw|^{1/2}
\]
implies
\begin{align*}
 \int_0^1|B(w_y(t))|^2\,dt
 &\leq C\int_0^1
       |w_y(t)|\,\|w_y(t)\|_{\mathcal V}^2
       |Aw_y(t)|\,dt\\
 &\leq C\|y\|_{\mathcal V}^3
       \left(\int_0^1|Aw_y(t)|^2\,dt\right)^{1/2}\\
 &\leq C_\nu\|y\|_{\mathcal V}^4.
\end{align*}
Here we used Poincar\'e's inequality, the preceding
supremum bound and Cauchy-Schwarz in time.

Define, for almost every $t\in(0,1)$,
\[
 f_y(t):=
 \frac{w_y'(t)+\nu Aw_y(t)+B(w_y(t))}
      {g(w_y(t))}.
\]
Since $g\geq g_*>0$, the preceding estimates give
\[
 \|f_y\|_{L^2(0,1;H)}
 \leq C\left(\|y\|_{\mathcal V}
                    +\|y\|_{\mathcal V}^2\right).
\]
By construction, $w_y$ solves \eqref{eq:NS-skeleton}
with initial value $0$ and control $f_y$. Hence
\[
 I_1^0(w_y)
 \leq\frac12\int_0^1|f_y(t)|^2\,dt
 \leq C\left(\|y\|_{\mathcal V}^2
                    +\|y\|_{\mathcal V}^4\right),
\]
which proves \eqref{eq:NS-small-connection}.

We now verify Assumption~\ref{ass:finite-path-approximation}.
Fix $\eta>0$ and an admissible infinite past path
$z\in C((-\infty,0];H)$ with finite action and
\[
 \lim_{t\to-\infty}|z(t)|=0.
\]
Choose a control $f\in L^2(-\infty,0;H)$ producing $z$
such that
\[
 \frac12\int_{-\infty}^0|f(t)|^2\,dt
 \leq I_{-\infty}(z)+\frac{\eta}{2}.
\]
For $T>1$, put
\[
 e_T:=|z(-T)|^2+
             \int_{-T}^{-T+1}|f(t)|^2\,dt.
\]
Then $e_T\to0$ as $T\to\infty$.

For all sufficiently large $T$, convergence of $z$ to $0$
gives $|z(t)|\leq1$ on $[-T,-T+1]$, and therefore
\[
 |g(z(t))|\leq |g(0)|+L_g
 \qquad\text{on }[-T,-T+1].
\]
The energy identity, Poincar\'e's inequality and Young's
inequality yield
\begin{align*}
 \frac12\frac{d}{dt}|z(t)|^2
       +\nu\|z(t)\|_{\mathcal V}^2
 &=g(z(t))\langle f(t),z(t)\rangle\\
 &\leq\frac{\nu}{2}\|z(t)\|_{\mathcal V}^2
       +\frac{(|g(0)|+L_g)^2}{2\nu\lambda_1}|f(t)|^2.
\end{align*}
Integrating over $[-T,-T+1]$ gives
\[
 \int_{-T}^{-T+1}\|z(t)\|_{\mathcal V}^2\,dt
 \leq C e_T.
\]
We can therefore choose
$t_T\in(-T+1/2,-T+1)$ such that $z(t_T)\in\mathcal V$ and
\[
 \|z(t_T)\|_{\mathcal V}^2
 \leq2\int_{-T+1/2}^{-T+1}
             \|z(t)\|_{\mathcal V}^2\,dt
 \leq2C e_T.
\]

Connect $0$ to $z(t_T)$ in time one using the path
constructed above, and then follow $z$ from $t_T$ to $0$.
More precisely, define
\[
 v_T(t):=
 \begin{cases}
  w_{z(t_T)}(t),&0\leq t\leq1,\\
  z(t_T+t-1),&1\leq t\leq1-t_T.
 \end{cases}
\]
Then $v_T\in C([0,1-t_T];H)$,
$v_T(0)=0$ and $v_T(1-t_T)=z(0)$.
Using $f_{z(t_T)}$ on the first interval and the shifted
control $f$ on the second, we obtain
\begin{align*}
 I_{1-t_T}^0(v_T)
 &\leq
 \frac12\int_0^1|f_{z(t_T)}(t)|^2\,dt
 +\frac12\int_{t_T}^0|f(t)|^2\,dt\\
 &\leq C(e_T+e_T^2)+I_{-\infty}(z)+\frac{\eta}{2}.
\end{align*}
Since $e_T\to0$, choosing $T$ sufficiently large gives
\[
 I_{1-t_T}^0(v_T)\leq I_{-\infty}(z)+\eta.
\]
This is Assumption~\ref{ass:finite-path-approximation}.
\end{proof}
Thus Proposition~\ref{prop:deterministic-compactness} and the
finite-path approximation Assumption~\ref{ass:finite-path-approximation} identify $V_{\rm NS}$ with the
infinite past quasipotential and give compactness of its
sublevels. The positive time compactness used there has
already been verified through
\eqref{eq:NS-enstrophy-skeleton}. Neither this step nor the
connection in Lemma~\ref{lem:NS-small-connection} requires
an $H^1$ estimate for the stochastic solution.

\begin{proof}[Proof of Theorem~\ref{thm:NS-stationary-LDP}]
The deterministic bounds, continuity, finite time compactness,
positive time spatial compactness and attraction required by
Assumption~\ref{ass:stationary-input} were verified above.
Lemma~\ref{lem:NS-small-connection} verifies
Assumption~\ref{ass:finite-path-approximation}.
Proposition~\ref{prop:NS-finite-time-LDP} supplies the finite
time uniform LDP over bounded initial sets, and
\eqref{eq:NS-invariant-H-tail} supplies the stationary norm
estimate \eqref{eq:stationary-norm-tail}. Theorem~\ref{thm:stationary-transfer}
therefore gives the stated LDP with good rate function
$V_{\rm NS}$.
\end{proof}

\begin{remark}\label{rem:NS-scope}
In the additive case, the stationary result of
Cerrai-Paskal~\cite{cerrai2022large} assumes $\beta>2$
and $\epsilon\delta(\epsilon)^{-2/\beta}\to0$.
Their compact exponential estimate uses the stochastic
$H^1$ energy. The present stationary theorem requires only
an $H$-norm tail. Its controlled estimate uses
$\epsilon\operatorname{Tr}C_\epsilon\to0$, giving
\eqref{eq:NS-correlation-scale} for $\beta>1$.

The example $g(u)=\sqrt{1+|u|^2}$ satisfies
Assumption~\ref{ass:NS}. Its diffusion is unbounded, and
the argument uses no contraction between two stochastic
solutions. Smallness is imposed on the noise parameter,
not on the viscosity. The covariance condition
\eqref{eq:NS-correlation-scale} is a sufficient condition;
no optimality is asserted. For a general Lipschitz
operator-valued diffusion, or a pointwise multiplicative
coefficient, the controlled comparison, spatial compactness
and the connection estimate would require separate
verification. They are not consequences of global
Lipschitz continuity alone.

For comparison, the additional radial noise bound in
\cite[Assumption~\textup{[A4]}]{liu2026stationary} is not satisfied
by $g(u)=\sqrt{1+|u|^2}$. If $e_1$ is a normalised Stokes
eigenvector and $u=ae_1$, then
\[
 |[g(u)R_\epsilon]^*u|^2
   =\frac{a^2(1+a^2)}{1+\delta(\epsilon)\lambda_1^\beta},
\]
which cannot be bounded by a constant multiple of $|u|^2$
for all $a$. Thus global Lipschitz continuity permits amplitudes
outside that radial bound even at a fixed covariance.
\end{remark}

\section{Reaction-diffusion systems with linear growth noise}
\label{sec:reaction-diffusion}

This section extends the stationary LDP of
Cerrai-R\"ockner~\cite{cerrai2005large} by removing their additional
growth condition for the unbounded-noise upper bound. We keep their
covariance, dissipativity and quasipotential hypotheses. In particular,
we do not replace their controllability assumptions by the statement
that the noise is merely Lipschitz.

\subsection{Setting and precise scope}

We use the dimension, elliptic operator and boundary conditions
of \cite{cerrai2005large}, specified in Assumption~\ref{ass:RD}
below. Write $H=L^2(D;\mathbb R^q)$, let $E$ be the product of
the boundary-compatible continuous-function spaces, and use the component maximum norm
$\|u\|_E=\max_i\sup_{\xi\in\overline D}|u_i(\xi)|$,
with the retained boundary conditions. Set $S(t)=e^{tA}$. We consider
\begin{equation}\label{eq:RD-system}
 du_\epsilon^x=(Au_\epsilon^x+F(u_\epsilon^x))\,dt
       +\sqrt\epsilon\,G(u_\epsilon^x)Q\,dW,
 \qquad u_\epsilon^x(0)=x\in E,
\end{equation}
where $W$ is cylindrical in $H$ and
$G(u)h(\xi)=g(\xi,u(\xi))h(\xi)$.
Here $Q$ is the \emph{noise amplitude operator}; the covariance
operator is $QQ^*$. This convention differs from the notation for
the covariance in the Burgers section.

\begin{condition}
\label{ass:RD}
We retain the elliptic setting and Hypotheses~1-5 of
\cite{cerrai2005large}, with reaction exponent $m>1$, as follows.
\begin{enumerate}[label=(\roman*)]
\item Let $d,q\geq1$, and let $D\subset\mathbb R^d$ be bounded
with smooth boundary. Let $A=\operatorname{diag}(A_1,\ldots,A_q)$,
where
\begin{equation}\label{eq:RD-elliptic-operator}
 A_i u=\sum_{h,k=1}^d\partial_{\xi_h}
       \bigl(a^i_{hk}(\xi)\partial_{\xi_k}u\bigr)-c_i u,
 \qquad i=1,\ldots,q.
\end{equation}
The coefficients $a^i_{hk}\in C^\infty(\overline D)$ form
symmetric, uniformly elliptic matrices. Each component has
homogeneous Dirichlet or co-normal boundary conditions, with
$c_i\geq0$ in the Dirichlet case and $c_i>0$ in the co-normal
case. The self-adjoint realisation of $A$ in $H$ has a complete
orthonormal eigenbasis satisfying
\[
 Ae_j=-\sigma_j e_j,\qquad \sigma_j>0,\qquad
 \sup_j\|e_j\|_\infty<\infty.
\]
\item The noise amplitude $Q$ is bounded, non-negative and
diagonal in this basis, $Qe_j=q_je_j$. Its covariance is $QQ^*$.
For $d=1$, set $\ell=\infty$ and $\chi=1/2$.
For $d\geq2$, assume
\begin{equation}\label{eq:RD-covariance-condition}
 \sum_jq_j^\ell<\infty\quad\text{for some }2\leq\ell<\infty,
 \qquad \chi:=\frac{d(\ell-2)}{2\ell}<1.
\end{equation}
\item The matrix coefficient
$g:\overline D\times\mathbb R^q\to\mathcal L(\mathbb R^q)$
is continuous and globally Lipschitz in its state variable,
uniformly in space. Consequently
\begin{equation}\label{eq:RD-g-linear}
 \sup_{\xi\in D}|g(\xi,z)|\leq C(1+|z|),\qquad z\in\mathbb R^q.
\end{equation}
\item The map $F:E\to E$ is locally Lipschitz,
$F(0)=0$, and $\|F(u)\|_E\leq C(1+\|u\|_E^m)$, for some $m> 1$.
For $x,h\in E$, suitable elements
$\delta_h\in\partial\|h\|_E$ satisfy
\begin{align}
 \langle F(x+h)-F(x),\delta_h\rangle&\leq0,
 \label{eq:RD-monotone}\\
 \langle F(x+h)-F(x),\delta_h\rangle
 &\leq-a\|h\|_E^m+C(1+\|x\|_E^m),\qquad a>0.
 \label{eq:RD-dissipative}
\end{align}
Here
\[
 \partial\|h\|_E:=\{\delta\in E^*:
       \|\delta\|_{E^*}=1,\ \langle h,\delta\rangle=\|h\|_E\}.
\]
\item Hypothesis~5 of \cite{cerrai2005large} is retained
in full: either $G(0)=0$, or the linearised range condition
(2.14) holds together with the two-sided spectral comparison
(2.15)-(2.16), the smoothness requirements, and the boundary
compatibility conditions (2.17)-(2.18).
\end{enumerate}
\end{condition}

Assumption~\ref{ass:RD} is in force throughout this section.
Hypothesis~6 of \cite{cerrai2005large} is not imposed.
Constants may depend on the fixed dimension, domain, operator
and coefficients.

We also write $A$ for its part in $E$. Under
Assumption~\ref{ass:RD}, $S(t)=e^{tA}$ is strongly continuous
and contractive on $E$. Standard elliptic semigroup estimates
give, for some $M_A\geq1$ and $\omega>0$,
\begin{equation}\label{eq:RD-semigroup-decay}
 \|S(t)\|_{\mathcal L(E)}\leq M_Ae^{-\omega t},
 \qquad t\geq0,
\end{equation}
and, for $0\leq\kappa\leq1$ and $1\leq r\leq p<\infty$,
\begin{equation}\label{eq:RD-semigroup-smoothing}
 \|S(t)h\|_{W^{\kappa,p}}
 \leq Ce^{-\omega t}(t\wedge1)^{
 -\kappa/2-\frac d2(\frac1r-\frac1p)}\|h\|_{L^r},
 \qquad t>0;
\end{equation}
see \cite[Section~2]{cerrai2005large}. Contractivity also gives
$\langle Au,\delta_u\rangle\leq0$ for $u\in D(A)$ and
$\delta_u\in\partial\|u\|_E$.

To identify the improvement, in dimension one with $Q=I$ their
additional condition reads $m>\max\{1+6\rho,2\}$ when
$|g(\xi,z)|\leq C(1+|z|^\rho)$. A cubic reaction therefore requires
$\rho<1/3$ in that argument. The result below allows linear growth
for any $m>1$ satisfying the other hypotheses above. It does not
assert that every polynomial system is monotone, or that the
controllability conditions hold for every multiplicative coefficient.

\begin{example}\label{ex:RD-cubic-linear}
On a bounded interval with zero Dirichlet boundary values, take
$A=\Delta$, $F(u)=-u^3$, and
$g(\xi,z)=\sqrt{1+z^2}$, with $Q=I$. The drift is monotone with
cubic dissipativity, the diffusion coefficient is smooth, globally
Lipschitz and linearly growing, and $g(0)=1$, $F'(0)=0$.
The spectral and boundary requirements in Hypothesis~5 hold with
$\delta=0$. Hence Theorem~\ref{thm:RD-stationary-LDP} gives the
stationary LDP for
\[
 du_\epsilon=(\Delta u_\epsilon-u_\epsilon^3)\,dt
          +\sqrt\epsilon\sqrt{1+u_\epsilon^2}\,dW.
\]
For its deterministic flow the comparison principle gives
$|z_0^x(t)|\leq S(t)|x|$, and hence
$\|z_0^x(t)\|_E\leq M_Ae^{-\omega t}\|x\|_E$.
This example has $m=3$ and linear noise growth, and therefore does
not satisfy the earlier additional inequality $m>7$.
\end{example}

\subsection{Covariance-sensitive convolution estimates}

Use the parameters $\ell$ and $\chi$ from
Assumption~\ref{ass:RD}(ii), namely
\begin{equation}\label{eq:RD-chi}
 \chi=\begin{cases}
       1/2,&d=1,\ \ell=\infty,\\
       d(\ell-2)/(2\ell),&d\geq2.
       \end{cases}
\end{equation}
The proof also works in dimension one with a finite summability
exponent and this latter value of $\chi$.

\begin{lemma}\label{lem:RD-kernel}
Under Assumption~\ref{ass:RD}, let $B$ be a bounded measurable
matrix multiplier on $D$.
For $0<t\leq1$, uniformly in $\xi\in\overline D$,
\begin{equation}\label{eq:RD-square-kernel}
 \sum_j|(S(t)BQe_j)(\xi)|^2
 \leq C\|B\|_\infty^2t^{-\chi}.
\end{equation}
In particular,
$\|S(t)BQ\|_{\mathcal L(H,E)}
 \leq C\|B\|_\infty t^{-\chi/2}$.
\end{lemma}
\begin{proof}
It suffices to estimate one output component. Set
$a_j=(S(t)Be_j)(\xi)$, including the finite component indices in
the index $j$. Uniform boundedness of the eigenfunctions and the
heat-kernel $L^1$ bound give
$\sup_j|a_j|\leq C\|B\|_\infty$.
Parseval's identity applied to the kernel row, and the heat-kernel
diagonal estimate, give
\[
 \sum_j|a_j|^2
 \leq C\|B\|_\infty^2t^{-d/2}.
\]
For $2<\ell<\infty$, H\"older's inequality and interpolation
between these two bounds imply
\begin{align*}
 \sum_jq_j^2|a_j|^2
 &\leq\left(\sum_jq_j^\ell\right)^{2/\ell}
       \left(\sum_j|a_j|^{2\ell/(\ell-2)}
       \right)^{(\ell-2)/\ell}\\
 &\leq C\|B\|_\infty^2t^{-d(\ell-2)/(2\ell)}.
\end{align*}
For $\ell=2$ use
$\sum_jq_j^2|a_j|^2\leq\sup_j|a_j|^2\sum_jq_j^2$.
For $\ell=\infty$ use boundedness of $Q$ and the Parseval bound.
Summing the finitely many output components proves the assertion.
The operator bound follows by Cauchy-Schwarz in the eigenbasis.
The multiplier has not been commuted through either $S(t)$ or $Q$.
\end{proof}

Under Assumption~\ref{ass:RD}, choose
$p>(d+2)/(1-\chi)$ and numbers $k,\alpha$ such that
\begin{equation}\label{eq:RD-factorisation-parameters}
 \frac dp<k<1-\chi-\frac2p,
 \qquad \frac k2+\frac1p<\alpha<\frac{1-\chi}{2}.
\end{equation}
This is possible because $\chi<1$. In particular
$W^{k,p}(D;\mathbb R^q)$ is compactly embedded in continuous
functions. Increasing $p$ if necessary makes $k<1$ possible.

\begin{lemma}
\label{lem:RD-damped-convolution}
Under Assumption~\ref{ass:RD}, let $p,k,\alpha$ satisfy
\eqref{eq:RD-factorisation-parameters}. Let
$I_s=[s,s+h]\subset[0,T]$, $0<h\leq1$, and let $v$ be a
predictable $E$-valued process. If $f\in\mathcal P_T^N$, define
\begin{align*}
 Z(t)={}&\int_s^t e^{-\lambda(t-r)}S(t-r)G(v(r))Qf(r)\,dr\\
 &+\sqrt\epsilon\int_s^t e^{-\lambda(t-r)}S(t-r)G(v(r))Q\,dW(r).
\end{align*}
For $\lambda\geq1$ and $0\leq\epsilon\leq1$,
\begin{equation}\label{eq:RD-damped-bound}
 \left(\mathbb E\sup_{I_s}\|Z(t)\|_E^p\right)^{1/p}
 \leq C_p\left(\sqrt{2N}\lambda^{-(1-\chi)/2}
           +\sqrt\epsilon\lambda^{-\eta}\right)
       \left[1+\left(\mathbb E\sup_{I_s}\|v(t)\|_E^p
              \right)^{1/p}\right],
\end{equation}
where $\eta=(1-\chi-k-2/p)/2>0$. Constants do not depend on
$s,h,T,N,\epsilon,\lambda$, or $f$.
\end{lemma}
\begin{proof}
Lemma~\ref{lem:RD-kernel}, linear growth and Cauchy-Schwarz give,
pathwise, for the deterministic convolution,
\begin{align*}
 &\sup_{I_s}\left\|\int_s^t e^{-\lambda(t-r)}S(t-r)
                  G(v(r))Qf(r)\,dr\right\|_E\\
 &\quad\leq C(1+\sup_{I_s}\|v\|_E)
       \left(\int_0^1e^{-2\lambda r}r^{-\chi}\,dr\right)^{1/2}
       \|f\|_{L^2(I_s;H)}\\
 &\quad\leq C\sqrt{2N}\lambda^{-(1-\chi)/2}
                     (1+\sup_{I_s}\|v\|_E).
\end{align*}
For the stochastic part, factorise it using
\[
 Y_\alpha(r)=\int_s^r(r-a)^{-\alpha}e^{-\lambda(r-a)}
                 S(r-a)G(v(a))Q\,dW(a).
\]
The outer convolution maps $L^p(I_s;L^p)$ to
$C(I_s;W^{k,p})$, with norm bounded by
$C\lambda^{-\alpha+k/2+1/p}$. Indeed its H\"older kernel
has exponent $(\alpha-1-k/2)p/(p-1)>-1$.
By pointwise BDG, Lemma~\ref{lem:RD-kernel}, and Fubini,
\begin{align*}
 \mathbb E\|Y_\alpha\|_{L^p(I_s;L^p)}^p
 &\leq C_p\left(\int_0^1r^{-2\alpha-\chi}
                           e^{-2\lambda r}\,dr\right)^{p/2}
          \mathbb E(1+\sup_{I_s}\|v\|_E)^p\\
 &\leq C_p\lambda^{(2\alpha+\chi-1)p/2}
          \mathbb E(1+\sup_{I_s}\|v\|_E)^p.
\end{align*}
Here $2\alpha+\chi<1$. Multiplying the two powers of $\lambda$
gives $\lambda^{-\eta}$. The embedding
$W^{k,p}\hookrightarrow C(\overline D;\mathbb R^q)$, with the
retained boundary conditions, and the deterministic estimate
prove \eqref{eq:RD-damped-bound}.
\end{proof}

\subsection{Uniform controlled moments and exponential estimates}

Let $u_{\epsilon,f}^x$ solve the controlled equation obtained by
adding $G(u_{\epsilon,f}^x)Qf$ to the drift in
\eqref{eq:RD-system}. We include $\epsilon=0$ and write
$z_f^x:=u_{0,f}^x$. For deterministic $f\in L^2(0,T;H)$,
this is the unique mild solution of
\begin{equation}\label{eq:RD-skeleton}
 \partial_t z_f^x=Az_f^x+F(z_f^x)+G(z_f^x)Qf,
 \qquad z_f^x(0)=x,
\end{equation}
with the boundary conditions encoded by $A$. Thus, the noise
coefficient $\Sigma$ in \eqref{eq:general-skeleton-mild} is
$\Sigma(z)=G(z)Q$ here.

\begin{lemma}\label{lem:RD-uniform-controlled}
Under Assumption~\ref{ass:RD}, for every $N>0$ and sufficiently
large $p$ there is $C_{N,p}$ such that
\begin{equation}\label{eq:RD-uniform-controlled}
 \sup_{T>0}\sup_{0\leq\epsilon\leq1}\sup_{f\in\mathcal P_T^N}
 \mathbb E\|u_{\epsilon,f}^x(T)\|_E^p
 \leq C_{N,p}(1+\|x\|_E^p).
\end{equation}
For every finite $T$, the corresponding
$\mathbb E\sup_{0\leq t\leq T}\|u_{\epsilon,f}^x(t)\|_E^p$
is bounded uniformly for $\epsilon\in[0,1]$, bounded $x$, and
$f\in\mathcal P_T^N$.
\end{lemma}
\begin{proof}
Restart the damped convolution of
Lemma~\ref{lem:RD-damped-convolution} on $I_s=[s,s+1]$, with
$v=u_{\epsilon,f}^x$, and write $w=u_{\epsilon,f}^x-Z$.
The remaining equation is
\[
 \partial_tw=Aw+F(w+Z)+\lambda Z,\qquad
 w(s)=u_{\epsilon,f}^x(s).
\]
The sup-norm subdifferential inequality, polynomial growth and
\eqref{eq:RD-dissipative} imply
\begin{align*}
 d^-\|w(t)\|_E
 &\leq-a\|w(t)\|_E^m+C(1+\|Z(t)\|_E^m)+\lambda\|Z(t)\|_E\\
 &\leq-a\|w(t)\|_E^m
     +C(1+\lambda^{1/(m-1)}+\|Z(t)\|_E)^m.
\end{align*}
The second line uses Young's inequality with conjugate exponents
$m$ and $m/(m-1)$. The comparison and convolution estimates can
first be applied up to the exit time from an $E$-ball. We use
this localisation below to establish finite moments before
applying the unstopped endpoint contraction.

For fixed $t\in I_s$ set
$b=1+\lambda^{1/(m-1)}+\sup_{s\leq r\leq t}\|Z(r)\|_E$.
Comparison with $y'=-ay^m+Cb^m$ gives
\begin{equation}\label{eq:RD-scalar-comparison}
 \|w(t)\|_E\leq Cb+
 \frac{\|u_{\epsilon,f}^x(s)\|_E}
 {(1+c(t-s)\|u_{\epsilon,f}^x(s)\|_E^{m-1})^{1/(m-1)}}.
\end{equation}
For example, a sufficiently large constant multiple of $b$ plus
the solution of $y'=-cy^m$ is a supersolution; this proves the
display also when the initial value is arbitrarily large.
Consequently,
\begin{align*}
 \sup_{I_s}\|u_{\epsilon,f}^x\|_E
 &\leq\|u_{\epsilon,f}^x(s)\|_E
          +C(1+\lambda^{1/(m-1)})+C\sup_{I_s}\|Z\|_E,\\
 \|u_{\epsilon,f}^x(s+1)\|_E
 &\leq C(1+\lambda^{1/(m-1)})+C\sup_{I_s}\|Z\|_E.
\end{align*}
In the second line we used $m>1$ to bound the last term of
\eqref{eq:RD-scalar-comparison} at $t-s=1$, independently of
the initial value.

Put $a_s=(\mathbb E\|u_{\epsilon,f}^x(s)\|_E^p)^{1/p}$ and
$b_s=(\mathbb E\sup_{I_s}\|u_{\epsilon,f}^x\|_E^p)^{1/p}$.
Choose $\lambda=\lambda(N,p)$ so large that
\[
 CC_p\bigl(\sqrt{2N}\lambda^{-(1-\chi)/2}
                              +\lambda^{-\eta}\bigr)\leq\frac14.
\]
Put $D=C(1+\lambda^{1/(m-1)})$.
To justify absorption without assuming that $b_s$ is already
finite, first stop at
$\tau_K=\inf\{t\geq s:\|u_{\epsilon,f}^x(t)\|_E\geq K\}
\wedge(s+1)$ and stop both integrands defining $Z$ there.
On $[s,\tau_K]$ the same comparison applies, and the stopped
convolution has the coefficient in
\eqref{eq:RD-damped-bound}, with the stopped path supremum
on its right. That supremum has a finite $p$th moment when
$a_s<\infty$, since it is bounded by $K+\|u_{\epsilon,f}^x(s)\|_E$.
Absorption therefore gives
\[
 \left(\mathbb E\sup_{s\leq t\leq\tau_K}
                  \|u_{\epsilon,f}^x(t)\|_E^p\right)^{1/p}
 \leq\frac43a_s+\frac43D+\frac13,
\]
uniformly in $K$. Letting $K\to\infty$ and using Fatou's lemma
proves that $b_s<\infty$ with this bound. Starting from
deterministic $x$, a finite iteration establishes these
moments on every finite interval. We can now apply the
unstopped endpoint comparison and take its $L^p(\Omega)$ norm.
The two estimates are
\begin{align}
 b_s&\leq\frac43a_s+\frac43D+\frac13,
 \label{eq:RD-unit-sup-bound}\\
 a_{s+1}&\leq D+\frac14(1+b_s)
          \leq\frac13a_s+\frac43D+\frac13.
 \label{eq:RD-endpoint-contraction}
\end{align}
All constants are independent of $s,T,\epsilon,f$.
At integer times iteration gives
\[
 a_j\leq3^{-j}\|x\|_E+
     \left(\frac43D+\frac13\right)\sum_{i=0}^{j-1}3^{-i}
 \leq\|x\|_E+2D+\frac12.
\]
On the final interval $[\lfloor T\rfloor,T]$, the same absorption
as in \eqref{eq:RD-unit-sup-bound} applies. This proves
\eqref{eq:RD-uniform-controlled}. Summing the $p$th moments
of the unit-interval suprema proves the assertion on every fixed
finite interval. The preliminary localisation above justifies
all moment absorptions; the endpoint contraction is applied
only after these finite moments have been established.
\end{proof}

\begin{remark}\label{rem:RD-deterministic-bound}
Under Assumption~\ref{ass:RD}, for $\epsilon=0$ and deterministic
controls the same calculation
is pathwise. In particular, controls of total energy at most $N$
satisfy
\begin{equation}\label{eq:RD-global-skeleton-bound}
 \sup_{0\leq t\leq T}\|z_f^x(t)\|_E\leq C_N(1+\|x\|_E),
 \qquad T>0.
\end{equation}
This is a bounded-initial-set estimate, not a uniform estimate
over every $x\in E$. The latter is unnecessary below.
In dimension one with bounded $Q$, the deterministic choice
$\lambda=C(1+N)^2$ is sufficient. For a cubic drift this gives
$\sup_t\|z_f^x(t)\|_E\leq C(1+\|x\|_E+N)$.
\end{remark}

\begin{proposition}\label{prop:RD-invariant-norm-tails}
Under Assumption~\ref{ass:RD}, for every $r>0$ there is $R_r$
such that
\[
 \sup_{T>0}\mathbb P(\|u_\epsilon^0(T)\|_E\geq R_r)
       \leq e^{-r/\epsilon},\qquad0<\epsilon\leq1.
\]
Every invariant measure $\nu_\epsilon$ obtained from the
Krylov-Bogoliubov averages starting at $0$ satisfies the same
estimate outside an $E$-ball, after enlarging the radius.
\end{proposition}
\begin{proof}
Lemma~\ref{lem:RD-uniform-controlled}, including its deterministic
case, gives
\[
 \sup_{T>0}\sup_{0\leq\epsilon\leq1}
 \sup_{f\in\mathcal P_T^N}
       \mathbb E\|u_{\epsilon,f}^0(T)\|_E\leq K_N<\infty.
\]
This is \eqref{eq:terminal-controlled-bound}, with
$x_0=0$ and $\epsilon_*=1$. Proposition~\ref{prop:terminal-criterion}
proves both assertions, uniformly over the selected
Krylov-Bogoliubov limit.

For completeness, ordinary spatial moments justify existence of
these limits under the present assumptions. Choose $k,p,\alpha$
as in \eqref{eq:RD-factorisation-parameters}. The endpoint and
unit interval bounds, with zero control and moment order $mp$,
give
\[
 \sup_{s\geq0}\sup_{0<\epsilon\leq1}
 \mathbb E\sup_{s\leq t\leq s+1}\|u_\epsilon^0(t)\|_E^{mp}<\infty.
\]
In the mild equation restarted at $s$, elliptic smoothing bounds
the initial and drift terms in $W^{k,p}$ at time $s+1$.
The kernel bound \eqref{eq:RD-square-kernel} and the same
factorisation calculation as in
Lemma~\ref{lem:RD-damped-convolution}, now keeping the
$W^{k,p}$ norm, bound the stochastic term. Consequently
\[
 \sup_{s\geq0}\sup_{0<\epsilon\leq1}
       \mathbb E\|u_\epsilon^0(s+1)\|_{W^{k,p}}^p<\infty.
\]
The compact embedding into $E$ gives tightness of the time
averages. Their invariant limits follow from the Feller property
and the Krylov-Bogoliubov argument, as in
\cite[Theorem~6.2]{cerrai2003stochastic}. This ordinary moment
estimate is used for existence; the LDP uses only the exponential
norm estimate above.
\end{proof}

\subsection{Verification of the stationary LDP}

For the skeleton \eqref{eq:RD-skeleton}, let
\[
 I_T^x(v):=\frac12\inf\left\{
       \|f\|_{L^2(0,T;H)}^2:
       f\in L^2(0,T;H),\ v=z_f^x\right\},
 \qquad v\in C([0,T];E),
\]
with $\inf\varnothing=+\infty$. Define the quasipotential by
\begin{equation}\label{eq:RD-quasipotential}
 V_{\rm RD}(y):=
 \inf_{\substack{T>0,\ v\in C([0,T];E)\\
                         v(0)=0,\ v(T)=y}}I_T^0(v),
 \qquad y\in E.
\end{equation}

\begin{theorem}[Reaction-diffusion stationary LDP]
\label{thm:RD-stationary-LDP}
Under Assumption~\ref{ass:RD}, every family of invariant measures
$\{\nu_\epsilon\}_{0<\epsilon\leq1}$ selected as subsequential
limits of the Krylov-Bogoliubov averages starting from $0$
satisfies an LDP in $E$ as $\epsilon\downarrow0$, with speed
$1/\epsilon$ and good rate function $V_{\rm RD}$ defined in
\eqref{eq:RD-quasipotential}. No condition corresponding to
Hypothesis~6 of \cite{cerrai2005large} is required.
\end{theorem}

\begin{proof}
We first consider the case $G(0)=0$ in Hypothesis~5 of
\cite{cerrai2005large}. Since $F(0)=0$, uniqueness gives
$u_\epsilon^0(t)=0$ for all $t\geq0$. Every Krylov-Bogoliubov
average starting from $0$ is therefore $\delta_0$.
The controlled equation starting from $0$ also remains at $0$,
so
\[
 V_{\rm RD}(y)=
 \begin{cases}
  0,&y=0,\\
  +\infty,&y\ne0.
 \end{cases}
\]
The asserted LDP follows immediately. Henceforth we assume
the range condition (2.14) and the additional spectral,
smoothness and boundary conditions in the other branch of
Hypothesis~5, all retained in Assumption~\ref{ass:RD}.

For each fixed initial value, compactness of the finite time
action sublevels follows from
\cite[Theorem~5.1]{cerrai2004large}.
Estimate \eqref{eq:RD-global-skeleton-bound} bounds the
skeletons uniformly over bounded initial values and controls
of bounded energy. Local Lipschitz continuity of $F$, global
Lipschitz continuity of $g$, and the kernel estimate
\eqref{eq:RD-square-kernel} give
\begin{equation}\label{eq:RD-skeleton-continuity}
 \|z_f^x-z_h^y\|_{C([0,T];E)}
 \leq C_{T,R,N}\left(\|x-y\|_E+
                       \|f-h\|_{L^2(0,T;H)}\right)
\end{equation}
for $\|x\|_E,\|y\|_E\leq R$ and $f,h\in\mathcal S_T^N$.
Indeed, on an interval of length $t\leq1$, the difference
of the controlled convolutions corresponding to the same
control $f$ is bounded by
\[
 Ct^{(1-\chi)/2}\|f\|_{L^2}
       \|z_f^x-z_h^y\|_{C([0,t];E)}.
\]
Choose $t$ sufficiently small, depending on $N$, that its
coefficient is at most $1/2$. The contribution involving
$f-h$ is bounded by the same kernel estimate. Absorbing the
path difference and applying Gronwall's inequality to the
drift term proves the estimate on this interval. Iteration
over finitely many such intervals gives
\eqref{eq:RD-skeleton-continuity} for each fixed $T$.

Estimate \eqref{eq:RD-global-skeleton-bound} also verifies
\eqref{eq:det-uniform-bound}, with a constant independent
of the duration. To verify
\eqref{eq:det-positive-time-compactness}, fix $t_0>0$ and
choose $k,p$ as in \eqref{eq:RD-factorisation-parameters}.
In particular, $k>d/p$ and $k+\chi<1$.
For $\|x\|_E\leq R$ and $f\in\mathcal S_{t_0}^N$, write
$M=M_{R,N}$ for the uniform path bound. Splitting the
semigroup into two factors and using
Lemma~\ref{lem:RD-kernel} gives
\[
 \|S(t)G(z)Qh\|_{W^{k,p}}
 \leq C_{t_0}t^{-(k+\chi)/2}
                  (1+\|z\|_E)\|h\|_H,
 \qquad0<t\leq t_0.
\]
Consequently, Cauchy-Schwarz yields
\begin{align*}
 &\left\|\int_0^{t_0}S(t_0-s)
                 G(z_f^x(s))Qf(s)\,ds\right\|_{W^{k,p}}\\
 &\qquad\leq C_{t_0}(1+M)
       \left(\int_0^{t_0}t^{-k-\chi}\,dt\right)^{1/2}
       \|f\|_{L^2(0,t_0;H)}
 \leq C_{t_0}(1+M)\sqrt N.
\end{align*}
The initial and drift terms satisfy
\[
 \|S(t_0)x\|_{W^{k,p}}\leq C_{t_0}R,
 \qquad
 \left\|\int_0^{t_0}S(t_0-s)F(z_f^x(s))\,ds\right\|_{W^{k,p}}
 \leq C_{t_0}(1+M^m),
\]
where the second bound uses the integrability of $t^{-k/2}$.
Thus
\[
 \sup_{\|x\|_E\leq R,\ f\in\mathcal S_{t_0}^N}
       \|z_f^x(t_0)\|_{W^{k,p}}
 \leq C_{t_0}\left[R+1+M^m+(1+M)\sqrt N\right].
\]
Since $kp>d$, the Sobolev embedding makes these endpoints
relatively compact in $C(\overline D;\mathbb R^q)$.
The Dirichlet boundary values are preserved under uniform
limits; for co-normal components the continuous state space
is $C(\overline D)$. Thus the closure lies in $E$, proving
\eqref{eq:det-positive-time-compactness}.

Under the retained operator and reaction hypotheses,
\cite[Proposition~3.6]{cerrai2005large} gives
\[
 \lim_{t\to\infty}\sup_{\|x\|_E\leq R}
       \|z_0^x(t)\|_E=0,
 \qquad R>0.
\]
In particular, attraction holds uniformly on compact initial
sets. All parts of Assumption~\ref{ass:stationary-input}
are therefore verified.

We next verify Assumption~\ref{ass:finite-path-approximation}
directly. Let $\delta$ be the exponent in the spectral
comparison (2.15) of \cite{cerrai2005large}, and put
\[
 \mathsf Y:=D((-A)^{\delta+1/2}).
\]
We first record the small cost connection supplied by the
construction in \cite[Theorem~4.2]{cerrai2005large}.
For a sufficiently small fixed $\tau>0$, the endpoint map
of the linearised equation is
\[
 \mathcal L_\tau f
 :=\int_0^\tau e^{(\tau-s)(A+DF(0))}G(0)Qf(s)\,ds.
\]
The retained range and spectral conditions give a bounded
linear operator
$R_\tau:\mathsf Y\to L^2(0,\tau;H)$ such that
$\mathcal L_\tau R_\tau=I_{\mathsf Y}$.
Under the retained smoothness and boundary assumptions, the
nonlinear endpoint map
\[
 \Gamma_\tau(p):=z_{R_\tau p}^0(\tau)
\]
is continuously differentiable near $0$ in $\mathsf Y$, with
\[
 \Gamma_\tau(0)=0,
 \qquad D\Gamma_\tau(0)=\mathcal L_\tau R_\tau=I_{\mathsf Y}.
\]
The inverse function theorem gives a local inverse satisfying
$\|\Gamma_\tau^{-1}(y)\|_{\mathsf Y}\leq2\|y\|_{\mathsf Y}$
for sufficiently small $y$. Hence, for some $\rho,C_\tau>0$,
the controls
$f_y:=R_\tau\Gamma_\tau^{-1}(y)$ satisfy
\begin{equation}\label{eq:RD-small-cost-connection}
 z_{f_y}^0(\tau)=y,
 \qquad
 \frac12\|f_y\|_{L^2(0,\tau;H)}^2
 \leq C_\tau\|y\|_{\mathsf Y}^2,
 \qquad \|y\|_{\mathsf Y}<\rho.
\end{equation}

Now let $z\in C((-\infty,0];E)$ satisfy
$z(t)\to0$ in $E$ as $t\to-\infty$ and
$I_{-\infty}(z)<\infty$. Its finite action representation
provides a control $f\in L^2(-\infty,0;H)$ generating $z$.
Exponential decay of $S(t)$ allows the initial term in the
mild equation to tend to zero as the initial time tends to
$-\infty$. Thus $z$ satisfies the infinite past mild equation
required in \cite[Lemma~3.5]{cerrai2005large}. That lemma gives
\[
 \lim_{t\to-\infty}\|z(t)\|_{\mathsf Y}=0.
\]
Fix $\eta>0$. Choose $S>0$ sufficiently large that
\[
 \|z(-S)\|_{\mathsf Y}<\rho,
 \qquad C_\tau\|z(-S)\|_{\mathsf Y}^2\leq\eta.
\]
By \eqref{eq:RD-small-cost-connection}, the path
$w:=z_{f_{z(-S)}}^0$ on $[0,\tau]$ satisfies
\[
 w(0)=0,\qquad w(\tau)=z(-S),\qquad I_\tau^0(w)\leq\eta.
\]
Concatenate this path with the restriction of $z$ to $[-S,0]$:
\[
 v_\eta(t):=
 \begin{cases}
  w(t),&0\leq t\leq\tau,\\
  z(t-\tau-S),&\tau\leq t\leq\tau+S.
 \end{cases}
\]
Then $v_\eta(0)=0$, $v_\eta(\tau+S)=z(0)$ and, by
additivity of the action,
\[
 I_{\tau+S}^0(v_\eta)
 \leq \eta+I_{-S,0}(z)
 \leq \eta+I_{-\infty}(z).
\]
Translating $v_\eta$ to $[-(\tau+S),0]$ verifies
\eqref{eq:finite-path-approximation}. In particular, this
verification uses the regularity and local control results,
without assuming the quasipotential identity.

Proposition~\ref{prop:deterministic-compactness} now identifies
$V_{\rm RD}$ with the infinite past quasipotential and proves
compactness of its sublevels. It also supplies the end point estimate \eqref{eq:stationary-endpoint-approximation}.

Finally, \cite[Theorems~6.2 and~6.3]{cerrai2004large} give
the finite time Freidlin-Wentzell uniform LDP on bounded
initial sets, with the action $I_T^x$ defined above.
Proposition~\ref{prop:RD-invariant-norm-tails} gives the
stationary norm estimate \eqref{eq:stationary-norm-tail},
uniformly over the selected Krylov-Bogoliubov limits.
Theorem~\ref{thm:stationary-transfer} therefore proves the
LDP with good rate function $V_{\rm RD}$ for every such
selection. The noise amplitude in \cite{cerrai2005large} is
written as $\epsilon$; replacing it by $\sqrt\epsilon$ gives
speed $1/\epsilon$ as used here.
\end{proof}

\begin{remark}
The proof uses the finite time uniform LDP on bounded initial
sets, stationary norm tails and deterministic regularisation.
The estimate uniform over all initial values in
\cite[Theorem~3.4]{cerrai2005large} is unnecessary, so the
additional growth restriction used to obtain that estimate
can be removed. The spectral, regularity and local control
requirements of Hypothesis~5 remain in force. The
quasipotential identity and compactness follow from the
abstract proposition after the approximation above has
been verified.
\end{remark}

\section{Allen-Cahn with inhomogeneous boundary conditions}
\label{sec:AC-application}

In \cite{bai2025allencahn}, we proved the stationary LDP for
the Allen-Cahn equation with opposite boundary values and
multiplicative noise satisfying
$|g(\xi,r)|\leq C(1+|r|^\rho)$, $0\leq\rho<1$.
We now remove this additional growth restriction and assume
global Lipschitz continuity in the state variable. Since the
reaction is not globally monotone, we use the deterministic
energy analysis from that paper to verify the abstract
hypotheses. The new input is a controlled moment estimate
uniform in time, which allows us to apply the abstract results
without proving stochastic exponential tightness.

\subsection{The model and its controlled equation}

Let $L>0$, $H=L^2(-L,L)$, and $\psi(\xi)=\xi/L$. Set
\[
 E=\{v\in C([-L,L]):v(-L)=v(L)=0\},
 \qquad \mathcal E=E+\psi,
\]
with the supremum norm on $E$ and the corresponding metric
on $\mathcal E$. We consider
\begin{equation}\label{ACE1}
\begin{cases}
 \partial_tu_\epsilon
 =\Delta u_\epsilon+u_\epsilon-u_\epsilon^3
       +\sqrt\epsilon\,g(\xi,u_\epsilon)\partial_tW,
       &\xi\in(-L,L),\\
 u_\epsilon(0)=x+\psi,\quad x\in E,\\
 u_\epsilon(t,\pm L)=\pm1,&t\geq0,
\end{cases}
\end{equation}
where $W$ is a cylindrical Wiener process on $H$.
Throughout this section,
$g:[-L,L]\times\mathbb R\to\mathbb R$ is continuous and
globally Lipschitz in its second variable, uniformly in $\xi$.
Consequently,
\begin{equation}\label{eq:linear-growth-g}
 |g(\xi,r)|\leq C_g(1+|r|),
 \qquad (\xi,r)\in[-L,L]\times\mathbb R.
\end{equation}
For the deterministic connection argument, we also assume
\begin{condition}\label{assumption 3}
There is $g_0>0$ such that
\begin{equation}\label{inf g}
 \inf_{(\xi,r)\in[-L,L]\times\mathbb R}|g(\xi,r)|=g_0.
\end{equation}
\end{condition}

Let $A=\Delta$ with homogeneous Dirichlet boundary conditions
and $S(t)=e^{tA}$. The shifted solution
$\bar u_\epsilon^x=u_\epsilon-\psi$ satisfies
\begin{equation}\label{ACE2}
 d\bar u_\epsilon^x
 =[A\bar u_\epsilon^x+F(\bar u_\epsilon^x)]\,dt
       +\sqrt\epsilon\,G(\bar u_\epsilon^x)\,dW,
 \qquad \bar u_\epsilon^x(0)=x,
\end{equation}
where
\[
 F(z)=(z+\psi)-(z+\psi)^3,\qquad
 [G(z)h](\xi)=g(\xi,z(\xi)+\psi(\xi))h(\xi).
\]
The results of \cite[Theorems~5.3, 5.5 and~6.2]
{cerrai2003stochastic} give a unique global mild solution,
a Feller Markov family, and invariant measures obtained
through Krylov-Bogoliubov averages. For each $\epsilon>0$,
let $\bar\mu_\epsilon$ be any weak limit of
\[
 \bar\mu_{\epsilon,T}
 =\frac1T\int_0^T\mathcal L(\bar u_\epsilon^0(t))\,dt
\]
along a sequence $T\to\infty$. Uniqueness of the invariant
measure is not needed. Define the corresponding invariant
measure on $\mathcal E$ by
\begin{equation}\label{mu_epsilon}
 \mu_\epsilon(B+\psi)=\bar\mu_\epsilon(B),
 \qquad B\in\mathcal B(E).
\end{equation}

For $f\in\mathcal P_T^N$, let $\bar u_{\epsilon,f}^x$ solve
\begin{equation}\label{eqn: u epsilon f}
\begin{aligned}
 d\bar u_{\epsilon,f}^x
 &=\bigl[A\bar u_{\epsilon,f}^x
       +F(\bar u_{\epsilon,f}^x)
       +G(\bar u_{\epsilon,f}^x)f\bigr]\,dt\\
 &\quad+\sqrt\epsilon\,G(\bar u_{\epsilon,f}^x)\,dW,
 \qquad \bar u_{\epsilon,f}^x(0)=x.
\end{aligned}
\end{equation}
The control classes are defined with noise Hilbert space
$H$. At $\epsilon=0$, write $z_f^x=\bar u_{0,f}^x$; thus
\begin{equation}\label{Skeleton equation-2}
 \partial_tz_f^x=Az_f^x+F(z_f^x)+G(z_f^x)f,
 \qquad z_f^x(0)=x.
\end{equation}
For deterministic $f\in L^2(0,T;H)$, this equation has a
unique mild solution in $C([0,T];E)$ by
\cite[Theorem~4.1]{cerrai2004large}.

We use the energy and stationary profile studied in
\cite[Section~3.1]{bai2025allencahn}:
\begin{equation}\label{Energy}
 \mathbf E_L(u)
 =\int_{-L}^L\left[\frac12|u'(\xi)|^2
                  +\frac14(u(\xi)^2-1)^2\right]\,d\xi,
 \qquad u\in\mathcal E\cap H^1(-L,L),
\end{equation}
with value $+\infty$ otherwise. Its unique minimiser $m_L$
is the unique stationary profile with $m_L(\pm L)=\pm1$.
Put
\begin{equation}\label{Energy 0}
 u_\ast=m_L-\psi,\qquad
 \mathbf E_L^*(v)=\mathbf E_L(v+\psi)-\mathbf E_L(m_L).
\end{equation}
Then $Au_\ast+F(u_\ast)=0$. The action and quasipotential are
\begin{equation}\label{Rate functional}
\begin{aligned}
 I_T^x(z)
 &=\frac12\inf\left\{\int_0^T\|f(t)\|_H^2\,dt:
             f\in L^2(0,T;H),\ z=z_f^x\right\},\\
 U_L(y)
 &=\inf\{I_T^{u_\ast}(z):T>0,\ z(T)=y\},
 \qquad y\in E,
\end{aligned}
\end{equation}
where the paths belong to $C([0,T];E)$ and
$\inf\varnothing=+\infty$. Constants below may depend on
$L$ and the fixed coefficients.

\subsection{Uniform controlled moments and exponential estimates}
\label{subsec:AC-exponential-estimates}
\label{subsec:AC-linear-growth}

The following estimate is the new stochastic input.
Nondegeneracy is not used in this subsection.

\begin{lemma}\label{lem:AC-linear-uniform-controlled}
For every $p>6$, there is $C_p>0$ such that, for all
$N,T>0$ and $x\in E$,
\begin{equation}\label{eq:AC-controlled-terminal-moment}
 \sup_{0\leq\epsilon\leq1}\sup_{f\in\mathcal P_T^N}
 \left(\mathbb E\|\bar u_{\epsilon,f}^x(T)\|_E^p\right)^{1/p}
 \leq C_p(1+\|x\|_E+N).
\end{equation}
On a finite interval,
\begin{equation}\label{inequality: expectation}
 \sup_{0\leq\epsilon\leq1}\sup_{f\in\mathcal P_T^N}
 \mathbb E\sup_{0\leq t\leq T}
                 \|\bar u_{\epsilon,f}^x(t)\|_E^p
 \leq C_p(1+T)(1+\|x\|_E+N)^p.
\end{equation}
For deterministic controls, one also has
\begin{equation}\label{ineqn: skeleton}
 \sup_{f\in\mathcal S_T^N}\sup_{0\leq t\leq T}
       \|z_f^x(t)\|_E\leq C(1+\|x\|_E+N),
\end{equation}
with a constant independent of $T$.
\end{lemma}

\begin{proof}
Write $u=\bar u_{\epsilon,f}^x$. On an interval $[s,s+1]$,
restart the damped convolution by setting
\[
\begin{aligned}
 Z(t)
 &=\int_s^t e^{-\lambda(t-r)}S(t-r)G(u(r))f(r)\,dr\\
 &\quad+\sqrt\epsilon\int_s^t
          e^{-\lambda(t-r)}S(t-r)G(u(r))\,dW(r),
 \qquad \lambda\geq1.
\end{aligned}
\]
Controls may be extended by zero beyond $T$.
Choose $3/(2p)<\vartheta<1/4$. The convolution calculation
in Lemma~\ref{lem:RD-damped-convolution}, with $d=q=1$,
$Q=I$, $\chi=1/2$ and $k=2\vartheta-2/p>1/p$, gives
\begin{align}
 &\left(\mathbb E\sup_{s\leq t\leq s+1}
                          \|Z(t)\|_E^p\right)^{1/p}
 \nonumber\\
 &\quad\leq C_{p,\vartheta}
 \left(\sqrt N\lambda^{-1/4}
       +\sqrt\epsilon\lambda^{-1/4+\vartheta}\right)
 \left[1+\left(\mathbb E\sup_{s\leq t\leq s+1}
                           \|u(t)\|_E^p\right)^{1/p}\right].
 \label{eq:AC-Z-linear}
\end{align}
That calculation uses only the heat semigroup and linear
growth of the multiplier, so applies to the shifted
coefficient $G$ here. Its factorisation parameters satisfy
$k/2+1/p=\vartheta<\alpha<1/4$; the resulting damping power is
\[
 \lambda^{-\alpha+k/2+1/p}\lambda^{\alpha-1/4}
 =\lambda^{-1/4+\vartheta}.
\]

The remainder $v=u-Z$ solves
\[
 \partial_tv=Av+F(v+Z)+\lambda Z,
 \qquad v(s)=u(s).
\]
For real $a,b$,
\[
 \operatorname{sgn}(a)\bigl[(a+b)-(a+b)^3\bigr]
 \leq-c|a|^3+C(1+|b|^3).
\]
The maximum principle for $A$, applied at points where
$|v|$ attains its maximum, therefore gives
\[
 D^+\|v(t)\|_E
 \leq-c\|v(t)\|_E^3+C(1+\|Z(t)\|_E^3)
                      +\lambda\|Z(t)\|_E.
\]
Here $\|\psi\|_\infty=1$ is included in the constant.
The argument for mild solutions follows by approximation,
as in \cite[Proposition~6.1]{cerrai2003stochastic}.
Young's inequality bounds the last term by
$\frac13\|Z(t)\|_E^3+\frac23\lambda^{3/2}$.
Comparison with the scalar differential equation yields
\begin{align}
 \|v(t)\|_E
 &\leq C\left(1+\sqrt\lambda+
                  \|Z\|_{C([s,s+1];E)}\right)
 \nonumber\\
 &\quad+\frac{\|u(s)\|_E}
        {\sqrt{1+c(t-s)\|u(s)\|_E^2}},
 \qquad s<t\leq s+1.
 \label{eq:AC-v-comparison}
\end{align}
Indeed, a sufficiently large multiple of
$1+\sqrt\lambda+\|Z\|_{C([s,s+1];E)}$, plus the displayed
decreasing term, is a supersolution.
Consequently, for one constant $C_\ast\geq1$,
\begin{align*}
 \sup_{s\leq t\leq s+1}\|u(t)\|_E
 &\leq\|u(s)\|_E+C_\ast(1+\sqrt\lambda)
                     +C_\ast\|Z\|_{C([s,s+1];E)},\\
 \|u(s+1)\|_E
 &\leq C_\ast(1+\sqrt\lambda)
                     +C_\ast\|Z\|_{C([s,s+1];E)}.
\end{align*}
The second inequality uses
$r/\sqrt{1+cr^2}\leq c^{-1/2}$.

Put $K_p=C_\ast C_{p,\vartheta}$ and choose
\begin{equation}\label{eq:choose-lambda-AC}
 \lambda=\max\left\{1,(8K_p\sqrt N)^4,
                (8K_p)^{1/(1/4-\vartheta)}\right\}.
\end{equation}
Then $\sqrt\lambda\leq C_p(1+N)$ and, uniformly for
$0\leq\epsilon\leq1$,
\[
 K_p\left(\sqrt N\lambda^{-1/4}
              +\sqrt\epsilon\lambda^{-1/4+\vartheta}\right)
 \leq\frac14.
\]
To justify absorption before assuming finite moments, stop
at the first exit of $u$ from an $E$-ball of radius $K$,
capped at $s+1$, and stop both convolution integrands.
Up to that time the same estimates hold, and the stopped
supremum has a finite $p$th moment whenever $u(s)$ does.
Absorption gives a bound independent of $K$; letting
$K\to\infty$ and using Fatou's lemma gives finite moments
on the whole interval. Iteration from the deterministic
initial value justifies the following unstopped estimates
on every finite horizon.

Write
\[
 a_s=\left(\mathbb E\|u(s)\|_E^p\right)^{1/p},
 \qquad
 b_s=\left(\mathbb E\sup_{s\leq t\leq s+1}
                         \|u(t)\|_E^p\right)^{1/p},
 \qquad D=C_\ast(1+\sqrt\lambda).
\]
The supremum estimate and \eqref{eq:AC-Z-linear} give
$b_s\leq a_s+D+\frac14(1+b_s)$, and hence
\begin{equation}\label{eq:AC-unit-sup}
 b_s\leq\frac43a_s+\frac43D+\frac13.
\end{equation}
The endpoint estimate is stronger:
\[
 a_{s+1}\leq D+\frac14(1+b_s)
          \leq\frac13a_s+\frac43D+\frac13.
\]
At integer times this yields
\[
 a_j\leq3^{-j}\|x\|_E+2D+\frac12,
 \qquad
 b_j\leq\frac43\,3^{-j}\|x\|_E+4D+1.
\]
The same convolution and supremum estimates hold on
intervals of length at most one, so they also control
the final interval $[\lfloor T\rfloor,T]$.
Since $D\leq C_p(1+N)$, this proves
\eqref{eq:AC-controlled-terminal-moment}.
Summing the $p$th moments of the interval suprema proves
\eqref{inequality: expectation}.
For deterministic controls and $\epsilon=0$, the same
interval bounds are pathwise and uniform in $j$;
taking their maximum proves \eqref{ineqn: skeleton}.
\end{proof}

\begin{theorem}\label{thm:AC-linear-invariant-exp}
\label{prop:exponential estimate}
For every $r>0$, there are $\widetilde R_r,R_r>0$ such that
\begin{align}
 \sup_{T>0}\mathbb P\bigl(
          \|\bar u_\epsilon^0(T)\|_E\geq\widetilde R_r\bigr)
 &\leq e^{-r/\epsilon},
 \label{eq:AC-uniform-time-exp}\\
 \bar\mu_\epsilon\bigl(\{x:\|x\|_E\geq R_r\}\bigr)
 &\leq e^{-r/\epsilon},
 \label{eq:AC-invariant-exponential-tail}
\end{align}
for $0<\epsilon\leq1$ and every selected Krylov-Bogoliubov
limit. The radii are independent of the selection.
These estimates do not require
Assumption~\ref{assumption 3}.
\end{theorem}

\begin{proof}
Lemma~\ref{lem:AC-linear-uniform-controlled}, with $x=0$,
verifies \eqref{eq:terminal-controlled-bound} with
$\epsilon_*=1$. Proposition~\ref{prop:terminal-criterion}
gives both conclusions uniformly over the selected limits.
\end{proof}

\subsection{Verification of the stationary LDP}
\label{sec:LDP Dynamics}
\label{sec:AC-energy}
\label{subsec:AC-abstract-verification}

We use three results from \cite{bai2025allencahn}:
the small energy estimate of Proposition~3.1,
convergence to $u_\ast$ in Lemma~3.3, and the connection
estimate of Proposition~3.4. Their proofs do not use
the additional restriction $\rho<1$, which enters the
stationary norm estimate in Proposition~4.2 of that paper.
In particular,
\begin{equation}\label{Energy lower bound}
 g_0^2U_L(y)\leq2\mathbf E_L^*(y),
 \qquad y\in E.
\end{equation}
We verify the abstract assumptions explicitly below.

\begin{proposition}\label{prop:AC-stationary-input}
\label{Prop:characterization of U_L}
Under the standing assumptions on $g$, including
Assumption~\ref{assumption 3}, the shifted skeleton satisfies
Assumptions~\ref{ass:stationary-input} and
\ref{ass:finite-path-approximation}, with equilibrium
$u_\ast=m_L-\psi$.
Moreover, $U_L=\widehat V$, the sublevels of $U_L$
are compact in $E$, and the endpoint estimate
\eqref{eq:stationary-endpoint-approximation} holds.
\end{proposition}

\begin{proof}
In the deterministic arguments cited below, we use
\eqref{ineqn: skeleton} in place of the controlled bound
established under sublinear growth in
\cite{bai2025allencahn}.
For bounded initial values and total control energy,
this estimate places the controlled trajectories in
a common bounded subset of $E$, independently of their
duration. The linear growth of $g$ therefore bounds
$G$ uniformly along these trajectories.
Together with the retained Lipschitz and nondegeneracy
assumptions, this allows the cited deterministic arguments
to apply under the present assumptions. We use the deterministic results of
\cite{bai2025allencahn}, whose proofs do not require
the additional strictly sublinear growth condition.

\emph{Finite time compactness and continuity.}
Section~2 of \cite{bai2025allencahn} gives the compactness of the action sublevels for each fixed initial value
and duration. Its estimate (2.8) gives
\begin{equation}\label{z_phi^x-z_varphi^y}
 \|z_f^x-z_h^y\|_{C([0,T];E)}
 \leq C_{T,R,N}
 \bigl(\|x-y\|_E+\|f-h\|_{L^2(0,T;H)}\bigr)
\end{equation}
for $\|x\|_E,\|y\|_E\leq R$ and
$f,h\in\mathcal S_T^N$.
This verifies Assumption~\ref{ass:stationary-input}(i).

\emph{Bounds and positive time compactness.}
Equation~\eqref{ineqn: skeleton} gives the bound
independent of duration required in
Assumption~\ref{ass:stationary-input}(ii).

Fix $t_0,R,N>0$.
Using \eqref{ineqn: skeleton}, the energy estimates in
\cite[proof of Lemma~3.5]{bai2025allencahn},
applied on $[0,t_0]$, give
\[
 \sup_{\|x\|_E\leq R}
 \sup_{f\in\mathcal S_{t_0}^N}
       \mathbf E_L^*(z_f^x(t_0))
 \leq C_{R,N,t_0}<\infty.
\]
That argument selects an intermediate time from the
integral $H_0^1$ estimate and then propagates the energy
bound to $t_0$.
The constants depend only on $t_0$, the common $E$-bound
and the control energy bound.

Since the sublevels of $\mathbf E_L^*$ are compact
in $E$, this proves
\eqref{eq:det-positive-time-compactness}
and completes the verification of
Assumption~\ref{ass:stationary-input}(ii).

\emph{Attraction.}
Lemma~3.3 of \cite{bai2025allencahn} gives convergence
to $u_\ast$ in $H_0^1$ and convergence of the excess
energy to zero for every initial value in $E$.
For each $t>0$, heat semigroup smoothing and
\eqref{z_phi^x-z_varphi^y} give continuity of
$x\mapsto z_0^x(t)$ from $E$ to $H_0^1$.
Hence, on any compact $K\subset E$, the functions
\[
 x\longmapsto\mathbf E_L^*(z_0^x(j)),
 \qquad j=1,2,\ldots,
\]
are continuous, decrease in $j$, and converge
pointwise to zero.
Dini's theorem gives uniform convergence on $K$.
Energy monotonicity extends this to all real times,
and the small energy estimate in
\cite[Proposition~3.1]{bai2025allencahn} yields
\[
 \lim_{t\to\infty}\sup_{x\in K}
       \|z_0^x(t)-u_\ast\|_E=0.
\]
This verifies Assumption~\ref{ass:stationary-input}(iii).
The equation is autonomous, and its controlled
solutions can be restarted and concatenated.

\emph{Finite path approximation and the quasipotential.}
The finite path construction in
\cite[proof of Proposition~3.8]{bai2025allencahn}
applies under the present Lipschitz and nondegeneracy
assumptions.
For every admissible infinite past path $z$ with
finite action and every $\eta>0$, it provides
$S_\eta>0$ and
$\bar v_\eta\in C([-S_\eta,0];E)$ such that
\[
 \bar v_\eta(-S_\eta)=u_\ast,\qquad
 \bar v_\eta(0)=z(0),\qquad
 I_{-S_\eta,0}(\bar v_\eta)
 \leq I_{-\infty}(z)+\eta.
\]
The duration may depend on $z$ and $\eta$.
This verifies
Assumption~\ref{ass:finite-path-approximation}.

Proposition~\ref{prop:deterministic-compactness}
now gives $U_L=\widehat V$, compactness of the
sublevels of $U_L$, and the endpoint estimate
\eqref{eq:stationary-endpoint-approximation}.
\end{proof}

\begin{theorem}\label{thm:Allen-Cahn-invariant-LDP}
Under the standing continuity and global Lipschitz assumptions
and Assumption~\ref{assumption 3}, every selected family
$\{\bar\mu_\epsilon\}_{\epsilon>0}$ satisfies an LDP in $E$
as $\epsilon\downarrow0$, with speed $\epsilon^{-1}$ and
good rate function $U_L$.
The corresponding family on $\mathcal E$ has good rate
function
\[
 \widetilde U_L(u)=U_L(u-\psi),
 \qquad u\in\mathcal E.
\]
\end{theorem}

\begin{proof}
Proposition~\ref{prop:AC-stationary-input} verifies the
deterministic hypotheses of
Theorem~\ref{thm:stationary-transfer}.
Theorems~6.2 and~6.3 of \cite{cerrai2004large} give the
finite time FWULDP in $C([0,T];E)$ uniformly over bounded
initial sets. These results require global Lipschitz
continuity of $g$ and apply without the restriction
$\rho<1$; see also \cite[Section~2]{bai2025allencahn}.
Theorem~\ref{thm:AC-linear-invariant-exp} supplies
\eqref{eq:stationary-norm-tail} for every choice of the
Krylov-Bogoliubov limits. The abstract stationary theorem
therefore gives both bounds with good rate function $U_L$.
Translation by $\psi$ is a homeomorphism from $E$ to
$\mathcal E$ and proves the final assertion.
\end{proof}

\begin{corollary}\label{cor:Allen-Cahn-concentration}
There are $c,\delta_0>0$ such that, for every
$0<\delta<\delta_0$, there is $\epsilon_\delta>0$ with
\[
 \mu_\epsilon\bigl(\{u\in\mathcal E:
                  \|u-m_L\|_\infty\geq\delta\}\bigr)
 \leq\exp\left(-\frac{c\delta^2}{\epsilon}\right),
 \qquad0<\epsilon\leq\epsilon_\delta.
\]
\end{corollary}

\begin{proof}
For a controlled path starting at $u_\ast$, the energy
chain rule and Young's inequality give
\begin{align*}
 \mathbf E_L^*(z_f^{u_\ast}(T))
 &\leq\frac12\int_0^T
              \|G(z_f^{u_\ast}(t))f(t)\|_H^2\,dt\\
 &\leq C(1+\|z_f^{u_\ast}\|_{C([0,T];E)}^2)
                         \frac12\|f\|_{L^2(0,T;H)}^2.
\end{align*}
The chain rule is justified by the regularity of the
forced heat equation: boundedness in $E$ gives
$F(z)+G(z)f\in L^2(0,T;H)$, and
$u_\ast\in H_0^1$.
Using \eqref{ineqn: skeleton} and taking the infimum over
controls and durations yields
\begin{equation}\label{Energy upper bound}
 \mathbf E_L^*(y)\leq C U_L(y)(1+U_L(y))^2.
\end{equation}
By \cite[Proposition~3.1]{bai2025allencahn}, there are
$\eta_0,C_E>0$ such that
\[
 \|y-u_\ast\|_E^2\leq C_E\mathbf E_L^*(y)
 \quad\text{whenever}\quad
 \mathbf E_L^*(y)\leq\eta_0.
\]
Choose $r_0>0$ so that
$Cr_0(1+r_0)^2\leq\eta_0$.
Then, for a constant $K>0$,
\[
 U_L(y)\leq r_0
 \quad\Longrightarrow\quad
 \|y-u_\ast\|_E^2\leq K U_L(y).
\]
If $0<\delta<\sqrt{Kr_0}$ and
$\|y-u_\ast\|_E\geq\delta$, either this inequality gives
$U_L(y)\geq\delta^2/K$, or
$U_L(y)>r_0>\delta^2/K$.
Thus the infimum of $U_L$ outside that ball is at least
$\delta^2/K$.
Applying the LDP upper bound to this closed set, taking
$c=1/(2K)$, and translating by $\psi$ proves the assertion.
\end{proof}

\section{Stochastic Burgers equations with rough additive noise}
\label{sec:burgers-application}

Let $H=L^2(0,1)$, let $A=-\Delta$ be the positive Dirichlet
Laplacian, and write $D=\partial_\xi$. We use the spectral spaces
$H^s=Dom(A^{s/2})$ for $s\geq0$, with negative exponents defined
by completion. Set $\mathcal V=H^1=H_0^1(0,1)$, equipped with
$\|v\|_{\mathcal V}=\|Dv\|_H$. Fix $\nu>0$ and
$0\leq\alpha<1/2$, and put $S(t)=e^{-\nu tA}$. We consider
\begin{equation}\label{eq:Burgers-rough}
 du_\epsilon^x
 =\left[-\nu Au_\epsilon^x+
              \frac12D((u_\epsilon^x)^2)\right]dt
       +\sqrt\epsilon\,A^{\alpha/2}\,dW,
 \qquad u_\epsilon^x(0)=x\in H,
\end{equation}
where $W$ is a cylindrical Wiener process in $H$.
The noise covariance $A^\alpha$ has infinite trace. Nevertheless,
\[
 \int_0^T\|S(t)A^{\alpha/2}\|_{\mathrm{HS}(H)}^2\,dt
 \leq\frac1{2\nu}\sum_{j\geq1}(j\pi)^{2\alpha-2}<\infty.
\]
Thus the stochastic convolution is $H$-valued. Its continuous
version and the well-posedness of \eqref{eq:Burgers-rough}
follow from the finite time results in
\cite[Section~3.1]{BAI2026111284}.

We prove controlled moment bounds and Gaussian norm estimates
for this fixed covariance, and then apply the abstract stationary
LDP theorem. No trace class approximation of the noise is used.
Constants below may depend on the fixed parameters $\nu$ and
$\alpha$; any other dependence is indicated where needed.

\subsection{A damped Ornstein-Uhlenbeck estimate}

\begin{lemma}\label{lem:Burgers-damped-OU}
Let $T>0$ and $0<\eta<(1-2\alpha)/4$, and set
$\gamma=(1-2\alpha)/4-\eta>0$. For $\lambda\geq1$, define
\[
 Z_{\epsilon,\lambda}(t)
 :=\sqrt\epsilon\int_0^t
       e^{-\lambda(t-s)}S(t-s)A^{\alpha/2}\,dW(s).
\]
There is a nonnegative random variable $\zeta_{\eta,T}$,
independent of $\epsilon$ and $\lambda$, with moments of every
finite order, such that, on one event of probability one,
\begin{equation}\label{eq:Burgers-OU-pathwise}
 \sup_{0\leq t\leq T}\|Z_{\epsilon,\lambda}(t)\|_{L^4}
 \leq\sqrt\epsilon\,\lambda^{-\gamma}\zeta_{\eta,T}
 \quad\text{for all }\epsilon\geq0,
 \quad\lambda\in\{1,2,4,\ldots\}.
\end{equation}
\end{lemma}

\begin{proof}
Write $\gamma_0=(1-2\alpha)/4$ and fix
$b\in(0,\gamma_0)$. For each $\lambda\geq1$, set
\[
 Y_{b,\lambda}(s)
 :=\int_0^s(s-r)^{-b}e^{-\lambda(s-r)}
                    S(s-r)A^{\alpha/2}\,dW(r).
\]
The Dirichlet eigenfunctions are uniformly bounded, and
\[
 \sum_{j\geq1}(j\pi)^{2\alpha}e^{-2\nu(j\pi)^2t}
                      |e_j(\xi)|^2
 \leq Ct^{-\alpha-1/2}.
\]
Consequently,
\[
 \sup_{s\leq T,\ \xi\in[0,1]}
       \mathbb E|Y_{b,\lambda}(s,\xi)|^2
 \leq C\int_0^\infty r^{-2b-\alpha-1/2}e^{-2\lambda r}\,dr
 \leq C_b\lambda^{2b-2\gamma_0}.
\]
Gaussian moment estimates and Minkowski's inequality therefore
give, for $P\geq4$,
\[
 \left(\mathbb E\int_0^T
          \|Y_{b,\lambda}(s)\|_{L^4}^P\,ds\right)^{1/P}
 \leq C_{b,P,T}\lambda^{b-\gamma_0}.
\]
The factorisation formula reads
\[
 Z_{1,\lambda}(t)
 =\frac{\sin(\pi b)}{\pi}
       \int_0^t(t-s)^{b-1}e^{-\lambda(t-s)}
                      S(t-s)Y_{b,\lambda}(s)\,ds.
\]
If $1/P<b$, H\"older's inequality for the outer kernel yields
\begin{equation}\label{eq:Burgers-OU-deterministic-lambda}
 \left(\mathbb E
       \|Z_{1,\lambda}\|_{C([0,T];L^4)}^P\right)^{1/P}
 \leq C_{b,P,T}\lambda^{-b+1/P}\lambda^{b-\gamma_0}
 =C_{b,P,T}\lambda^{-\gamma_0+1/P}.
\end{equation}
The constants are uniform for $\lambda\geq1$.

Define
\[
 \zeta_{\eta,T}
 :=\sum_{j=0}^\infty
       2^{j\gamma}\|Z_{1,2^j}\|_{C([0,T];L^4)}.
\]
For any $q\geq1$, choose
$P>\max\{4,q,1/b,2/\eta\}$. By Minkowski's inequality and
\eqref{eq:Burgers-OU-deterministic-lambda},
\[
 \|\zeta_{\eta,T}\|_{L^q(\Omega)}
 \leq C_{b,P,T}
       \sum_{j=0}^\infty2^{-j(\eta-1/P)}<\infty.
\]
The definition gives \eqref{eq:Burgers-OU-pathwise}
simultaneously for all dyadic $\lambda$, first for
$\epsilon=1$ and then for every $\epsilon\geq0$ by scaling.
\end{proof}

\subsection{Uniform controlled moments and exponential estimates}

For $f\in\mathcal P_T^N$, let $u_{\epsilon,f}^x$ solve
\[
 du_{\epsilon,f}^x
 =\left[-\nu Au_{\epsilon,f}^x
       +\frac12D((u_{\epsilon,f}^x)^2)
       +A^{\alpha/2}f\right]dt
       +\sqrt\epsilon\,A^{\alpha/2}\,dW,
 \qquad u_{\epsilon,f}^x(0)=x.
\]
These are the controlled solutions associated with the
Brownian shifts in Section~\ref{sec:abstract-exponential-estimate};
for $\epsilon>0$, the identification follows from Girsanov's
theorem and pathwise uniqueness. For $\epsilon=0$, write $z_f^x=u_{0,f}^x$.
For deterministic $f\in L^2(0,T;H)$, this is the unique mild
solution of
\begin{equation}\label{eq:Burgers-skeleton}
 \partial_t z_f^x
 =-\nu Az_f^x+\frac12D((z_f^x)^2)+A^{\alpha/2}f,
 \qquad z_f^x(0)=x.
\end{equation}
Here the noise Hilbert space is $H$, and the coefficient
$\Sigma$ in \eqref{eq:general-skeleton-mild} is the constant
operator $A^{\alpha/2}$, interpreted through the semigroup
in the mild equation.

\begin{lemma}\label{lem:Burgers-controlled}
For every $T>0$, there is $C_T>0$ such that
\begin{equation}\label{eq:Burgers-controlled-finite}
 \sup_{0\leq\epsilon\leq1}
 \sup_{f\in\mathcal P_T^N}
 \mathbb E\sup_{0\leq t\leq T}
          \|u_{\epsilon,f}^x(t)\|_H^2
 \leq C_T(1+\|x\|_H^2+N).
\end{equation}
There is also a constant $C>0$, independent of $T$, such that
\begin{equation}\label{eq:Burgers-controlled-terminal}
 \sup_{T>0}\sup_{0\leq\epsilon\leq1}
 \sup_{f\in\mathcal P_T^N}
 \mathbb E\|u_{\epsilon,f}^x(T)\|_H^2
 \leq C(1+\|x\|_H^2+N).
\end{equation}
Both constants are independent of $x$ and $N$.
\end{lemma}

\begin{proof}
Fix $\eta\in(0,(1-2\alpha)/4)$ and the corresponding
$\gamma>0$ in Lemma~\ref{lem:Burgers-damped-OU}.
For each deterministic dyadic $\lambda$, put
$Y=u_{\epsilon,f}^x-Z_{\epsilon,\lambda}$. Then
\[
 \partial_tY
 =-\nu AY+\frac12D(Y+Z_{\epsilon,\lambda})^2
       +A^{\alpha/2}f+\lambda Z_{\epsilon,\lambda},
 \qquad Y(0)=x.
\]
Write $Z=Z_{\epsilon,\lambda}$ temporarily.
The pure Burgers term cancels in the $H$ energy identity.
The remaining nonlinear terms satisfy
\[
 \left|\int_0^1ZY\,DY\,d\xi\right|
 \leq\|Z\|_{L^4}\|Y\|_{L^4}\|DY\|_H
 \leq\frac\nu8\|Y\|_{\mathcal V}^2
       +C_\nu\|Z\|_{L^4}^{8/3}\|Y\|_H^2
\]
and
\[
 \frac12\left|\int_0^1Z^2DY\,d\xi\right|
 \leq\frac\nu8\|Y\|_{\mathcal V}^2
       +C_\nu\|Z\|_{L^4}^4.
\]
Since $\alpha<1$, the control term is bounded by
\[
 |\langle f,A^{\alpha/2}Y\rangle_H|
 \leq C\|f\|_H\|Y\|_{\mathcal V}
 \leq\frac\nu8\|Y\|_{\mathcal V}^2+C_\nu\|f\|_H^2.
\]
Young's and Poincar\'e's inequalities give the same absorption
for $\lambda\langle Z,Y\rangle_H$, with remainder
$C_\nu\lambda^2\|Z\|_H^2$.
It follows that there are $\kappa,C>0$, independent of
$T$, $\lambda$, $\epsilon$, $x$ and $f$, such that
\begin{align}
 \frac d{dt}\|Y\|_H^2
 &\leq-\kappa\|Y\|_H^2
       +C\|Z_{\epsilon,\lambda}\|_{L^4}^{8/3}\|Y\|_H^2
       \nonumber\\
 &\quad+C\left(
       \|Z_{\epsilon,\lambda}\|_{L^4}^4
       +\|f\|_H^2
       +\lambda^2\|Z_{\epsilon,\lambda}\|_H^2
       \right).
 \label{eq:Burgers-Y-energy}
\end{align}
The calculation is justified by approximation and the energy
estimate for $Y$; It\^o's formula is not applied to the rough
stochastic solution.

Choose $b_0>0$ so that $Cb_0^{8/3}\leq\kappa/2$.
On the event in Lemma~\ref{lem:Burgers-damped-OU}, choose
the least dyadic $\lambda$ satisfying
\[
 \lambda\geq
 \max\left\{1,
       (\sqrt\epsilon\,\zeta_{\eta,T}/b_0)^{1/\gamma}
       \right\}.
\]
Then
\[
 \|Z_{\epsilon,\lambda}\|_{C([0,T];L^4)}\leq b_0,
 \qquad
 \lambda\leq2\max\left\{1,
       (\sqrt\epsilon\,\zeta_{\eta,T}/b_0)^{1/\gamma}
       \right\},
\]
and hence
\[
 \sup_{t\leq T}\left(
       \|Z_{\epsilon,\lambda}(t)\|_{L^4}^4
       +\lambda^2\|Z_{\epsilon,\lambda}(t)\|_H^2
       \right)
 \leq C\left[1+
       (\sqrt\epsilon\,\zeta_{\eta,T})^{2/\gamma}\right].
\]
All convolutions and energy estimates were established for
the deterministic dyadic values before making this pathwise
selection. No stochastic integral with a random damping
parameter is required.

Applying Gronwall's inequality to
\eqref{eq:Burgers-Y-energy}, and then using $u=Y+Z$, gives
\begin{equation}\label{eq:Burgers-pathwise-controlled}
 \begin{aligned}
 \|u_{\epsilon,f}^x(t)\|_H^2
 &\leq2e^{-\kappa t/2}\|x\|_H^2
 \\
 &\quad+C\left[1+\int_0^t\|f(s)\|_H^2\,ds
       +(\sqrt\epsilon\,\zeta_{\eta,T})^{2/\gamma}
       \right],
 \qquad0\leq t\leq T.
 \end{aligned}
\end{equation}
The moment bounds for $\zeta_{\eta,T}$ prove
\eqref{eq:Burgers-controlled-finite}, including
$\epsilon=0$.

To obtain a bound independent of $T$, choose $\tau>0$ so that
$\theta:=2e^{-\kappa\tau/2}<1$ and restart the convolution
on each interval $[s,s+\tau]$. The restarted version of
$\zeta_{\eta,\tau}$ has the same distribution for every
starting time $s$. Thus \eqref{eq:Burgers-pathwise-controlled}
gives, whenever $s+\tau\leq T$,
\[
 \mathbb E\|u_{\epsilon,f}^x(s+\tau)\|_H^2
 \leq\theta\,\mathbb E\|u_{\epsilon,f}^x(s)\|_H^2
       +C_\tau(1+N).
\]
This uses the almost sure total control-energy bound and
does not require independence of the control and the Brownian
increments. Iteration yields
\[
 \mathbb E\|u_{\epsilon,f}^x(n\tau)\|_H^2
 \leq\theta^n\|x\|_H^2
       +\frac{C_\tau(1+N)}{1-\theta},
 \qquad n\tau\leq T.
\]
For $T=n\tau+r$, $0\leq r<\tau$, apply the same pathwise
estimate on the remaining interval, using the convolution
bound on an interval of length $\tau$. This gives
\[
 \mathbb E\|u_{\epsilon,f}^x(T)\|_H^2
 \leq2\mathbb E\|u_{\epsilon,f}^x(n\tau)\|_H^2
       +C_\tau(1+N),
\]
which proves \eqref{eq:Burgers-controlled-terminal}.
\end{proof}

\begin{theorem}[Gaussian norm tails]
\label{thm:Burgers-gaussian-tail}
For every $T>0$, there is $C_T>0$ such that
\begin{equation}\label{eq:Burgers-finite-gaussian}
 \mathbb P\left(
       \sup_{0\leq t\leq T}\|u_\epsilon^x(t)\|_H
       \geq C_T(1+\|x\|_H+r)\right)
 \leq e^{-r^2/\epsilon}.
\end{equation}
There is also a constant $C>0$, independent of $T$, such that
\begin{equation}\label{eq:Burgers-terminal-gaussian}
 \sup_{T>0}\mathbb P\left(
       \|u_\epsilon^x(T)\|_H
       \geq C(1+\|x\|_H+r)\right)
 \leq e^{-r^2/\epsilon}.
\end{equation}
Both estimates hold for every $x\in H$, $r>0$ and
$0<\epsilon\leq1$, with constants independent of $x$, $r$
and $\epsilon$.
\end{theorem}

\begin{proof}
Apply Theorem~\ref{thm:abstract-exponential-estimate} first
to the solution paths in $C([0,T];H)$.
The deterministic instance of
\eqref{eq:Burgers-controlled-finite} and the triangle
inequality give its controlled bounds with
\[
 C_N+\bar C_N\leq C_T(1+\|x\|_H+\sqrt N).
\]
The quantitative conclusion proves
\eqref{eq:Burgers-finite-gaussian}.
For \eqref{eq:Burgers-terminal-gaussian}, use the terminal
value maps with target $H$.
Estimate~\eqref{eq:Burgers-controlled-terminal} gives the
same square-root dependence on $N$, with constants uniform
in the terminal time.
\end{proof}

\begin{proposition}\label{prop:Burgers-invariant-norm-tails}
For every $0<\epsilon\leq1$, equation
\eqref{eq:Burgers-rough} has an invariant measure obtained
as a weak limit of Krylov-Bogoliubov averages starting at
$0$. For every $r>0$, there is $R_r>0$, independent of
$\epsilon$ and of the selected weak limit, such that
\begin{equation}\label{eq:Burgers-stationary-norm-tail}
 \mu_\epsilon^{\rm B}
       \bigl(\{x\in H:\|x\|_H>R_r\}\bigr)
 \leq e^{-r/\epsilon},
 \qquad0<\epsilon\leq1.
\end{equation}
\end{proposition}

\begin{proof}
We establish existence using an ordinary spatial moment.
Fix $0<\sigma<1/2-\alpha$. The heat kernel estimate
\[
 \|S(t)Dh\|_{H^\sigma}
 \leq Ct^{-3/4-\sigma/2}\|h\|_{L^1},
 \qquad0<t\leq1,
\]
has an integrable time singularity since $\sigma<1/2$.
Writing
\[
 Z_\epsilon(t)
 :=\sqrt\epsilon\int_0^tS(t-s)A^{\alpha/2}\,dW(s),
\]
the mild equation gives
\begin{equation}\label{eq:Burgers-spatial-smoothing}
 \|u_\epsilon^x(1)\|_{H^\sigma}
 \leq C\|x\|_H
       +C\sup_{0\leq t\leq1}\|u_\epsilon^x(t)\|_H^2
       +\|Z_\epsilon(1)\|_{H^\sigma}.
\end{equation}
Moreover,
\[
 \mathbb E\|Z_\epsilon(1)\|_{H^\sigma}^2
 \leq\frac\epsilon{2\nu}
       \sum_{j\geq1}(j\pi)^{2(\sigma+\alpha-1)}<\infty.
\]
Together with \eqref{eq:Burgers-controlled-finite}, this
implies
\[
 \mathbb E\|u_\epsilon^x(1)\|_{H^\sigma}
 \leq C(1+\|x\|_H^2),
 \qquad0<\epsilon\leq1.
\]
Restarting at time $t-1$ and using
\eqref{eq:Burgers-controlled-terminal} therefore yields
\[
 \sup_{t\geq1}\sup_{0<\epsilon\leq1}
       \mathbb E\|u_\epsilon^0(t)\|_{H^\sigma}<\infty.
\]
Since bounded subsets of $H^\sigma$ have compact closure
in $H$, the laws for $t\geq1$ are tight. In the averages
\[
 \mu_{\epsilon,T}
 :=\frac1T\int_0^T\mathcal L(u_\epsilon^0(t))\,dt,
\]
the initial interval $[0,1]$ contributes mass at most $1/T$.
Thus the averages are tight as $T\to\infty$.
Continuous dependence on the initial value, as in
\cite[Section~3.1]{BAI2026111284}, gives the Feller
property. The Krylov-Bogoliubov construction then provides
an invariant measure along each convergent subsequence.

Finally, \eqref{eq:Burgers-controlled-terminal} verifies
\eqref{eq:terminal-controlled-bound} with $x_0=0$ and
$\epsilon_*=1$. Proposition~\ref{prop:terminal-criterion}
therefore gives \eqref{eq:Burgers-stationary-norm-tail}
for every selected weak limit. The spatial moment is used
for existence; the stationary LDP uses the exponential
estimate in the $H$ norm.
\end{proof}

\subsection{Verification of the stationary LDP}

For $x\in H$ and $z\in C([0,T];H)$, define
\begin{equation}\label{eq:Burgers-path-action}
 I_T^x(z)
 :=\frac12\inf\left\{
       \int_0^T\|f(t)\|_H^2\,dt:
       f\in L^2(0,T;H),\ z=z_f^x\right\},
 \qquad\inf\varnothing=+\infty.
\end{equation}
In particular, $I_T^x(z)=+\infty$ if $z(0)\ne x$.
The quasipotential based at the equilibrium $0$ is
\begin{equation}\label{eq:Burgers-quasipotential}
 \begin{aligned}
 U_\alpha(y)
 &:=\inf_{\substack{T>0,\ z\in C([0,T];H)\\
                    z(0)=0,\ z(T)=y}}I_T^0(z)
 \\
 &=\frac12\inf\left\{
       \int_0^T\|f(t)\|_H^2\,dt:
       T>0,\ f\in L^2(0,T;H),\ z_f^0(T)=y\right\},
 \qquad y\in H.
 \end{aligned}
\end{equation}
It takes values in $[0,+\infty]$ and satisfies
$U_\alpha(0)=0$. Its dependence on the fixed viscosity
$\nu$ is suppressed in the notation.

\begin{theorem}\label{thm:Burgers-stationary-LDP}
Let $\nu>0$ and $0\leq\alpha<1/2$. Every family
$\{\mu_\epsilon^{\rm B}\}_{0<\epsilon\leq1}$ selected as in
Proposition~\ref{prop:Burgers-invariant-norm-tails}
satisfies an LDP in $H$ as $\epsilon\downarrow0$, with
speed $1/\epsilon$ and good rate function $U_\alpha$.
\end{theorem}

\begin{proof}
We verify the hypotheses of
Theorem~\ref{thm:stationary-transfer} with state space
$H$, noise Hilbert space $H$, and equilibrium $0$.

\emph{The finite time uniform LDP and action compactness.}
For the fixed amplitude $A^{\alpha/2}$, the noise
eigenvalues in \cite[Assumption~1]{BAI2026111284}
are $\sigma_j=(j\pi)^\alpha$. They satisfy
\[
 \sum_{j\geq1}\frac{\sigma_j^2}{j^2}<\infty,
 \qquad
 \epsilon\sum_{j\geq1}\frac{\sigma_j^2}{j^2}\longrightarrow0.
\]
The relative noise operator is the identity, and the
remaining spectral bound holds with exponent
$2\alpha<1<3/2$. Hence
\cite[Theorem~3.8]{BAI2026111284} gives the
Freidlin-Wentzell uniform LDP in $C([0,T];H)$,
uniformly over bounded initial sets, with good rate
function \eqref{eq:Burgers-path-action}.
The cited proof applies to $S(t)=e^{-\nu tA}$ for any
fixed $\nu>0$, with constants depending on $\nu$.
Thus \eqref{eq:stationary-finite-lower} and
\eqref{eq:stationary-finite-upper} hold, and each set
$K_T^x(r)=\{z:I_T^x(z)\leq r\}$ is compact.

\emph{Energy bounds, continuity and attraction.}
The skeleton energy identity, using
$\langle D(z^2),z\rangle=0$, gives
\[
 \frac12\frac d{dt}\|z\|_H^2
       +\nu\|z\|_{\mathcal V}^2
 =\langle f,A^{\alpha/2}z\rangle_H
 \leq\frac\nu2\|z\|_{\mathcal V}^2+C\|f\|_H^2.
\]
The identity and the following estimates are justified
by approximation. In particular,
\begin{equation}\label{eq:Burgers-skeleton-energy}
 \frac d{dt}\|z(t)\|_H^2
       +\nu\|z(t)\|_{\mathcal V}^2
 \leq C\|f(t)\|_H^2,
\end{equation}
and
\begin{equation}\label{eq:Burgers-global-skeleton-bound}
 \sup_{0\leq t\leq T}\|z_f^x(t)\|_H^2
       +\int_0^T\|z_f^x(t)\|_{\mathcal V}^2\,dt
 \leq C\left(
       \|x\|_H^2+\|f\|_{L^2(0,T;H)}^2\right),
\end{equation}
where $C$ is independent of $T$.
This proves \eqref{eq:det-uniform-bound}.
For zero control, the exact energy identity and
Poincar\'e's inequality give
\[
 \|z_0^x(t)\|_H\leq e^{-\nu\pi^2t}\|x\|_H.
\]
Thus attraction to $0$ is uniform over bounded initial sets.

For $w=z_f^x-z_h^y$, the difference equation and
$\mathcal V\hookrightarrow L^\infty(0,1)$ yield
\[
 \frac d{dt}\|w\|_H^2+\nu\|w\|_{\mathcal V}^2
 \leq C\left(
       \|z_f^x\|_{\mathcal V}^2
       +\|z_h^y\|_{\mathcal V}^2\right)\|w\|_H^2
       +C\|f-h\|_H^2.
\]
By Gronwall's inequality and
\eqref{eq:Burgers-global-skeleton-bound},
\begin{equation}\label{eq:Burgers-skeleton-continuity}
 \|z_f^x-z_h^y\|_{C([0,T];H)}
 \leq C_{R,N}\left(
       \|x-y\|_H+\|f-h\|_{L^2(0,T;H)}\right)
\end{equation}
when $\|x\|_H,\|y\|_H\leq R$ and
$f,h\in\mathcal S_T^N$.
This verifies the continuity requirement in
Assumption~\ref{ass:stationary-input}(i).

\emph{Compactness at a positive time.}
Fix $t_0>0$ and choose $0<\sigma<1/2-\alpha$.
Let $M=M_{R,N}$ be the uniform bound for
$\|z_f^x\|_{C([0,t_0];H)}$ supplied by
\eqref{eq:Burgers-global-skeleton-bound}.
The same heat kernel estimate used in
\eqref{eq:Burgers-spatial-smoothing} bounds the nonlinear
convolution at $t_0$ by $C_{t_0}M^2$.
For the control term,
\[
 \left\|\int_0^{t_0}S(t_0-s)A^{\alpha/2}f(s)\,ds
       \right\|_{H^\sigma}
 \leq C\left(\int_0^{t_0}t^{-\sigma-\alpha}\,dt\right)^{1/2}
          \|f\|_{L^2(0,t_0;H)}
 \leq C_{t_0}\sqrt N.
\]
The initial term is bounded by $C_{t_0}R$, so
\[
 \sup_{\|x\|_H\leq R,\ f\in\mathcal S_{t_0}^N}
       \|z_f^x(t_0)\|_{H^\sigma}
 \leq C_{t_0}(R+M^2+\sqrt N).
\]
The compact embedding $H^\sigma\Subset H$ proves
\eqref{eq:det-positive-time-compactness}.
All parts of Assumption~\ref{ass:stationary-input}
are therefore verified.

\emph{A connection with small action.}
For $y\in\mathcal V=H_0^1(0,1)$, define
\[
 w_y(t):=tS(1-t)y,\qquad0\leq t\leq1.
\]
Then $w_y(0)=0$ and $w_y(1)=y$. We first estimate the
regularity of this path. For $t>0$,
\[
 \frac{d}{dt}\|S(t)y\|_{\mathcal V}^2
 =-2\nu\|AS(t)y\|_H^2.
\]
Integrating this identity and using continuity in
$\mathcal V$ at $t=0$ gives
\[
 \int_0^1\|AS(t)y\|_H^2\,dt
 =
 \frac{\|y\|_{\mathcal V}^2
       -\|S(1)y\|_{\mathcal V}^2}{2\nu}
 \leq\frac1{2\nu}\|y\|_{\mathcal V}^2.
\]
Consequently,
\begin{align*}
 \int_0^1\|Aw_y(t)\|_H^2\,dt
 &=\int_0^1t^2\|AS(1-t)y\|_H^2\,dt\\
 &\leq\int_0^1\|AS(1-t)y\|_H^2\,dt\\
 &=\int_0^1\|AS(t)y\|_H^2\,dt
 \leq\frac1{2\nu}\|y\|_{\mathcal V}^2.
\end{align*}
Semigroup contractivity in $\mathcal V$ also gives
\[
 \sup_{0\leq t\leq1}\|w_y(t)\|_{\mathcal V}
 \leq\|y\|_{\mathcal V}.
\]
For $0<t<1$,
\[
 w_y'(t)=S(1-t)y+\nu Aw_y(t).
\]
Thus, by semigroup contractivity in $H$ and Poincar\'e's
inequality,
\begin{align*}
 \|w_y'\|_{L^2(0,1;H)}
 &\leq
 \|S(1-\cdot)y\|_{L^2(0,1;H)}
 +\nu\|Aw_y\|_{L^2(0,1;H)}\\
 &\leq
 \|y\|_H+\sqrt{\frac{\nu}{2}}\|y\|_{\mathcal V}
 \leq C_\nu\|y\|_{\mathcal V}.
\end{align*}
In particular,
\[
 w_y\in W^{1,2}(0,1;H)\cap L^2(0,1;D(A))
             \cap C([0,1];\mathcal V).
\]

The embedding $\mathcal V\hookrightarrow L^\infty(0,1)$
yields
\begin{align*}
 \|D(w_y(t)^2)\|_H
 &\leq2\|w_y(t)\|_{L^\infty}\|Dw_y(t)\|_H\\
 &\leq C\|w_y(t)\|_{\mathcal V}^2.
\end{align*}
Hence
\[
 \|D(w_y^2)\|_{L^2(0,1;H)}
 \leq C\|y\|_{\mathcal V}^2.
\]
Since $A^{-\alpha/2}$ is bounded on $H$, the control
\[
 f_y(t):=
 A^{-\alpha/2}
 \left(
 w_y'(t)+\nu Aw_y(t)-\frac12D(w_y(t)^2)
 \right)
\]
belongs to $L^2(0,1;H)$ and satisfies
\[
 \|f_y\|_{L^2(0,1;H)}
 \leq C\left(\|y\|_{\mathcal V}
                    +\|y\|_{\mathcal V}^2\right).
\]
By construction, $w_y$ solves
\eqref{eq:Burgers-skeleton} with initial value $0$ and
control $f_y$. Therefore
\begin{equation}\label{eq:Burgers-small-connection}
 I_1^0(w_y)
 \leq\frac12\|f_y\|_{L^2(0,1;H)}^2
 \leq C\left(\|y\|_{\mathcal V}^2
                    +\|y\|_{\mathcal V}^4\right).
\end{equation}

We now verify the criterion in
Proposition~\ref{prop:local-connection-identification}.
Let $0<\delta\leq1$ and let $w\in C([0,1];H)$ satisfy
\[
 \|w(0)\|_H<\delta,\qquad
 I_1^{w(0)}(w)<\delta.
\]
By the definition of the action, there is
$f\in L^2(0,1;H)$ such that
\[
 w=z_f^{w(0)},\qquad
 \frac12\int_0^1\|f(t)\|_H^2\,dt\leq2\delta.
\]
Integrating \eqref{eq:Burgers-skeleton-energy} gives
\[
 \int_0^1\|w(t)\|_{\mathcal V}^2\,dt
 \leq C\left(
       \|w(0)\|_H^2+\int_0^1\|f(t)\|_H^2\,dt
       \right)
 \leq C(\delta^2+\delta).
\]
There is therefore a time $s\in(1/2,1)$ such that
$w(s)\in\mathcal V$ and
\[
 \|w(s)\|_{\mathcal V}^2
 \leq2\int_{1/2}^1\|w(t)\|_{\mathcal V}^2\,dt
 \leq C_1\delta.
\]
Set $y=w(s)$. The path $w_y$ connects $0$ to this exact
state in time one, and \eqref{eq:Burgers-small-connection}
implies
\[
 I_1^0(w_y)
 \leq C\left(C_1\delta+C_1^2\delta^2\right)
 \leq C_2\delta.
\]
The constants are independent of $w$ and $\delta$.
For every $\eta>0$, choosing
\[
 0<\delta\leq\min\{1,\eta/C_2\}
\]
therefore gives a connection from $0$ to $w(s)$ with
action at most $\eta$.
This verifies the local connection criterion of
Proposition~\ref{prop:local-connection-identification},
which establishes
Assumption~\ref{ass:finite-path-approximation}.

\emph{Application of the stationary theorem.}
Proposition~\ref{prop:deterministic-compactness} now
identifies $U_\alpha$ with the infinite past quasipotential
and gives compactness of its sublevels.
Proposition~\ref{prop:Burgers-invariant-norm-tails}
provides \eqref{eq:stationary-norm-tail}.
Together with the finite time uniform LDP,
Theorem~\ref{thm:stationary-transfer} gives the stated
LDP with good rate function $U_\alpha$.
\end{proof}

\begin{remark}
The stochastic input to the stationary theorem is the
uniform controlled terminal estimate
\eqref{eq:Burgers-controlled-terminal}, which gives an
exponential estimate outside bounded subsets of $H$
through the variational criterion.
Compactness and the finite path approximation are verified
for the deterministic skeleton. Thus the proof does not
require a stochastic spatial exponential estimate or
It\^o's formula for the full rough noise solution.
\end{remark}

\section*{Acknowledgements}
\addcontentsline{toc}{section}{Acknowledgements}
We acknowledge financial support from EPSRC grant EP/S005293/2
and Royal Society Newton Fund grant NIF\textbackslash R1\textbackslash221003.

\phantomsection
\addcontentsline{toc}{section}{References}
\begingroup
\footnotesize
\bibliographystyle{siam}
\bibliography{Ref}
\endgroup
\end{document}